\documentclass[9pt,reqno,final]{amsart}
\usepackage{amsfonts}
\usepackage{amsmath,version}
\usepackage{amssymb,graphicx,fancybox,mathrsfs,stmaryrd,relsize}
\usepackage[notcite,notref]{showkeys}
\usepackage{ushort} 
\usepackage{url,hyperref}
\usepackage{subfigure}
\usepackage{color}
\usepackage{cases}
\usepackage{subeqnarray} 
\usepackage{cases}
\usepackage{algorithm}
\usepackage{algorithmicx}
\usepackage{algpseudocode}

\usepackage{multirow}
\usepackage{makecell}

\catcode`\@=11 \theoremstyle{plain}
\@addtoreset{equation}{section}
\renewcommand{\theequation}{\arabic{section}.\arabic{equation}}

\newtheorem{prop}{\bf Proposition}[section]
\newtheorem{rem}{\bf Remark}[section]
\newtheorem{thm}{\bf Theorem}[section]
\newtheorem{lm}{\bf Lemma}[section]

\newtheorem{exm}{\bf Example} 
\newenvironment{example}{\begin{exm}} {\end{exm}} 

\newcommand{\bs}[1]{\boldsymbol{#1}}

\graphicspath{{../figs/}}

\begin{document}

{\title[Spectral algorithms for Maxwell's equations] {Highly efficient Gauss's law-preserving spectral algorithms for Maxwell's double-curl source and eigenvalue problems based on eigen-decomposition}
\author[
S. Lin, \,    H. Li \,  $\&$ \,  Z. Yang
]{
	Sen Lin${}^{1}$,   \,\,   Huiyuan Li${}^{2}$ \,\, and \,\,  Zhiguo Yang${}^{*,3}$
	}
\thanks{
\noindent ${}^{1}$ School of Mathematical Sciences, Shanghai Jiao Tong University, Shanghai 200240, China. 
Email: sjtu\_Einsteinlin7@sjtu.edu.cn (S. Lin).
\\
\indent \,\,\, ${}^{2}$ State Key Laboratory of Computer Science/Laboratory of Parallel Computing, Institute of
Software, Chinese Academy of Sciences, Beijing 100190, China. 
Email: huiyuan@iscas.ac.cn (H. Li).
\\
\indent \,\,\, ${}^{3}$ Corresponding author. School of Mathematical Sciences, MOE-LSC and CMA-Shanghai, Shanghai Jiao Tong University, Shanghai 200240, China. 
Email: yangzhiguo@sjtu.edu.cn (Z. Yang)
}
}


\begin{abstract} 
	In this paper, we present Gauss's law-preserving spectral methods and their efficient solution algorithms for curl-curl source and eigenvalue problems 
	in two and three dimensions arising from Maxwell's equations. 
	Arbitrary order $\bs{H}({\rm curl})$-conforming spectral basis functions in two and three dimensions are firstly proposed using compact combination of Legendre polynomials. 
	A mixed formulation involving a Lagrange multiplier is then adopted to preserve the Gauss's law in the weak sense. 
	To overcome the bottleneck of computational efficiency caused by the saddle-point nature of the mixed scheme, we present highly efficient solution algorithms based on 
	reordering and decoupling of the resultant linear algebraic system and numerical eigen-decomposition of one dimensional mass matrix. 
	The proposed solution algorithms are direct methods requiring only several matrix-matrix or matrix-tensor products of $N$-by-$N$ matrices, where $N$ is the highest polynomial order in each direction. 
	Compared with other direct methods, the computational complexities are reduced from $\mathcal{O}(N^{6})$ and $\mathcal{O}(N^{9})$ 
	to $\mathcal{O}(N^{3})$ and $\mathcal{O}(N^{4})$ with small and constant pre-factors for 2D and 3D cases, respectively, 
	and can further be accelerated to $\mathcal{O}(N^{\log_{2}{7}})$ and $\mathcal{O}(N^{1+\log_{2}{7}})$  (${\log_{2}{7}}\approx{2.807}$), 
	when boosted with the Strassen's matrix multiplication algorithm. 
	Moreover, these algorithms strictly obey the Helmholtz-Hodge decomposition, thus totally eliminate the spurious eigen-modes of non-physical zero eigenvalues. 
	Extensions of the proposed methods and algorithms to problems in complex geometries with variable coefficients and inhomogeneous boundary conditions are discussed to deal with more general situations.  
	Ample numerical examples for solving Maxwell's source and eigenvalue problems are presented to demonstrate the accuracy and efficiency of the proposed methods.
\end{abstract}

\maketitle

\section{Introduction}
\label{Sec1}
This work concerns accurate and efficient spectral methods for two typical model problems arising from time-harmonic Maxwell's equations in two and three dimensions, 
i.e., the double-curl source problem:
\begin{equation}
\label{eq:MaxwSouProb}
\nabla\times\nabla\times\bs{u} + \kappa\,\bs{u}=\bs{f}, \quad \nabla\cdot\bs{u} = \rho \;\; \text{in} \;\; \Omega, \quad 
\bs{n}\times\bs{u} = \bs{0} \;\; \text{on} \;\; \partial\Omega,
\end{equation}
and the eigenvalue problem: find $\lambda\in\mathbb{R}$ and $\bs{u}\neq\bs{0}$ such that 
\begin{equation}
\label{eq:MaxwEigProb}
\nabla\times\nabla\times\bs{u} =\lambda\bs{u}, \quad \nabla\cdot\bs{u} = 0 \;\;  \text{in} \;\; \Omega, \quad 
\bs{n}\times\bs{u} = \bs{0} \;\; \text{on} \;\; \partial\Omega.
\end{equation}
Here, $\bs{x}\in\Omega\subset\mathbb{R}^d$, $d=2,3$, $\kappa\in\mathbb{C}$, $\bs{u}(\bs{x})$ usually represents the electric field, 
$\bs{f}(\bs{x})$ is the current density and $\rho(\bs{x})$ is the charge density satisfying $\rho=\nabla\cdot\bs{f}/{\kappa}$, 
$\bs{n}$ is the unit outward normal vector to the boundary $\partial\Omega$ and this boundary condition corresponds to the perfect electric conducting boundary condition.  
The Gauss's law $\nabla\cdot\bs{u}=\rho$ (resp. $\nabla\cdot\bs{u}=0$) is implicitly implied by the double-curl equation and serves as intrinsic constraint for the model problems.

Accurate numerical simulation of these double-curl problems with Gauss's law constraint plays a crucial role in understanding physical mechanisms 
and solving practical engineering problems related to electromagnetism, such as astrophysics, inertial confinement fusion, telecommunication, 
electromagnetic materials and devises, semi-conductor manufacturing, quantum chromodynamics, quantum field theory and more  
\cite{Buchholz-1986,ChenFF1984,Davidson2001,JoanJohnsonWinn2008,Kijowski-Rudolph-2002,MarkowichRinghofer2012,
	Mund-Rehren-Schroer-2020,Orfanidis2002,Sanders-Dappiaggi-Hack-2014, Zhang-Wang-Yang-2020}. 
Nevertheless, the double-curl problems \eqref{eq:MaxwSouProb} and \eqref{eq:MaxwEigProb} are prohibitively difficult to solve numerically, due to the involvement of Gauss's law constraint.
On the one hand, it poses significant challenge to construct suitable approximation spaces that 
guarantee the well-posedness of discrete problems and maintain Gauss's law at the discrete level (either strongly or weakly). 
Failure to satisfy this constraint may result in severe numerical inaccuracies or the occurrence of spurious solutions. 
On the other hand, the saddle-point nature induced by the constraint makes them more expensive to solve and more challenging to construct efficient solution algorithms, 
compared with their self-adjoint coercive counterparts.

These challenges and difficulties have attracted increasing interest from the community in the past six decades,
and a great many studies have been devoted to solving the Maxwell's double-curl problem and associated eigenvalue problem numerically. 
We restrict our review of literature to methods based on weak or variational formulation, 
which leaves out a vast of algorithms based on staggered-grid finite difference Yee scheme (see e.g. \cite{Taflove-1980,Yee-1966}) 
or divergence-correction technique (see e.g. \cite{Balsara-Kim-2004,Dedner-Kemm-Kroner-2002,Munz-Omnes-Schneider-2000}). 
The endeavor traces back to N\'{e}d\'{e}lec, who firstly proposed two families of edge elements to discretize the Maxwell's equation (see \cite{Nedelec-1980,Nedelec-1986}, 
and comprehensive discussions in \cite{Monk2003,SunZhou2017}). 
Kikuchi in \cite{Kikuchi-1987} proposed a mixed formulation by introducing a Lagrange multiplier to maintain the Gauss's law in a weak sense. 
In \cite{Bramble-Kolev-Pasciak-2005}, Bramble et. al. approximated the Maxwell’s eigenvalue problem based on two very weak div-curl formulation 
with complementary boundary conditions. 
In \cite{Ciarlet-Wu-Zou-2014}, Ciarlet et. al. proposed and investigated some novel edge element approximations to Maxwell's equations, 
which ensure the optimal strong convergence of the Gauss's laws in some appropriate norm and avoid the solution of saddle-point discrete systems. 
This method was further extend in \cite{Yousept-Zou-2017} to solve the optimal control of the stationary Maxwell system with a non-vanishing charge density. 
In \cite{Liu-Tobon-Tang-Liu-2015}, the authors proposed a higher order nodal spectral element method for 2D Maxwell's equations 
based on Gauss-Lobatto-Legendre interpolating basis functions and Kikuchi's mixed formulation. 
Duan et. al. in \cite{Duan-Ma-Zou-2021} proposed a modification of Kikuchi's formulation to enforce the Gauss's law locally at the discrete level. 
As for efficient solution algorithms, several efficient iterative methods for the solution of Maxwell's equations discretized with edge elements 
have been proposed in \cite{Hiptmair-Xu-2007,Hu-Zou-2004,Liang-Xu-2022,Liang-Xu-2023,Pazner-Kolev-Dohrmann-2023}.

To the best of the author's knowledge, there is currently no spectral methods that simultaneously gaurantee the Gauss's law at the discrete level 
and offer fast solution algorithms for Maxwell's equations in both two and three dimensions.  
In order to fill this gap in the existing literature, we propose a family of $\bs{H}({\rm curl})$-conforming vectorial spectral basis functions of arbitrary order 
using compact combination of Legendre polynomials for both two and three dimensional cases. 
By combining the proposed spectral approximation and the mixed formulation in \cite{Kikuchi-1987}, 
we arrive at a Gauss's law preserving spectral method for these double-curl problems.  
The novelty of this work primarily resides in the algorithmic aspects, as listed in the following three points:
\begin{itemize}
\item 
A highly efficient solution algorithm is proposed for the double-curl source problem, which is based on a sophisticated decoupling procedure 
via reordering of the degrees of freedom as well as a matrix diagonalisation technique. 
The proposed algorithm has a significant merit in that it requires only several matrix-matrix or matrix-tensor products of 1D matrices 
with complexities of $\mathcal{O}(N^{\log_{2}7})$ and $\mathcal{O}(N^{1+\log_{2}7})$ with small and constant pre-factors for 2D and 3D cases, respectively, 
thus is highly efficient and well-suited for modern computing architecture. 
Numerical results illustrate the superior accuracy and efficiency of the proposed spectral method and corresponding algorithm, even for indefinite systems and highly oscillatory solutions.

\item 
The proposed method is further extended to the Maxwell's eigenvalue problem, which shares the same merits with that of the source problem 
such that the eigenvalue and eigenvector can be obtained by a semi-analytical manner. 

\item 
Moreover, the solution algorithm adheres to the discrete Helmholtz-Hodge decomposition exactly, which totally eliminates spurious eigen-modes of non-physical zero eigenvalues.  
As a notable outcome, we conclude through numerical experiments that the percentage of reliable eigenvalues among the numerical eigenvalues obtained 
by the proposed Gauss's law preserving spectral method for 2D and 3D Maxwell's eigenvalue problems are asymptotically $(2/{\pi})^2$ and $(2/{\pi})^3$, respectively. 
These findings serve as compelling evidence that this spectral eigen-solver is highly advantageous for accurate and large-scale computations. 
\end{itemize}

It is also necessary to point out that the proposed $\bs{H}({\rm curl})$-conforming spectral basis is equivalent to 
that one constructed using integrated Legendre polynomials in the realm of high-order finite element methods, up to normalization constants (\cite{Fuentes-Keith-2015,Zaglmayr2006}). 
Thus, the proposed method may also shed light on developing efficient solution algorithms for high-order finite element methods and spectral element methods.

The rest of this paper is organized as follows. 
In Section \ref{Sec2}, we propse a new $\bs{H}({\rm curl})$-conforming vectorial spectral basis functions 
and present the Gauss's law preserving spectral method for Maxwell's double-curl source problems in two dimensions. 
Highly efficient matrix-free solution algorithm based on numerical eigen-decomposition of the 1D mass matrix is then developed for the proposed spectral method.  
Following the same spirt, we present the approximate scheme and associated solution algorithm for the double-curl eigenvalue problem in 2D. 
We then extend the proposed spectral method and corresponding solution algorithm in three dimensions in Section \ref{Sec3}. 
Discussions of the proposed method to complex geometries with variable coefficients and inhomogeneous boundary conditions are presented in Section \ref{Sec4}. 
In Section \ref{Sec5}, we present several numerical results to demonstrate the accuracy and efficiency of the proposed methods and algorithms. 
We conclude with closing remarks in the last section.

\vspace{5pt} 
\noindent{\bf Notations and conventions}: 
Throughout this paper, we denote scalars, vectors and matrices by italic, bold italic and bold non-italic letters, such as $f$, $\bs{u}=(u_i)$ and $\mathbf{A}=(A_{ij})$, respectively. 
Let $\mathbb{R}$, $\mathbb{Z}$ and $\mathbb{C}$ denote the spaces of real, integer and complex numbers, respectively. 
Let $L^2(\Omega)$ be the square-integral function space with inner product $(\cdot,\cdot)_{\Omega}$ and norm $\|\cdot\|_{\Omega}$ defined on domain $\Omega$.  
The Sobolev spaces $H^m(\Omega)$ ($m\geq0$, $m\in\mathbb{Z}$) and $\bs{H}({\rm curl};\Omega)$ follow standard definitions (see e.g. \cite{Adam1975,Monk2003}) 
with norms $\|\cdot\|_{m,\Omega}$ and $\|\cdot\|_{{\rm curl},\Omega}$, respectively. 
The subscript $\Omega$ is omitted, if no ambiguity occurs. $\nabla\times$ and $\nabla\cdot$ denote the usual curl and divergence operators.

Besides, we denote $\{ \bs{e}_k \}_{k=1}^{N-1}$ and $\{\bs{\ushort{e}}_k \}_{k=1}^N$ be the canonical bases for vector spaces $\mathbb{R}^{(N-1)\times{1}}$ and $\mathbb{R}^{N\times{1}}$ 
and $\mathbf{I}_{N}\in\mathbb{R}^{N\times{N}}$ and $\mathbf{0}_{N}\in\mathbb{R}^{N\times{N}}$ be the identity and zero matrices, respectively. 
The following decomposition of identity matrix $\mathbf{I}_N\in\mathbb{R}^{N\times{N}}$ finds useful:
\begin{equation}
\label{eq:MatINE}
\mathbf{I}_N = 
\begin{pmatrix}
\bs{\ushort{e}}_1^{\intercal} \\
\mathbf{\ushort{E}} 
\end{pmatrix}
, \quad 
\mathbf{\ushort{E}} = 
\begin{pmatrix}
\bs{\ushort{e}}_2^{\intercal} \\
\vdots \\
\bs{\ushort{e}}_N^{\intercal} 
\end{pmatrix}
\in\mathbb{R}^{(N-1)\times{N}},
\quad
\mathbf{I}_{N-1} = \mathbf{\ushort{E}}\mathbf{\ushort{E}}^{\intercal}.
\end{equation} 
We further introduce the notation of  Kronecker product of matrices $\otimes$ and vectorization operation ``${\rm vec}$" to transform matrix or tensor $\mathbf{U}$ into vector $\bs{\vec{U}}$, 
i.e. $ \bs{\vec{U}} ={\rm vec}(\mathbf{U})$ following standard column-wise linear index order (see \cite[p.~27-28]{GolubCharles2013} for more details). 
Here, we denote vectors obtained from vectorization operator by $\bs{\vec{\cdot}}$ to distinguish from the usual vector notation. 
It is also necessary to define the inverse operation ``${\rm ivec}$" such that $\mathbf{U} ={\rm ivec}(\bs{\vec{U}})$ iff $ \bs{\vec{U}} ={\rm vec}(\mathbf{U})$, following the same indexing rule.

\section{Maxwell's double-curl problems in two dimensions} 
\label{Sec2}
In this section, we focus our attention on the Maxwell's double-curl source problem \eqref{eq:MaxwSouProb} and eigenvalue problem \eqref{eq:MaxwEigProb} in two dimensions. 
The starting point of our proposed spectral approximation and corresponding algorithm is the mixed weak formulation proposed by Kikuchi in \cite{Kikuchi-1987}, 
where a Lagrange multiplier $p$ is introduced to impose the Gauss's law in the weak sense.  
The weak formulation reads: find $\bs{u}\in\bs{H}_{0}({\rm curl};\Omega)$ and $p\in{H}_{0}^1(\Omega)$ such that
\begin{equation}
\label{eq:MaxwWeak}
\begin{aligned}
& ( \nabla\times\bs{u},\nabla\times\bs{v}) + \kappa ( \bs{u},\bs{v} ) + ( \nabla{p},\bs{v}) = ( \bs{f},\bs{v} ),\quad \forall\bs{v}\in\bs{H}_{0}({\rm curl};\Omega), \\
& ( \bs{u},\nabla{q} ) = - (\rho,q), \quad \forall{q} \in{H}_{0}^1(\Omega),
\end{aligned}
\end{equation} 
where $\bs{H}_{0}({\rm curl};\Omega):=\big\{\bs{v}\in\bs{H}({\rm curl};\Omega),\;\bs{n}\times\bs{v}=\bs{0}\;\text{on}\;\partial\Omega\big\}$
and ${H}_{0}^1(\Omega):=\big\{ {v}\in{H}^1(\Omega),\; v= 0 \; \text{on}\;\partial \Omega \big\}$.
Though numerical methods based on this formulation has been studied extensively under the framework of mixed finite element methods (\cite{Duan-Ma-Zou-2021,Monk2003,SunZhou2017}), 
it is less explored in the realm of spectral methods and spectral-accurate approximation and efficient solution algorithm are highly demanded.

\subsection{Spectral approximation scheme for system \texorpdfstring{\eqref{eq:MaxwWeak}}{Lg}}
\label{Sec2-1}
In what follows, let us fix $\Omega:={\Lambda}^2=(-1,1)^2$ and propose the approximation pair $\bs{H}_{N,0}({\rm curl};{\Lambda}^2) \times{H}_{N,0}^1({\Lambda}^2)$ for $(\bs{u},p)$ as follows. 
\begin{prop}
\label{Prop1-DiscSpace2d}
Denote the index sets $\mathbb{N}=\{0,1,\cdots,N-1\}$, $\mathbb{N}^{+}=\mathbb{N}\setminus\{0\}$ and functions
\begin{equation}
\label{eq:phipsi}
\phi_m(\xi) = \sqrt{\frac{2m+1}{2}}L_m(\xi),\quad 
\psi_{m+1}(\xi) = \frac{1}{\sqrt{2(2m+1)}} \big(L_{m+1}(\xi)-L_{m-1}(\xi) \big), \quad \xi\in \Lambda,
\end{equation}
where $L_m(\xi)$ is the Legendre polynomial of order $m\geq0$.  
The conforming spectral approximation spaces $\bs{H}_{N,0}({\rm curl};{\Lambda}^2)$ and ${H}_{N,0}^1({\Lambda}^2)$ of order $N$ for $\bs{H}_{0}({\rm curl};{\Lambda}^2)$ 
and ${H}_{0}^{1}({\Lambda}^2)$, 
respectively, are defined by
\begin{align}
	& \label{eq:HN0curl2d}
	\bs{H}_{N,0}({\rm curl};{\Lambda}^2) = {\rm span} \Big( \big\{ \bs\Phi_{m,n}^1(\bs{x})\big\}_{(m,n)\in(\mathbb{N}, \mathbb{N}^{+}) },
	\big\{\bs\Phi_{m,n}^2(\bs{x}) \big\}_{(m,n)\in(\mathbb{N}^{+}, \mathbb{N}) }  \Big),
	\\
	& 
	{H}_{N,0}^1({\Lambda}^2) ={\rm span} \Big( \big\{ \Psi_{m,n}(\bs{x}) \big\}_{(m,n)\in(\mathbb{N}^{+}, \mathbb{N}^{+}) } \Big),
	\label{eq:HN012d}
\end{align}
where 
\begin{equation}
	\label{eq:PhiPsi2d}
	\begin{aligned}
	& {\bs\Phi}_{m,n}^1(\bs{x}) = \big( \phi_{m}(x_1) \psi_{n+1}(x_2), 0 \big)^{\intercal},  \;\;
	\\
	& 
	{\bs\Phi}_{m,n}^2(\bs{x}) = \big( 0, \psi_{m+1}(x_1) \phi_{n}(x_2) \big)^{\intercal}, \;\;
	\\
	& 
	{\Psi}_{m,n}(\bs{x}) = \psi_{m+1}(x_1)\psi_{n+1}(x_2).
	\end{aligned}
\end{equation}
\end{prop}

Correspondingly, the approximation scheme for weak formulation \eqref{eq:MaxwWeak} in two dimensions takes the form: 
find $\bs{u}_N\in\bs{H}_{N,0}({\rm curl};{\Lambda}^2)$ and $p_N\in{H}_{N,0}^1({\Lambda}^2)$ such that
\begin{equation}
\label{eq:MaxwDis2d}
\begin{aligned}
& ( \nabla\times\bs{u}_N,\nabla\times\bs{v} ) + {\kappa} ( \bs{u}_N,\bs{v} ) + ( \nabla{p_N},\bs{v} ) = ( \bs{f},\bs{v} ), \quad \forall \bs{v}\in\bs{H}_{N,0}({\rm curl};{\Lambda}^2),
\\
& (\bs{u}_N,\nabla{q}) = -(\rho,q), \quad \forall{q}\in {H}_{N,0}^1({\Lambda}^2).
\end{aligned}
\end{equation}

\begin{lm}[see \cite{Shen-1994b}]
	\label{Lem1-MIMat} 
	Denote two symmetric matrices $\mathbf{M}=(M_{mn})\in\mathbb{R}^{(N-1)\times(N-1)}$ and $\mathbf{I}=(I_{mn})\in\mathbb{R}^{(N-1)\times(N-1)}$, 
	where the entries are given by
	\begin{equation}
	\label{eq:MIdef}
	M_{mn} = \big( \psi_{n+1}, \psi_{m+1} \big) =  \int_{-1}^{1} \psi_{n+1}(\xi)\psi_{m+1}(\xi) \,\mathrm{d}\xi,  \;\;
	I_{mn} =  \big( \phi_{n}, \phi_{m} \big) = \int_{-1}^{1} \phi_n(\xi)\phi_m(\xi) \,\mathrm{d}\xi. \quad  
	\end{equation}
	Then there hold
	\begin{align*}
	& M_{mn}=M_{nm}=
	\begin{cases}
	 \dfrac{1}{2n+1}\big( \dfrac{1}{2n-1}+\dfrac{1}{2n+3} \big),         \quad & m=n,  \\[5pt] 
	 - \dfrac{1}{\sqrt{2n+1}}\dfrac{1}{\sqrt{2n+5}}\dfrac{1}{2n+3} ,  \quad & m=n+2, \\[5pt]
	0,\quad & {\rm otherwise};
	\end{cases}
	\qquad \; I_{mn}=\delta_{mn}.
	\end{align*}
\end{lm}

Let us expand the numerical solution $\bs{u}_N(\bs{x})$ and $p_N(\bs{x})$ by
\begin{equation}
\label{eq:uNpN2d}
\bs{u}_N(\bs{x})
= \sum_{m=0}^{N-1}\sum_{n=1}^{N-1} u_{mn}^1 \bs\Phi_{m,n}^1(\bs{x}) + \sum_{m=1}^{N-1}\sum_{n=0}^{N-1} u_{mn}^2 \bs\Phi_{m,n}^2(\bs{x}) ,
\quad \; p_N(\bs{x}) = \sum_{m,n=1}^{N-1} p_{mn} {\Psi}_{m,n}(\bs{x}),
\end{equation}
and define source matrices $\mathbf{F}=(f_{mn})\in\mathbb{R}^{N\times(N-1)}$, $\mathbf{G}=(g_{mn})\in\mathbb{R}^{(N-1)\times{N}}$ and $\mathbf{R}=(r_{mn})\in\mathbb{R}^{(N-1)\times(N-1)}$ with
\begin{equation}
\label{eq:MatFGR2d}
f_{mn}=\big( \bs{f}, \bs\Phi_{m,n}^1 \big), \;\; g_{mn}=\big(\bs{f}, \bs\Phi_{m,n}^2 \big), \;\; r_{mn}=-\big( \rho,{\Psi}_{m,n} \big).
\end{equation}
Taking $\bs{v}= \bs{\Phi}_{\hat{m},\hat{n}}^i(\bs{x})$, $i=1,2$ and $q({\bs{x}})=\Psi_{\hat{m},\hat{n}}(\bs{x})$ into the numerical scheme \eqref{eq:MaxwDis2d} and denoting the unknowns
\begin{equation}
\label{eq:upVec2d}
\begin{aligned}
& \mathbf{U}=(u_{mn}^1)\in\mathbb{R}^{N\times(N-1)},\quad \mathbf{V}=(u_{mn}^2)\in\mathbb{R}^{(N-1)\times{N}}, \quad \mathbf{P}=(p_{mn})\in\mathbb{R}^{(N-1)\times(N-1)},
\\
& 
\bs{\vec{u}}={\rm vec}( \mathbf{U} ), \quad \bs{\vec{v}}={\rm vec}( \mathbf{V} ), \quad \bs{\vec{p}}={\rm vec}( \mathbf{P} ), \quad  
\bs{\vec{f}}={\rm vec}( \mathbf{F} ), \quad \bs{\vec{g}}={\rm vec}( \mathbf{G} ), \quad \bs{\vec{r}}={\rm vec}( \mathbf{R}), 
\end{aligned}
\end{equation}
one arrives at the following algebraic system
\begin{equation}
\label{eq:MaxwAlgSys2d}
\begin{pmatrix}
( \mathbf{I}_{N-1} + {\kappa}\mathbf{M} )\otimes{\mathbf{I}_N}
& - \mathbf{\ushort{E}}\otimes\mathbf{\ushort{E}}^{\intercal}
& \mathbf{M}\otimes\mathbf{\ushort{E}}^{\intercal}
\\
-\mathbf{\ushort{E}}^{\intercal}\otimes\mathbf{\ushort{E}}
& {\mathbf{I}_N}\otimes( \mathbf{I}_{N-1} + {\kappa}\mathbf{M} )
& \mathbf{\ushort{E}}^{\intercal}\otimes\mathbf{M}
\\
\mathbf{M}\otimes\mathbf{\ushort{E}}
& \mathbf{\ushort{E}}\otimes\mathbf{M}
& \mathbf{0}_{N-1}\otimes\mathbf{0}_{N-1} 
\end{pmatrix}
\begin{pmatrix}
\bs{\vec{u}} \\
\bs{\vec{v}} \\
\bs{\vec{p}}
\end{pmatrix}
=\begin{pmatrix}
\bs{\vec{f}} \\
\bs{\vec{g}} \\
\bs{\vec{r}}
\end{pmatrix},
\end{equation}
where $\mathbf{\ushort{E}}\in\mathbb{R}^{(N-1)\times{N}}$ is given by equation \eqref{eq:MatINE}.

\subsection{Numerical algorithm} 
\label{Sec2-2}
It is well-known that due to the indefiniteness and poor spectral property caused by the saddle-point nature of the mixed formulation, 
linear system \eqref{eq:MaxwAlgSys2d} poses a significant numerical obstacle for efficient solution algorithms. 
Here, we propose a novel semi-analytic approach to solve this algebraic system efficiently. 
The derivation of the algorithm can be divided into the following two steps:

\textbf{Step 1: Reordering and decoupling linear system \eqref{eq:MaxwAlgSys2d}.} 
In fact, one can reorder linear system \eqref{eq:MaxwAlgSys2d} by left multiplying with the permutation matrix 
\begin{equation}
\label{eq:MaxwPemMat2d}
\begin{pmatrix}
\mathbf{I}_{N-1}\otimes\bs{\ushort{e}}_{1}^{\intercal}  &  &   \\
\mathbf{I}_{N-1}\otimes\mathbf{\ushort{E}} &  &  \\
&   \bs{\ushort{e}}_{1}^{\intercal}\otimes\mathbf{I}_{N-1} &  \\
&   \mathbf{\ushort{E}}\otimes\mathbf{I}_{N-1}    &  \\
&  &    \mathbf{I}_{N-1}\otimes\mathbf{I}_{N-1}
\end{pmatrix},
\end{equation}
which directly leads to
\begin{equation}
\label{eq:MaxwDAlgSys2d}
\begin{pmatrix}
\mathbf{I}+{\kappa}\mathbf{M}   &  &   &  & \\
& ( \mathbf{I} + {\kappa}\mathbf{M} )\otimes\mathbf{I} &   & -\mathbf{I}\otimes\mathbf{I}  & \mathbf{M}\otimes\mathbf{I}
\\
&  & \mathbf{I} + {\kappa}\mathbf{M}     &  & 
\\
&  -\mathbf{I}\otimes\mathbf{I}  & & \mathbf{I}\otimes(\mathbf{I}+{\kappa}\mathbf{M})     & \mathbf{I}\otimes\mathbf{M} 
\\		
&  \mathbf{M}\otimes\mathbf{I}  &  & \mathbf{I}\otimes\mathbf{M}  &  \mathbf{0}\otimes\mathbf{0} 
\end{pmatrix}
\begin{pmatrix}
(\mathbf{I}\otimes\bs{\ushort{e}}_{1}^{\intercal}) \bs{\vec{u}} \\
(\mathbf{I}\otimes\mathbf{\ushort{E}}) \bs{\vec{u}}  \\
(\bs{\ushort{e}}_{1}^{\intercal}\otimes\mathbf{I}) \bs{\vec{v}}  \\
(\mathbf{\ushort{E}}\otimes\mathbf{I}) \bs{\vec{v}}  \\
(\mathbf{I}\otimes\mathbf{I}) \bs{\vec{p}} 
\end{pmatrix}
=\begin{pmatrix}
(\mathbf{I}\otimes\bs{\ushort{e}}_{1}^{\intercal}) \bs{\vec{f}} \\
(\mathbf{I}\otimes\mathbf{\ushort{E}})  \bs{\vec{f}}  \\
(\bs{\ushort{e}}_{1}^{\intercal}\otimes\mathbf{I}) \bs{\vec{g}} \\
(\mathbf{\ushort{E}}\otimes\mathbf{I}) \bs{\vec{g}} \\
(\mathbf{I}\otimes\mathbf{I}) \bs{\vec{r}}
\end{pmatrix},
\end{equation}
where $\mathbf{I}$ and $\mathbf{0}$ are short for $\mathbf{I}_{N-1}$ and $\mathbf{0}_{N-1}$ respectively.
Let us denote 
\begin{equation}
\label{eq:UFVec2d}
\begin{aligned}
	& 
	\bs{\vec{U}^{x}} = ( \mathbf{I}_{N-1}\otimes\bs{\ushort{e}}_{1}^{\intercal} ) \bs{\vec{u}}, \quad  
	\bs{\vec{U}^{y}} = ( \bs{\ushort{e}}_{1}^{\intercal}\otimes\mathbf{I}_{N-1} ) \bs{\vec{v}}, \quad 
	\bs{\vec{F}^{x}} = ( \mathbf{I}_{N-1}\otimes\bs{\ushort{e}}_{1}^{\intercal} )  \bs{\vec{f}}, \quad 
	\bs{\vec{F}^{y}} = ( \bs{\ushort{e}}_{1}^{\intercal}\otimes\mathbf{I}_{N-1} ) \bs{\vec{g}}, 
	\\
	& 
	\bs{\vec{U}^{1}} = ( \mathbf{I}_{N-1}\otimes\mathbf{\ushort{E}} ) \bs{\vec{u}}, \quad 
	\bs{\vec{U}^{2}} = ( \mathbf{\ushort{E}}\otimes\mathbf{I}_{N-1} ) \bs{\vec{v}}, \quad 
	\bs{\vec{F}^{1}} = ( \mathbf{I}_{N-1}\otimes\mathbf{\ushort{E}} ) \bs{\vec{f}}, \quad 
	\bs{\vec{F}^{2}} = ( \mathbf{\ushort{E}}\otimes\mathbf{I}_{N-1} ) \bs{\vec{g}}.
\end{aligned}
\end{equation}
As a result, the coupled linear system \eqref{eq:MaxwAlgSys2d} can be decomposed into the following form
		\begin{equation}
		\label{eq:MaxwFinSys2d}
		\begin{pmatrix}
		\mathbf{I}+{\kappa}\mathbf{M}   &  &   &  & 
		\\
		&  \mathbf{I} + {\kappa}\mathbf{M} &    &  & 
		\\
		&  & ( \mathbf{I} + {\kappa}\mathbf{M} )\otimes\mathbf{I} & -\mathbf{I}\otimes\mathbf{I}  & \mathbf{M}\otimes\mathbf{I}
		\\
		&  & -\mathbf{I}\otimes\mathbf{I}   & \mathbf{I}\otimes(\mathbf{I}+{\kappa}\mathbf{M})     & \mathbf{I}\otimes\mathbf{M} 
		\\		
		&  & \mathbf{M}\otimes\mathbf{I}  &  \mathbf{I}\otimes\mathbf{M}  &   \mathbf{0}\otimes\mathbf{0}
		\end{pmatrix}
		\begin{pmatrix}
		\bs{\vec{U}^{x}} \\
		\bs{\vec{U}^{y}}  \\
		\bs{\vec{U}^{1}} \\
		\bs{\vec{U}^{2}} \\
		\bs{\vec{p}}
		\end{pmatrix}
		=\begin{pmatrix}
		\bs{\vec{F}^{x}} \\
		\bs{\vec{F}^{y}}  \\
		\bs{\vec{F}^{1}} \\
		\bs{\vec{F}^{2}} \\
		\bs{\vec{r}}
		\end{pmatrix},
		\end{equation}
where $\bs{\vec{U}^{x}}$ and $\bs{\vec{U}^{y}}$ are fully decoupled and can be computed independently. 
To facilitate the design of efficient algorithm for the sub-system involving $\bs{\vec{U}^{1}}$, $\bs{\vec{U}^{2}}$ and $\bs{\vec{p}}$, we adopt the notations
\begin{equation}
\label{eq:UsVec2Mat2d}
\mathbf{U}^s={\rm ivec}(\bs{\vec{U}^{s}} ), \quad \mathbf{F}^s={\rm ivec}(\bs{\vec{F}^{s}} ),\;\; s=1,2, \quad 
\mathbf{P}={\rm ivec}(\bs{\vec{p}} ),\quad  \mathbf{R}={\rm ivec}(\bs{\vec{r}} )
\end{equation}
and reformulate the sub-problem into an equivalent form of matrix equations 
\begin{subequations}
\label{eq:2dSubMatEq}
\begin{align}
		&  \mathbf{U}^1 - \mathbf{U}^2   + ( {\kappa} \mathbf{U}^1 + \mathbf{P} ) \mathbf{M} = \mathbf{F}^1,
		\label{eq:2dSubMatEq1} 
		\\
		&  -\mathbf{U}^1 + \mathbf{U}^2  + \mathbf{M} ( {\kappa} \mathbf{U}^2 + \mathbf{P} )  = \mathbf{F}^2,  
		\label{eq:2dSubMatEq2}
		\\
		&  \mathbf{U}^1\mathbf{M} + \mathbf{M}\mathbf{U}^2  = \mathbf{R}, 
		\label{eq:2dSubMatEq3}
\end{align}		
\end{subequations} 
where we have used the identity of Kronecker product of matrices that $(\mathbf{A}\otimes \mathbf{B})\bs{\vec{u}}={\rm vec}(\mathbf{B} \mathbf{U} \mathbf{A}^{\intercal})$ 
if $\bs{\vec{u}}={\rm vec} (\mathbf{U})$.

\textbf{Step 2: Solving the decoupled sub-problems via eigen-decompostion.} 
Next, we propose a highly efficient solution algorithm for the coupled sub-system \eqref{eq:2dSubMatEq} with a semi-analytic approach based on eigen-decomposition of matrix $\mathbf{M}$.
\begin{lm}[cf. \cite{Shen-1994b}]
\label{Lem2-MIdiag}
The symmetric penta-diagonal matrix $\mathbf{M}\in\mathbb{R}^{(N-1)\times(N-1)}$ can be diagonalized with $\mathcal{O}(N^2)$ operations into
\begin{equation}
\label{eq:MatMIdiag}
\mathbf{M} =  \mathbf{Q}\mathbf{D}\mathbf{Q}^{\intercal},
\quad
\mathbf{Q}\mathbf{Q}^{\intercal} = \mathbf{Q}^{\intercal}\mathbf{Q}=\mathbf{I},
\end{equation}
where $\mathbf{D}={\rm diag}(d_1,\cdots,d_{N-1})$ is a diagonal matrix with its main diagonal entries the eigenvalues of matrix $\mathbf{M}$, 
and $\mathbf{Q}$ is the corresponding orthonormal eigenvector matrix. 
Here, an operation denotes as either floating-point addition or multiplication.
\end{lm}

One can insert the third equation into the summation of the first two equations of system \eqref{eq:2dSubMatEq} to obtain
\begin{equation}
\label{eq:PMatEq2d}
\mathbf{P}\mathbf{M} + \mathbf{M}\mathbf{P} = \mathbf{F}^{1}+\mathbf{F}^{2}-\kappa \mathbf{R}.  
\end{equation}
And then, by introducing auxiliary matrices
\begin{equation}
\label{eq:MatPFRhat}
\mathbf{\hat{P}}=\mathbf{Q}^{\intercal}\mathbf{P}\mathbf{Q},\quad \mathbf{\hat{F}}^1=\mathbf{Q}^{\intercal}\mathbf{F}^1\mathbf{Q},
\quad \mathbf{\hat{F}}^2=\mathbf{Q}^{\intercal}\mathbf{F}^2\mathbf{Q},\quad \mathbf{\hat{R}}=\mathbf{Q}^{\intercal}\mathbf{R}\mathbf{Q}
\end{equation}
into the above equation, one arrives at
\begin{equation*}
\mathbf{\hat{P}} \mathbf{D} + \mathbf{D} \mathbf{\hat{P}} = \mathbf{\hat{F}}^1 +\mathbf{\hat{F}}^2-\kappa \mathbf{\hat{R}} \;\; \Longleftrightarrow \;\; 
\hat{P}_{ij}d_j + d_i\hat{P}_{ij} = \hat{F}^1_{ij}+\hat{F}^2_{ij}-\kappa \hat{R}_{ij},
\end{equation*}
which directly leads to semi-analytic formula for $\mathbf{P}$:
\begin{equation}
\label{eq:PSol2d}
\mathbf{P}=\mathbf{Q}\mathbf{\hat P}\mathbf{Q}^{\intercal},\quad \hat{P}_{ij}=\frac{ {\hat{F}}_{ij}^1 +{\hat{F}}_{ij}^2-\kappa {\hat{R}}_{ij}   }{d_i+d_j}.
\end{equation}
Once $\mathbf{P}$ is obtained, one can express $\mathbf{U}^2$ in terms of $\mathbf{U}^1$ using equation \eqref{eq:2dSubMatEq3}, 
and insert it into equation \eqref{eq:2dSubMatEq1} to obtain the following fully decoupled equation for $\mathbf{U}^1$:
\begin{equation}
\label{eq:U1MatEq2d}
\mathbf{U}^1 + \mathbf{M}^{-1}\mathbf{U}^1\mathbf{M} + {\kappa} \mathbf{U}^1\mathbf{M} 
= \mathbf{F}^1 - \mathbf{P}\mathbf{M} + \mathbf{M}^{-1}\mathbf{R}.
\end{equation}
With a similar manner, one can also derive the decoupled equation for  $\mathbf{U}^2$:
\begin{equation}
\label{eq:U2MatEq2d}
\mathbf{U}^2 + \mathbf{M}\mathbf{U}^2\mathbf{M}^{-1} + {\kappa}\mathbf{M}\mathbf{U}^2 = \mathbf{F}^2 - \mathbf{M}\mathbf{P} + \mathbf{R}\mathbf{M}^{-1}.
\end{equation}
Let us introduce the auxiliary matrices $\mathbf{\hat{U}}^{s}=\mathbf{Q}^{\intercal}\mathbf{U}^{s}\mathbf{Q}$ with $s=1,2$ 
such that equations \eqref{eq:U1MatEq2d}-\eqref{eq:U2MatEq2d} can be transformed into
\begin{equation*}
\begin{aligned}
& \mathbf{\hat{U}}^1 + \mathbf{D}^{-1}\mathbf{\hat{U}}^1\mathbf{D} + {\kappa} \mathbf{\hat{U}}^1\mathbf{D} 
=  \mathbf{\hat{F}}^1 -\mathbf{\hat{P}} \mathbf{D} +\mathbf{D}^{-1}  \mathbf{\hat{R}} \;\; \Longleftrightarrow  \;\;
\hat{U}^1_{ij}    + d_i^{-1}\hat{U}^1_{ij}d_j +\kappa \hat{U}^1_{ij}d_j  = \hat{F}^1_{ij}-\hat{P}_{ij}d_j+  d_i^{-1}\hat{R}_{ij},
\\
& \mathbf{\hat{U}}^2 + \mathbf{D} \mathbf{\hat{U}}^2\mathbf{D}^{-1} + {\kappa}\mathbf{D} \mathbf{\hat{U}}^2 
=  \mathbf{\hat{F}}^2 -\mathbf{D}  \mathbf{\hat{P}} + \mathbf{\hat{R}}\mathbf{D}^{-1}  \;\;  \Longleftrightarrow \;\;
\hat{U}^2_{ij}    + d_i\hat{U}^2_{ij}d_j^{-1} +\kappa d_i \hat{U}^2_{ij}  = \hat{F}^2_{ij}-d_i\hat{P}_{ij}+  \hat{R}_{ij}d_j^{-1},
\end{aligned}
\end{equation*}
which directly give
\begin{equation}
\label{eq:U12Sol2d}
\mathbf{U}^{s}=\mathbf{Q}\mathbf{\hat{U}}^{s}\mathbf{Q}^{\intercal}, \;\; s=1,2, \quad 
\hat{U}_{ij}^1 = \frac{\hat{F}_{ij}^1 - {\hat{P}}_{ij}d_j  + d_i^{-1}{\hat{R}}_{ij}}{1+d_i^{-1} d_j + {\kappa}d_j}, \quad 
\hat{U}_{ij}^2 = \frac{\hat{F}_{ij}^2 - d_i\hat{P}_{ij} + \hat{R}_{ij}d_j^{-1}}{1+d_i d_j^{-1} + {\kappa}d_i}.
\end{equation}
As for $\bs{\vec{U}^{x}}$-$\bs{\vec{U}^{y}}$ in equation \eqref{eq:MaxwFinSys2d}, each of them satisfies a fully decoupled equation and can be expressed analytically as
\begin{equation}
\label{eq:UxySol2d}
\bs{\vec{U}^{x}} = \mathbf{Q}\bs{{\hat{U}}^{x}}, \; \;  \bs{\vec{U}^{y}} = \mathbf{Q}\bs{{\hat{U}}^{y}}, \; \; 
\hat{U}_{j}^{x} = \frac{(\mathbf{Q}^{\intercal}\bs{\vec{F}^{x}})_j}{(1+{\kappa}d_j)}, \; \;
\hat{U}_{j}^{y} = \frac{(\mathbf{Q}^{\intercal}\bs{\vec{F}^{y}})_j}{(1+{\kappa}d_j)}, \; \; j=1,\cdots, N-1.
\end{equation}

To conclude, let us gather the semi-analytic expressions of the solution for system \eqref{eq:MaxwAlgSys2d}  and summarize the solution procedure in Algorithm \ref{Algo1}.
\begin{algorithm}[ht!]
	\caption{Matrix-free, semi-analytic solution algorithm for $\bs{\vec{u}}$, $\bs{\vec{v}}$ and $\bs{\vec{p}}$ in \eqref{eq:MaxwAlgSys2d} with Gauss' law constraint}
	\label{Algo1}
	\begin{algorithmic}
		\item  
		[{\bf Input:}] 
		The polynomial order $N$, parameter $\kappa$,  $\mathbf{F}$, $\mathbf{G}$ and $\mathbf{R}$ in \eqref{eq:MatFGR2d}, 
		$\ushort{\bs{e}}_1$ and $\ushort{\mathbf{E}}$ in \eqref{eq:MatINE}, $\mathbf{D}$, $\mathbf{Q}$ in \eqref{eq:MatMIdiag}.
		\vspace{3pt} 
		
		\item 
		[{\bf Output:}] 
		The solutions $\bs{\vec{u}}$, $\bs{\vec{v}}$ and $\bs{\vec{p}}$.
		\vspace{3pt} 
		\item[1:] Compute $\bs{\vec{F}^{x}}$-$\bs{\vec{F}^{y}}$ and $\mathbf{F}^{1}$-$\mathbf{F}^{2}$ by \eqref{eq:UFVec2d} in a matrix-free manner, i.e.,
		\begin{equation*}
		\bs{\vec{F}^{x}}=\mathbf{F}^{\intercal}\bs{\ushort{e}}_{1}, \quad \bs{\vec{F}^{y}}=\mathbf{G}\bs{\ushort{e}}_{1},
		\quad \mathbf{F}^1=\mathbf{\ushort{E}}\mathbf{F},\quad \mathbf{F}^2= \mathbf{G}\mathbf{\ushort{E}}^{\intercal}.
		\end{equation*}
		
		\item [2:] Compute $\mathbf{\hat{F}}^1$, $\mathbf{\hat{F}}^2$ and $\mathbf{\hat{R}}$ by \eqref{eq:MatPFRhat}. 
		
		\item [3:] Compute $\mathbf{P}$ by \eqref{eq:PSol2d}.
		
		\item [4:] Compute $\mathbf{U}^1$ and $\mathbf{U}^2$ by \eqref{eq:U12Sol2d}.
		
		\item [5:] Compute $\bs{\vec{U}^{x}}$ and $\bs{\vec{U}^{y}}$ by \eqref{eq:UxySol2d}.
		
		\item [6:] Assemble $\mathbf{U}$ and $\mathbf{V}$ by $\mathbf{U}=[(\bs{\vec{U}^{x}})^{\intercal};\mathbf{U}^1]$ and $\mathbf{V}=[\bs{\vec{U}^{y}}, \mathbf{U}^2]$.
		
		\item [7:] Return $\bs{\vec{u}}={\rm vec}(\mathbf{U})$, $\bs{\vec{v}}={\rm vec}(\mathbf{V})$ and $\bs{\vec{p}}={\rm vec}(\mathbf{P})$.
	\end{algorithmic}
\end{algorithm}

\begin{rem}
	\label{Remk1-Algo1}
	The above algorithm takes advantage of tensor-product nature of the proposed arbitrary order $\bs{H}({\rm curl})$-conforming spectral basis functions, 
	which avoids the assembly of global coefficient matrix and requires only several 1D matrix-matrix multiplications to obtain the numerical solution. 
	Note that if matrices $\mathbf{A}$, $\mathbf{B}$, $\mathbf{X}\in\mathbb{R}^{N\times{N}}$, 
	matrix-matrix product $\mathbf{Y}=\mathbf{B}\mathbf{X}\mathbf{A}^{\intercal}$ requires $2N^2(2N-1)$ operations or $\mathcal{O}(N^3)$ to evaluate, 
	while the disregard of Kronecker structure in $\bs{\vec{y}}=(\mathbf{A}\otimes\mathbf{B})\bs{\vec{x}}$ with $\bs{\vec{y}}={\rm vec}(\mathbf{Y})$ and $\bs{\vec{x}}={\rm vec}(\mathbf{X})$
	leads to the cost of $N^2(2N^2-1)$ operations.  
	Thus, the computational cost is reduced from $\mathcal{O}(N^6)$ to $\mathcal{O}(N^3)$, when compared with direct solution algorithm based on LU factorization, 
	and can further be reduced to $\mathcal{O}(N^{\log_{2}7})$ 
	with the help of Strassen's matrix multiplication algorithm \cite{Strassen-1969}, as verified by Example \ref{Exmp5-1} in Section \ref{Sec5}.
\end{rem}

\subsection{Double-curl eigenvalue problems in 2D}
\label{Sec2-3} 
Now, with a slight deviation from the solution algorithm of the double-curl source problem \eqref{eq:MaxwSouProb}, 
one may first ignore the divergence-free condition in the eigenvalue problem  \eqref{eq:MaxwEigProb} and directly resort to the weak formulation: 
find $\lambda\in\mathbb{R}$ and $\bs{u}\in\bs{H}_0({\rm curl};\Omega)$ with $\bs{u}\neq\bs{0}$ such that 
\begin{equation}
\label{eq:EigWeak}
\begin{aligned}
& (\nabla\times\bs{u},\nabla\times\bs{v})  = \lambda (\bs{u},\bs{v}), \quad \forall\bs{v}\in\bs{H}_{0}({\rm curl};\Omega). 
\end{aligned}
\end{equation} 
Then the approximation scheme for the weak formulation \eqref{eq:EigWeak} in two dimensions with $\Omega={\Lambda}^2$ takes the form: 
find ${\lambda}^{N}\in\mathbb{R}$ and $\bs{u}_{N}\in\bs{H}_{N,0}({\rm curl};{\Lambda}^2)$ with $\bs{u}_N\neq\bs{0}$ such that
\begin{equation} 
\label{eq:EigDis2d}
\begin{aligned}
& (\nabla\times\bs{u}_N,\nabla\times\bs{v}) = {\lambda}^N (\bs{u}_N,\bs{v}) ,\quad \forall \bs{v}\in\bs{H}_{N,0}({\rm curl};{\Lambda}^2). 
\end{aligned}
\end{equation}

Let us adopt the spectral approximations of $\bs{u}_N$ in \eqref{eq:uNpN2d} and 
the notations $\bs{\vec{u}}={\rm vec}(\mathbf{U})$ and $\bs{\vec{v}}={\rm vec}(\mathbf{V})$ as the unknown vectors in \eqref{eq:upVec2d}. 
Following the same derivation of \eqref{eq:MaxwAlgSys2d}, one obtains the following algebraic eigenvalue problem of \eqref{eq:EigDis2d}:
\begin{equation}
\label{eq:EigAlgSys2d} 
\begin{pmatrix}
\mathbf{I}_{N-1}\otimes{\mathbf{I}_N}
& - \mathbf{\ushort{E}}\otimes\mathbf{\ushort{E}}^{\intercal}
\\
-\mathbf{\ushort{E}}^{\intercal}\otimes\mathbf{\ushort{E}}
& {\mathbf{I}_N}\otimes\mathbf{I}_{N-1}
\end{pmatrix}
\begin{pmatrix}
\bs{\vec{u}} 
\\
\bs{\vec{v}} 
\end{pmatrix}
= {\lambda}^{N}
\begin{pmatrix}
\mathbf{M}\otimes{\mathbf{I}_N}
& 
\\
& {\mathbf{I}_N}\otimes\mathbf{M}
\end{pmatrix}
\begin{pmatrix}
\bs{\vec{u}} 
\\
\bs{\vec{v}} 
\end{pmatrix}.
\end{equation}

The proposed semi-analytic computational approach is still valid for solving the above algebraic system efficiently. 
By left multiplying the first 4-by-2 block submatrix of the permutation matrix given in \eqref{eq:MaxwPemMat2d} 
and adopting the notations $\mathbf{U}^1$, $\mathbf{U}^2$, $\bs{\vec{U}^x}$ and $\bs{\vec{U}^y}$ in \eqref{eq:UFVec2d} and \eqref{eq:UsVec2Mat2d}, 
one obtains an equivalent form of system \eqref{eq:EigAlgSys2d} as follows: 
\begin{equation}
\label{eq:EigMatEq2d} 
\begin{aligned}
& \mathbf{U}^1 - \mathbf{U}^2  = {\lambda}^{N} \mathbf{U}^1\mathbf{M},  \quad 
-\mathbf{U}^1 + \mathbf{U}^2  = {\lambda}^{N} \mathbf{M}\mathbf{U}^2,  
\\
&  \bs{\vec{U}^{x}} = {\lambda}^{N}  \mathbf{M} \bs{\vec{U}^{x}},  \qquad  \qquad
\bs{\vec{U}^{y}} = {\lambda}^{N}  \mathbf{M} \bs{\vec{U}^{y}}.
\end{aligned}
\end{equation}

By summing up the first two equations of \eqref{eq:EigMatEq2d}, one has 
\begin{equation}
\label{eq:EigU12Rel2d}
{\lambda}^{N} ( \mathbf{U}^1\mathbf{M} + \mathbf{M}\mathbf{U}^2 ) = \mathbf{0}
\;\; \Longleftrightarrow  \;\; 
{\lambda}^{N} = 0 \;\; \text{or} \;\;  \mathbf{U}^1\mathbf{M} + \mathbf{M}\mathbf{U}^2  = \mathbf{0}.
\end{equation}

On the one hand, when ${\lambda}^{N}=0$, one obtains the zero eigenvalues and the coresponding eigenfunctions of \eqref{eq:EigMatEq2d} as follows:
\begin{align*}
{\lambda}_{ij}^{N} = 0, \quad 
\mathbf{U}_{ij}^{1} = \mathbf{U}_{ij}^{2} = \mathbf{I}_{ij}, \quad
\bs{\vec{U}^{x}} = \bs{\vec{U}^{y}} = \bs{\vec{0}}, 
\end{align*}
where $\mathbf{I}_{ij}=\bs{e}_{i}\otimes\bs{e}_{j}^{\intercal}$ denotes a matrix of size $(N-1)$-by-$(N-1)$ with $(i,j)$-th element being $1$ and the rest elements being $0$.

On the other hand, when ${\lambda}^{N}\neq0$,  one has $\mathbf{U}^1\mathbf{M} + \mathbf{M}\mathbf{U}^2  = \mathbf{0}$ by \eqref{eq:EigU12Rel2d}.
Then the first two equations of \eqref{eq:EigMatEq2d} yields the following decoupled equations for $\mathbf{U}^{1}$ and $\mathbf{U}^{2}$:
\begin{align}
\label{eq:UsMatEq2d}
\mathbf{M}\mathbf{U}^{1} + \mathbf{U}^{1}\mathbf{M} = {\lambda}^{N}\mathbf{M}\mathbf{U}^{1}\mathbf{M} \quad \text{and} \quad
\mathbf{M}\mathbf{U}^{2} + \mathbf{U}^{2}\mathbf{M} = {\lambda}^{N}\mathbf{M}\mathbf{U}^{2}\mathbf{M},
\end{align}
respectively. 
Introducing the auxiliary matrices $\mathbf{\hat{U}}^{s}=\mathbf{Q}^{\intercal}\mathbf{U}^{i}\mathbf{Q}$ with $s=1,2$ 
such that equations \eqref{eq:UsMatEq2d} can be transformed into 
\begin{align}
\label{eq:UshatMatEq2d}
& \mathbf{D}\mathbf{\hat{U}}^{s} + \mathbf{\hat{U}}^{s}\mathbf{D} = {\lambda}^{N} \mathbf{D}\mathbf{\hat{U}}^{s}\mathbf{D}, \; s=1,2
\;\; \Longleftrightarrow  \;\;
d_i \hat{U}_{ij}^{s} + \hat{U}_{ij}^{s}d_j = {\lambda}_{ij}^{N} d_i\hat{U}_{ij}^{s}d_j, \; s=1,2,
\end{align}
which directly gives rise to the non-zero eigenvalues and the corresponding eigenfunctions of the first two equations of \eqref{eq:EigMatEq2d} as follows:
\begin{align*}
\lambda_{ij}^{N} =  \frac{1}{d_i} + \frac{1}{d_j}, \quad 
& \mathbf{U}_{ij}^{1} = \mathbf{Q}\mathbf{I}_{ij}\mathbf{Q}^{\intercal}, \quad 
\mathbf{U}_{ij}^{2} = -\mathbf{Q}({d_i}/{d_j}\mathbf{I}_{ij})\mathbf{Q}^{\intercal}.
\end{align*}

As for $\bs{\vec{U}^{x}}$-$\bs{\vec{U}^{y}}$ in equation \eqref{eq:EigMatEq2d}, both of them satisfy decoupled equations of the form:
\begin{equation*}
( {\lambda}^{N}\mathbf{M}-\mathbf{I} ) \bs{\vec{U}^{x}} = \bs{\vec{0}}, \quad ( {\lambda}^{N}\mathbf{M}-\mathbf{I} ) \bs{\vec{U}^{y}} = \bs{\vec{0}}. 
\end{equation*}
It is straightforward to find that the  ${\lambda}_{i}^{N}$ ($i=1,2,\cdots,N-1$) is the reciprocal of the $i$-th eigenvalue of matrix $\mathbf{M}$, 
and the associated egienfunction solution $\bs{\vec{U}^{x}}$(resp. $\bs{\vec{U}^{y}}$) is just the $i$-th column of matrix $\mathbf{Q}$.

To conclude, we summarize the semi-analytic solution algorithm for the eigenvalue problem \eqref{eq:EigAlgSys2d} in Theorem \ref{Thm1-LamU2dSol} as follows.
\begin{thm}
\label{Thm1-LamU2dSol}
The exact eigen-solutions $({\lambda}^{N},\bs{\vec{u}},\bs{\vec{v}})$ of system \eqref{eq:EigAlgSys2d} shall be divided into two categories:
\begin{itemize}
	\item The eigen-solution $(\lambda_{ij}^{N},\mathbf{U}_{ij},\mathbf{V}_{ij}) $ associated with the non-zero eigenvalues  includes:
	\begin{itemize}
	\item[---] 
	For $1\leq{i,j}\leq{N-1}$, 
	\begin{equation}
		\label{eq:LamU2dSolIn}
		\lambda_{ij}^{N} = \frac{1}{d_i} + \frac{1}{d_j}, \quad  
		\mathbf{U}_{ij} = [\bs{\vec{0}}^{\intercal}; \mathbf{U}_{ij}^{1}], \quad  
		\mathbf{V}_{ij} = [\bs{\vec{0}},\mathbf{U}_{ij}^{2}], \quad 
	\end{equation}
	with
	\begin{align}
	\label{eq:U2dSolIn}
	&\mathbf{U}_{ij}^{1} = \mathbf{Q}\mathbf{I}_{ij}\mathbf{Q}^{\intercal}, \quad 
	\mathbf{U}_{ij}^{2} = -\mathbf{Q}({d_i}/{d_j}\mathbf{I}_{ij})\mathbf{Q}^{\intercal}.
	\end{align}
	
	\item[---]
	For $i=0$ and $1\leq{j}\leq{N-1}$, 
	\begin{equation}
	\label{eq:LamU2dSolFy}
		\lambda_{0j}^{N} = \frac{1}{d_j}, \quad 
		\mathbf{U}_{0j} = [(\mathbf{U}_{0j}^{1})^{\intercal}; \mathbf{0}_{N-1}],\quad 
		\mathbf{V}_{0j} = [\mathbf{U}_{0j}^{2},\mathbf{0}_{N-1}],\quad
		\mathbf{U}_{0j}^{1} = \mathbf{Q}(:,j), \quad
		\mathbf{U}_{0j}^{2} = \bs{\vec{0}}. 
	\end{equation}
	
	\item[---]
	For $j=0$ and $1\leq{i}\leq{N-1}$, 
	\begin{equation}
	\label{eq:LamU2dSolFx}
	\lambda_{i0}^{N} = \frac{1}{d_i}, \quad 
	\mathbf{U}_{i0} = [(\mathbf{U}_{i0}^{1})^{\intercal}; \mathbf{0}_{N-1}], \quad 
	\mathbf{V}_{i0} = [\mathbf{U}_{i0}^{2},\mathbf{0}_{N-1}],\quad
	\mathbf{U}_{i0}^{1} = \bs{\vec{0}}, \quad
	\mathbf{U}_{i0}^{2} = \mathbf{Q}(:,i).
	\end{equation}
	\end{itemize}

	\item  
	The eigen-solution $({\lambda}_{ij}^{N},\mathbf{U}_{ij},\mathbf{V}_{ij})$ associated with the zero eigenvalues takes the form
	\begin{equation}
	\label{eq:Lam0U2dSol}
	\lambda_{ij}^{N} = 0, \quad 
	\mathbf{U}_{ij} = [\bs{\vec{0}}^{\intercal}; \mathbf{U}_{ij}^{1}], \quad 
	\mathbf{V}_{ij} = [\bs{\vec{0}},\mathbf{U}_{ij}^{2}], \quad
	\mathbf{U}_{ij}^{1}= \mathbf{U}_{ij}^{2} = \mathbf{I}_{ij}. 
	\end{equation}
\end{itemize}
\end{thm}

In fact, similar to the double-curl source problem \eqref{eq:MaxwSouProb}, one also can treat the eigenvalue problem \eqref{eq:MaxwEigProb} with mixed formulation 
and guarantee the divergence-free condition $\nabla\cdot\bs{u}=0$ explicitly. 
The corresponding approximation scheme takes the form: 
find ${\lambda}^{N}\in\mathbb{R}$, $\bs{u}_{N}\in\bs{H}_{N,0}({\rm curl};{\Lambda}^2)$ with $\bs{u}_N\neq\bs{0}$ and $p_N\in{H}_{N,0}^{1}({\Lambda}^2)$ such that
\begin{equation}
\label{eq:EigMDis2d}
\begin{aligned}
& (\nabla\times\bs{u}_N,\nabla\times\bs{v}) + (\nabla{p_N},\bs{v}) = {\lambda}^N (\bs{u}_N,\bs{v}) ,\quad \forall \bs{v}\in\bs{H}_{N,0}({\rm curl};{\Lambda}^2),
\\
& (\bs{u}_N,\nabla{q}) = 0, \quad \forall q\in {H}_{N,0}^{1}({\Lambda}^2).
\end{aligned}
\end{equation}
Following the same derivation of \eqref{eq:EigAlgSys2d} and \eqref{eq:EigMatEq2d}, one obtains the following algebraic eigenvalue problem of \eqref{eq:EigMDis2d}:
\begin{equation}
\label{eq:EigMAlgSys2d}
\begin{pmatrix}
\mathbf{I}_{N-1}\otimes{\mathbf{I}_N}
& - \mathbf{\ushort{E}}\otimes\mathbf{\ushort{E}}^{\intercal}
& \mathbf{M}\otimes\mathbf{\ushort{E}}^{\intercal}
\\
-\mathbf{\ushort{E}}^{\intercal}\otimes\mathbf{\ushort{E}}
& {\mathbf{I}_N}\otimes\mathbf{I}_{N-1}
& \mathbf{\ushort{E}}^{\intercal}\otimes\mathbf{M}
\\
\mathbf{M}\otimes\mathbf{\ushort{E}}
& \mathbf{\ushort{E}}\otimes\mathbf{M}
& \mathbf{0}\otimes\mathbf{0} 
\end{pmatrix}
\begin{pmatrix}
\bs{\vec{u}} \\
\bs{\vec{v}} \\
\bs{\vec{p}}
\end{pmatrix}
= {\lambda}^{N}
\begin{pmatrix}
\mathbf{M}\otimes{\mathbf{I}_N}
& &
\\
& {\mathbf{I}_N}\otimes\mathbf{M}
&
\\
& 
& \mathbf{0}\otimes\mathbf{0} 
\end{pmatrix}
\begin{pmatrix}
\bs{\vec{u}} \\
\bs{\vec{v}} \\
\bs{\vec{p}}
\end{pmatrix},
\end{equation}
and the equivalent form of system \eqref{eq:EigMAlgSys2d}:
\begin{equation}
\label{eq:EigMMatEq2d} 
\begin{aligned}
& \mathbf{U}^1 - \mathbf{U}^2   + \mathbf{P}\mathbf{M} = {\lambda}^{N} \mathbf{U}^1\mathbf{M},  \quad
-\mathbf{U}^1 + \mathbf{U}^2  + \mathbf{M}\mathbf{P} = {\lambda}^{N} \mathbf{M}\mathbf{U}^2,  \quad 
\mathbf{U}^1\mathbf{M} + \mathbf{M}\mathbf{U}^2 = \mathbf{0},
\\
& \bs{\vec{U}^{x}} = {\lambda}^{N}  \mathbf{M} \bs{\vec{U}^{x}},  \quad
\bs{\vec{U}^{y}} = {\lambda}^{N}  \mathbf{M} \bs{\vec{U}^{y}}.
\end{aligned}
\end{equation}
By summing up the first two equations and utlizing the third equation of system \eqref{eq:EigMMatEq2d},  one has 
\begin{equation*}
\mathbf{P}\mathbf{M} + \mathbf{M}\mathbf{P} = {\lambda}^{N}( \mathbf{U}^1\mathbf{M}+\mathbf{M}\mathbf{U}^2) = \mathbf{0}  
\;\; \Longleftrightarrow \;\; ( \mathbf{M}\otimes\mathbf{I} + \mathbf{I}\otimes\mathbf{M} ) \bs{\vec{p}} = \bs{\vec{0}}.
\end{equation*}
Since $\mathbf{M}\otimes\mathbf{I} + \mathbf{I}\otimes\mathbf{M} $ is full-rank, one directly obtains $\bs{\vec{p}}=\bs{0}$, i.e., $\mathbf{P}=\mathbf{0}$.
This eliminates $\mathbf{P}$ from the first two equations of system \eqref{eq:EigMatEq2d}, which gives 
\begin{align*}
& \mathbf{U}^1 - \mathbf{U}^2  = {\lambda}^{N} \mathbf{U}^1\mathbf{M}, \quad  
-\mathbf{U}^1 + \mathbf{U}^2  = {\lambda}^{N} \mathbf{M}\mathbf{U}^2.
\end{align*}
Due to the above equations is just the first two equations of \eqref{eq:EigMatEq2d}, 
the eigen-solutions of \eqref{eq:EigMAlgSys2d} coresponding to non-zero eigenvalues are equivalent to the ones of \eqref{eq:EigAlgSys2d}.

\begin{rem}
\label{Remk2-Thm1}
The proposed algorithm strictly obeys the Helmholtz-Hodge decomposition (cf. \cite{BoffiBrezziFortin2013,Monk2003}), 
i.e., for any $\bs{u}\in \bs{H}_{0}({\rm curl};\Omega)$, there exists a unique scalar potential $\phi\in H_0^1(\Omega)$ 
and a unique vector potential $\bs{\psi}\in\bs{H}({\rm curl};\Omega)\cap\bs{H}_0({\rm div};\Omega)$ with $\nabla\cdot\bs\psi=0$ such that
\begin{equation*}
\bs{u} = \nabla\phi + \nabla\times\bs{\psi}.
\end{equation*}
This can be seen clearly with the following observations:
\begin{itemize}
\item
When ${\lambda}^{N}\neq0$, the constraint $\mathbf{U}^1\mathbf{M} + \mathbf{M}\mathbf{U}^2 = \mathbf{0}$ in \eqref{eq:EigU12Rel2d} or \eqref{eq:EigMMatEq2d} 
is exactly the divergence-free condition satisfied by the numerical eigen-functions.
As a result, all eigen-solutions associated with non-zero eigenvalues in Theorem \ref{Thm1-LamU2dSol} are exactly divergence-free solutions at the discrete level.

\item 
When ${\lambda}^{N}=0$, the eigenfunctions in \eqref{eq:Lam0U2dSol} correspond to
\begin{equation*}
\bs{u} = \bs{\Phi}_{m,n}^1(\bs{x}) + \bs{\Phi}_{m,n}^2(\bs{x}) =
\big( \phi_{m}(x_1) \psi_{n+1}(x_2), \psi_{m+1}(x_1) \phi_{n}(x_2) \big)^{\intercal}=\nabla\big(\psi_{m+1}(x_1) \psi_{n+1}(x_2)\big),
\end{equation*}
where $\big\{{\Psi}_{m,n}(\bs{x})=\psi_{m+1}(x_1) \psi_{n+1}(x_2)\big\}_{m,n=1}^{+\infty}$ is the complete basis for $H_0^1(\Omega)$.
\end{itemize}	
This conclusion is also valid for the double-curl eigenvalue problem in 3D, we refer to Remark \ref{Remk4-Thm2} in Subsection \ref{Sec3-3} for more details.
\end{rem}

\section{Maxwell's double-curl problems in three dimensions} 
\label{Sec3}
In this section, we shall extend the proposed Gauss's law-preserving spectral methods and their efficient solution algorithms to solve the double-curl source and eigenvalue problems in three dimensions. 
Though the design philosophy resembles that of the 2D case, the approximation scheme, solution algorithm and implementation technique are undoubtedly more involved and worth detailed explanation.

\subsection{Spectral approximation scheme for double-curl source problem \texorpdfstring{\eqref{eq:MaxwWeak}}{Lg}}
\label{Sec3-1}
Hereinafter, we fix $\Omega:={\Lambda}^3=(-1,1)^3$ and propose the approximation pair $\bs{H}_{N,0}({\rm curl};{\Lambda}^3) \times{H}_{N,0}^1({\Lambda}^3)$ for $(\bs{u}, p)$ as follows. 
\begin{prop}
	\label{Prop2-DiscSpace3d}
	The conforming spectral approximation spaces $\bs{H}_{N,0}({\rm curl};{\Lambda}^3)$ and ${H}_{N,0}^1({\Lambda}^3)$ of order $N$ for 
	$\bs{H}_{0}({\rm curl};{\Lambda}^3)$ and ${H}_{0}^{1}({\Lambda}^3)$ are defined respectively by
	\begin{equation}
	\label{eq:HN0curl3d}
	\bs{H}_{N,0}({\rm curl};{\Lambda}^3) = {\rm span} \big( 
	\{ \bs\Phi_{m,n,l}^1 \}, 
	\{ \bs\Phi_{m,n,l}^2 \}, 
	\{ \bs\Phi_{m,n,l}^3 \} \big), 
	\quad  {H}_{N,0}^1({\Lambda}^3) = {\rm span}
	\big(  \big\{ \Psi_{m,n,l} \big\}  \big), %
	\end{equation}
	where
	\begin{equation}
	\label{eq:PhiPsi3d}
	\begin{aligned}
	&\bs \Phi_{m,n,l}^1(\bs{x}) = \big( \phi_{m}(x_1)\psi_{n+1}(x_2)\psi_{l+1}(x_3),\;\; 0, \;\;0 \big)^{\intercal}, \quad 
	(m,n,l)\in(\mathbb{N},\mathbb{N}^{+},\mathbb{N}^{+}), 
	\\
	&\bs \Phi_{m,n,l}^2(\bs{x}) = \big( 0,\;\;  \psi_{m+1}(x_1)\phi_{n}(x_2)\psi_{l+1}(x_3),\;\; 0 \big)^{\intercal},\quad 
	(m,n,l)\in(\mathbb{N}^{+},\mathbb{N},\mathbb{N}^{+}), 
	\\
	&\bs \Phi_{m,n,l}^3(\bs{x}) = \big( 0,\;\;  0,\;\; \psi_{m+1}(x_1)\psi_{n+1}(x_2)\phi_{l}(x_3) \big)^{\intercal},\quad 
	(m,n,l)\in(\mathbb{N}^{+},\mathbb{N}^{+},\mathbb{N}), 
	\\
	&\Psi_{m,n,l}(\bs{x}) = \psi_{m+1}(x_1)\psi_{n+1}(x_2)\psi_{l+1}(x_3), \qquad \qquad \;  
	(m,n,l)\in(\mathbb{N}^{+},\mathbb{N}^{+},\mathbb{N}^{+}).
	\end{aligned}
	\end{equation}
\end{prop}

Correspondingly, the approximation scheme for weak formulation \eqref{eq:MaxwWeak} in three dimensions reads: 
find $\bs{u}_N\in\bs{H}_{N,0}({\rm curl};{\Lambda}^3)$ and $p_N\in{H}_{N,0}^1(\Lambda^3)$ such that
\begin{equation}
\label{eq:MaxwDis3d}
\begin{aligned}
& (\nabla\times\bs{u}_N,\nabla\times\bs{v}) + {\kappa} (\bs{u}_N,\bs{v}) + (\nabla{p}_N,\bs{v})=(\bs{f},\bs{v}),\quad \forall \bs{v}\in \bs{H}_{N,0}({\rm curl};{\Lambda}^3),
\\
& (\bs{u}_N,\nabla{q}) = -(\rho,q), \quad \forall {q}\in {H}_{N,0}^1(\Lambda^3).
\end{aligned}
\end{equation}
Let us expand the numerical solution $\bs{u}_N$ and $p_N$ by 
\begin{align}
\bs{u}_N(\bs{x}) 
&= 
\sum_{m=0}^{N-1}\sum_{n,l=1}^{N-1} u_{mnl}^1\bs\Phi_{m,n,l}^1(\bs{x})
+ \sum_{n=0}^{N-1}\sum_{m,l=1}^{N-1} u_{mnl}^2\bs\Phi_{m,n,l}^2(\bs{x})
+ \sum_{l=0}^{N-1}\sum_{m,n=1}^{N-1} u_{mnl}^3\bs\Phi_{m,n,l}^3(\bs{x}),
\label{eq:uN3d}
\\
p_N(\bs{x})  &= \sum_{m,n,l=1}^{N-1} p_{mnl} \Psi_{m,n,l}(\bs{x}), 
\label{eq:pN3d}
\end{align}
and take the test functions $\bs{v}= \bs{\Phi}_{\hat{m},\hat{n},\hat{l}}^i(\bs{x})$ with $i=1,2,3$ and $q(\bs{x})=\Psi_{\hat{m},\hat{n},\hat{l}}(\bs{x})$ into scheme \eqref{eq:MaxwDis3d}. 
Denote the source and known tensors by
\begin{equation}\label{eq:MatFGHR3d}
\begin{aligned}
& \mathbf{F}=(f_{mnl}), \;\;  f_{mnl}=\big( \bs{f}, \bs\Phi_{m,n,l}^1 \big);    \quad \;
\mathbf{G}=(g_{mnl}),\;\; g_{mnl}=\big(\bs{f}, \bs\Phi_{m,n,l}^2 \big);  
\\
& \mathbf{H}=(h_{mnl}), \;\; h_{mnl}=\big(\bs{f}, \bs\Phi_{m,n,l}^3 \big);    \quad 
\mathbf{R}=(r_{mnl}),\;\; r_{mnl}=-\big( \rho,{\Psi}_{m,n,l} \big);  
\\
& \mathbf{U}=(u_{mnl}^1),\quad \mathbf{V}=(u_{mnl}^2),    \quad \mathbf{W}=(u_{mnl}^3),   \quad  \mathbf{P}=(p_{mnl}),
\end{aligned}
\end{equation}
and reformulate them into vector forms
\begin{equation}
\label{eq:upVec3d}
\begin{aligned}
& \bs{\vec{f}}={\rm vec}( \mathbf{F} ), \quad \bs{\vec{g}}={\rm vec}( \mathbf{G} ), \; \quad \bs{\vec{h}}={\rm vec}( \mathbf{H} ), \quad  \bs{\vec{r}}={\rm vec}( \mathbf{R}), 
\\
& \bs{\vec{u}}={\rm vec}( \mathbf{U} ), \quad  \bs{\vec{v}}={\rm vec}( \mathbf{V} ), \quad  \bs{\vec{w}}={\rm vec}( \mathbf{W} ),  \quad \bs{\vec{p}}={\rm vec}( \mathbf{P} ),
\end{aligned}
\end{equation}
one arrives at the following algebraic system
\begin{equation}
\label{eq:MaxwAlgSys3d}
\begin{pmatrix}
\mathbf{A}_{11}  &  \mathbf{A}_{12}  &  \mathbf{A}_{13}  &  \mathbf{A}_{14}
\\
\mathbf{A}_{21}  &  \mathbf{A}_{22}  &  \mathbf{A}_{23}  &  \mathbf{A}_{24}
\\  
\mathbf{A}_{31}  &  \mathbf{A}_{32}  &  \mathbf{A}_{33}  &  \mathbf{A}_{34}
\\
\mathbf{A}_{41}  &  \mathbf{A}_{42}  &  \mathbf{A}_{43}  &  \mathbf{A}_{44} 
\end{pmatrix}
\begin{pmatrix}
\bs{\vec{u}} \\
\bs{\vec{v}} \\
\bs{\vec{w}} \\
\bs{\vec{p}}
\end{pmatrix}
=\begin{pmatrix}
\bs{\vec{f}} \\
\bs{\vec{g}} \\
\bs{\vec{h}} \\
\bs{\vec{r}}
\end{pmatrix}.
\end{equation}
In the above, the components of the coefficient matrix are defined as 
\begin{align*}
\mathbf{A}_{11} &= \mathbf{I}_{N-1}\otimes\mathbf{M}\otimes\mathbf{I}_{N} + \mathbf{M}\otimes\mathbf{I}_{N-1}\otimes\mathbf{I}_{N} + {\kappa}\mathbf{M}\otimes\mathbf{M}\otimes\mathbf{I}_{N},
\\
\mathbf{A}_{22} &= \mathbf{I}_{N-1}\otimes\mathbf{I}_{N}\otimes\mathbf{M} + \mathbf{M}\otimes\mathbf{I}_{N}\otimes\mathbf{I}_{N-1} + {\kappa}\mathbf{M}\otimes\mathbf{I}_{N}\otimes\mathbf{M},
\\
\mathbf{A}_{33} &= \mathbf{I}_{N}\otimes\mathbf{I}_{N-1}\otimes\mathbf{M} + \mathbf{I}_{N}\otimes\mathbf{M}\otimes\mathbf{I}_{N-1} + {\kappa}\mathbf{I}_{N}\otimes\mathbf{M}\otimes\mathbf{M},
\\
\mathbf{A}_{12} &=\mathbf{A}_{21}^{\intercal} = - \mathbf{M}\otimes\mathbf{\ushort{E}}^{\intercal}\otimes\mathbf{\ushort{E}}, 
\quad 
\mathbf{A}_{13} = \mathbf{A}_{31}^{\intercal} = - \mathbf{\ushort{E}}^{\intercal}\otimes\mathbf{M}\otimes\mathbf{\ushort{E}}, 
\\
\mathbf{A}_{23} &= \mathbf{A}_{32}^{\intercal} = - \mathbf{\ushort{E}}^{\intercal}\otimes\mathbf{\ushort{E}}\otimes\mathbf{M},
\quad
\mathbf{A}_{14} = \mathbf{A}_{41}^{\intercal} = \mathbf{M}\otimes\mathbf{M}\otimes\mathbf{\ushort{E}}^{\intercal},
\\
\mathbf{A}_{24} &= \mathbf{A}_{42}^{\intercal} = \mathbf{M}\otimes\mathbf{\ushort{E}}^{\intercal}\otimes\mathbf{M}, \;\;
\quad
\mathbf{A}_{34} = \mathbf{A}_{43}^{\intercal} = \mathbf{\ushort{E}}^{\intercal}\otimes\mathbf{M}\otimes\mathbf{M},
\\
\mathbf{A}_{44} &= \mathbf{0}_{N-1}\otimes\mathbf{0}_{N-1}\otimes\mathbf{0}_{N-1}.
\end{align*}

\subsection{Numerical algorithm} 
\label{Sec3-2}
Due to the saddle-point nature of the mixed variational formulation,  the resultant linear system is fully-coupled and poses great numerical challenge. 
Following the same spirit, we decompose linear system \eqref{eq:MaxwAlgSys3d} into decoupled sub-problems by left multiplying with the permutation matrix 
\begin{equation}
\label{eq:MaxwPemMat3d}
\begin{pmatrix}
\mathbf{I}_{N-1}\otimes\mathbf{I}_{N-1}\otimes\bs{\ushort{e}}_{1}^{\intercal}  &  &  &  
\\
\mathbf{I}_{N-1}\otimes\mathbf{I}_{N-1}\otimes\mathbf{\ushort{E}} &  & &  
\\
&   \mathbf{I}_{N-1}\otimes\bs{\ushort{e}}_{1}^{\intercal}\otimes\mathbf{I}_{N-1} & &  
\\
&   \mathbf{I}_{N-1}\otimes\mathbf{\ushort{E}}\otimes\mathbf{I}_{N-1}    & & 
\\
&  & \bs{\ushort{e}}_{1}^{\intercal}\otimes\mathbf{I}_{N-1}\otimes\mathbf{I}_{N-1} & 
\\
&   & \mathbf{\ushort{E}}\otimes\mathbf{I}_{N-1}\otimes\mathbf{I}_{N-1}    & \\
&  &   &  \mathbf{I}_{N-1}\otimes\mathbf{I}_{N-1}\otimes\mathbf{I}_{N-1}
\end{pmatrix}.
\end{equation}
Denote
\begin{equation}
\label{eq:UFVec3d}
\begin{aligned}
& 
\bs{\vec{U}^{x}} = ( \mathbf{I}_{N-1}\otimes\mathbf{I}_{N-1}\otimes\bs{\ushort{e}}_{1}^{\intercal} ) \bs{\vec{{u}}}, \quad 
\bs{\vec{U}^{y}} = ( \mathbf{I}_{N-1}\otimes\bs{\ushort{e}}_{1}^{\intercal}\otimes\mathbf{I}_{N-1} ) \bs{\vec{{v}}}, \quad 
\bs{\vec{U}^{z}} = ( \bs{\ushort{e}}_{1}^{\intercal}\otimes\mathbf{I}_{N-1}\otimes\mathbf{I}_{N-1} ) \bs{\vec{{w}}}, \quad 
\\
& 
\bs{\vec{U}^{1}} = ( \mathbf{I}_{N-1}\otimes\mathbf{I}_{N-1}\otimes\mathbf{\ushort{E}} ) \bs{\vec{{u}}}, \quad 
\bs{\vec{U}^{2}} = ( \mathbf{I}_{N-1}\otimes\mathbf{\ushort{E}}\otimes\mathbf{I}_{N-1} ) \bs{\vec{{v}}}, \quad  
\bs{\vec{U}^{3}} = ( \mathbf{\ushort{E}}\otimes\mathbf{I}_{N-1}\otimes\mathbf{I}_{N-1} ) \bs{\vec{{w}}}, \quad 
\\
& 
\bs{\vec{F}^{x}} = ( \mathbf{I}_{N-1}\otimes\mathbf{I}_{N-1}\otimes\bs{\ushort{e}}_{1}^{\intercal} ) \bs{\vec{f}}, \quad 
\bs{\vec{F}^{y}} = ( \mathbf{I}_{N-1}\otimes\bs{\ushort{e}}_{1}^{\intercal}\otimes\mathbf{I}_{N-1} ) \bs{\vec{g}}, \quad 
\bs{\vec{F}^{z}} = ( \bs{\ushort{e}}_{1}^{\intercal}\otimes\mathbf{I}_{N-1}\otimes\mathbf{I}_{N-1} ) \bs{\vec{h}}, \quad 
\\
&
\bs{\vec{F}^{1}} = ( \mathbf{I}_{N-1}\otimes\mathbf{I}_{N-1}\otimes\mathbf{\ushort{E}} ) \bs{\vec{f}}, \quad 
\bs{\vec{F}^{2}} = ( \mathbf{I}_{N-1}\otimes\mathbf{\ushort{E}}\otimes\mathbf{I}_{N-1} ) \bs{\vec{g}}, \quad 
\bs{\vec{F}^{3}} = ( \mathbf{\ushort{E}}\otimes\mathbf{I}_{N-1}\otimes\mathbf{I}_{N-1} ) \bs{\vec{h}},
\end{aligned}
\end{equation}
one can obtain the sub-systems involving $\bs{\vec{U}^{1}}$-$\bs{\vec{U}^{3}}$ and $\bs{\vec{p}}$, and $\bs{\vec{U}^{x}}$-$\bs{\vec{U}^{y}}$, respectively:
\begin{equation}
\label{eq:3dSub1MatEq}
\begin{aligned}
& \mathbf{B}_{11} \bs{\vec{U}^{1}} + \mathbf{B}_{12} \bs{\vec{U}^{2}} + \mathbf{B}_{13} \bs{\vec{U}^{3}} + \mathbf{B}_{14} \bs{\vec{p}} = \bs{\vec{F}^{1}}, 
\\
& \mathbf{B}_{21} \bs{\vec{U}^{1}} + \mathbf{B}_{22} \bs{\vec{U}^{2}} + \mathbf{B}_{23} \bs{\vec{U}^{3}} + \mathbf{B}_{24} \bs{\vec{p}} = \bs{\vec{F}^{2}}, 
\\
& \mathbf{B}_{31} \bs{\vec{U}^{1}} + \mathbf{B}_{32} \bs{\vec{U}^{2}} + \mathbf{B}_{33}\bs{\vec{U}^{3}} + \mathbf{B}_{34} \bs{\vec{p}} = \bs{\vec{F}^{3}}, 
\\
& \mathbf{B}_{41} \bs{\vec{U}^{1}} + \mathbf{B}_{42} \bs{\vec{U}^{2}} + \mathbf{B}_{43}\bs{\vec{U}^{3}} + \mathbf{B}_{44} \bs{\vec{p}} = \bs{\vec{r}}, 
\end{aligned}
\end{equation}
and 
\begin{align}
\label{eq:3dSub2MatEq}
& \mathbf{B}_{xx} \bs{\vec{U}^{x}} = \bs{\vec{F}^{x}}, \;
\mathbf{B}_{yy} \bs{\vec{U}^{y}} = \bs{\vec{F}^{y}}, \;
\mathbf{B}_{zz} \bs{\vec{U}^{z}} = \bs{\vec{F}^{z}}.
\end{align}
Here, the components of the coefficient matrix are given by
\begin{align*}
\mathbf{B}_{11} &= \mathbf{I}_{N-1}\otimes\mathbf{M}\otimes\mathbf{I}_{N-1} + \mathbf{M}\otimes\mathbf{I}_{N-1}\otimes\mathbf{I}_{N-1} + {\kappa}\mathbf{M}\otimes\mathbf{M}\otimes\mathbf{I}_{N-1},
\\
\mathbf{B}_{22} &= \mathbf{I}_{N-1}\otimes\mathbf{I}_{N-1}\otimes\mathbf{M} + \mathbf{M}\otimes\mathbf{I}_{N-1}\otimes\mathbf{I}_{N-1} + {\kappa}\mathbf{M}\otimes\mathbf{I}_{N-1}\otimes\mathbf{M},
\\
\mathbf{B}_{33} &= \mathbf{I}_{N-1}\otimes\mathbf{I}_{N-1}\otimes\mathbf{M} + \mathbf{I}_{N-1}\otimes\mathbf{M}\otimes\mathbf{I}_{N-1} + {\kappa}\mathbf{I}_{N-1}\otimes\mathbf{M}\otimes\mathbf{M},
\\
\mathbf{B}_{44} &= \mathbf{0}_{N-1}\otimes\mathbf{0}_{N-1}\otimes\mathbf{0}_{N-1},
\\
\mathbf{B}_{12} &=\mathbf{B}_{21}^{\intercal} = - \mathbf{M}\otimes\mathbf{I}_{N-1}\otimes\mathbf{I}_{N-1},  \quad
\mathbf{B}_{13} = \mathbf{B}_{31}^{\intercal} = - \mathbf{I}_{N-1}\otimes\mathbf{M}\otimes\mathbf{I}_{N-1}, 
\\
\mathbf{B}_{23} &= \mathbf{B}_{32}^{\intercal} = - \mathbf{I}_{N-1}\otimes\mathbf{I}_{N-1}\otimes\mathbf{M}, \quad 
\mathbf{B}_{14} = \mathbf{B}_{41}^{\intercal} = \mathbf{M}\otimes\mathbf{M}\otimes\mathbf{I}_{N-1},
\\
\mathbf{B}_{24} &= \mathbf{B}_{42}^{\intercal} = \mathbf{M}\otimes\mathbf{I}_{N-1}\otimes\mathbf{M}, \quad \qquad 
\mathbf{B}_{34} = \mathbf{B}_{43}^{\intercal} = \mathbf{I}_{N-1}\otimes\mathbf{M}\otimes\mathbf{M},
\\
\mathbf{B}_{xx} &= \mathbf{B}_{yy} = \mathbf{B}_{zz} = \mathbf{I}_{N-1}\otimes\mathbf{M} + \mathbf{M}\otimes\mathbf{I}_{N-1} + {\kappa}\mathbf{M}\otimes\mathbf{M}.
\end{align*}

Next, we can follow the same idea presented in Subsection \ref{Sec2-2} to arrive at  a highly efficient, semi-analytic, matrix-free approach to solve sub-systems \eqref{eq:3dSub1MatEq}-\eqref{eq:3dSub2MatEq}. 
The derivation of the algorithm follows the same manner with the 2D case but much more complex and tedious, we postpone the details in Appendix \ref{Appx1-Algo2}
and give the solution procedure for system \eqref{eq:MaxwAlgSys3d} in Algorithm \ref{Algo2}.

\begin{algorithm}[ht!]
	\caption{Matrix-free, semi-analytic solution algorithm for  \eqref{eq:MaxwAlgSys3d} with Gauss's law constraint}
	\label{Algo2}
	\begin{algorithmic}
		\item  
		[{\bf Input:}] 
		The polynomial order $N$, parameter $\kappa$,  $\mathbf{F}$, $\mathbf{G}$, $\mathbf{H}$ and $\mathbf{R}$ in \eqref{eq:MatFGHR3d},  $\ushort{\bs{e}}_1$ 
		and $\ushort{\mathbf{E}}$ in \eqref{eq:MatINE}, $\mathbf{D}$, $\mathbf{Q}$ in \eqref{eq:MatMIdiag}.
		\vspace{3pt} 
		
		\item 
		[{\bf Output:}] 
		The solutions $\bs{\vec{u}}$, $\bs{\vec{v}}$, $\bs{\vec{w}}$ and $\bs{\vec{p}}$.
		\vspace{2pt} 
		\item[1:] 
		Compute $\bs{\vec{F}^{x}}$-$\bs{\vec{F}^{z}}$ and $\bs{\vec{F}^{1}}$-$\bs{\vec{F}^{3}}$ by \eqref{eq:UFVec3d}.
	
		\item [2:] 
		Compute $\mathbf{\hat{F}}^1$, $\mathbf{\hat{F}}^2$, $\mathbf{\hat{F}}^3$ and $\mathbf{\hat{R}}$ by
		\begin{equation}
		\label{eq:MatFsRhat}
		\mathbf{\hat{F}}^{s}={\rm ivec}(\bs{\vec{\hat{F}}^{s}}), \quad 
		\bs{\vec{\hat{F}}^{s}} = ( \mathbf{Q}^{\intercal}\otimes\mathbf{Q}^{\intercal}\otimes\mathbf{Q}^{\intercal} ) \bs{\vec{F}^{s}}, \;\; s=1,2,3; \quad 
		\mathbf{\hat{R}}={\rm ivec}(\bs{\vec{\hat{r}}}), \;\;  \bs{\vec{\hat{r}}} = ( \mathbf{Q}^{\intercal}\otimes\mathbf{Q}^{\intercal}\otimes\mathbf{Q}^{\intercal} ) \bs{\vec{r}}. 
		\end{equation}
	
		\item [3:] 
		Compute $\bs{\vec{p}}$ by 
		\vspace{-6pt}
		\begin{equation}
		\vspace{-4pt}
		\label{eq:PSol3d}
		\bs{\vec{p}} = (\mathbf{Q}\otimes\mathbf{Q}\otimes\mathbf{Q}) \bs{\vec{\hat{p}}}, \quad
		\bs{\vec{\hat{p}}}={\rm vec}(\mathbf{\hat{P}}), \quad \mathbf{\hat{P}}=(\hat{P}_{ijk}), \quad
		\hat{P}_{ijk} = \frac{ \hat{F}_{ijk}^1 +  \hat{F}_{ijk}^2 +  \hat{F}_{ijk}^3 -  {\kappa}\hat{R}_{ijk} }{d_{j}d_{k}+ d_{i}d_{k} + d_{i}d_{j}}.
		\end{equation}
		
		\item [4:] 
		Compute $\bs{\vec{U}^{1}}$, $\bs{\vec{U}^{2}}$ and $\bs{\vec{U}^{3}}$ by 
		\begin{align}
		& 
		\bs{\vec{U}^{1}} = ( \mathbf{Q}\otimes\mathbf{Q}\otimes\mathbf{Q} ) \bs{\vec{\hat{U}}^{1}}, \quad
		\bs{\vec{\hat{U}}^{1}} = {\rm vec}(\mathbf{\hat{U}}^{1}), \quad \mathbf{\hat{U}^1}=(\hat{U}_{ijk}^{1} ), \quad
		\hat{U}_{ijk}^{1} = \frac{\hat{R}_{ijk} + d_i \hat{F}_{ijk}^1 - d_i d_j d_k \hat{P}_{ijk}}{d_jd_k + d_id_k + d_id_j + {\kappa} d_id_jd_k},
		\label{eq:U1Sol3d}
		\\
		& 
		\bs{\vec{U}^{2}} = ( \mathbf{Q}\otimes\mathbf{Q}\otimes\mathbf{Q} ) \bs{\vec{\hat{U}}^{2}}, \quad
		\bs{\vec{\hat{U}}^{2}} = {\rm vec}(\mathbf{\hat{U}}^{2}), \quad \mathbf{\hat{U}^2}=(\hat{U}_{ijk}^{2} ), \quad
		\hat{U}_{ijk}^{2} = \frac{\hat{R}_{ijk} + d_j \hat{F}_{ijk}^2 - d_i d_j d_k \hat{P}_{ijk}}{d_jd_k + d_id_k + d_id_j + {\kappa} d_id_jd_k},
		\label{eq:U2Sol3d}
		\\
		&
		\bs{\vec{U}^{3}} = ( \mathbf{Q}\otimes\mathbf{Q}\otimes\mathbf{Q} ) \bs{\vec{\hat{U}}^{3}}, \quad
		\bs{\vec{\hat{U}}^{3}} = {\rm vec}(\mathbf{\hat{U}}^{3}), \quad \mathbf{\hat{U}^3}=(\hat{U}_{ijk}^{3} ), \quad
		\hat{U}_{ijk}^{3} = \frac{\hat{R}_{ijk} + d_k \hat{F}_{ijk}^3 - d_i d_j d_k \hat{P}_{ijk}}{d_jd_k + d_id_k + d_id_j + {\kappa} d_id_jd_k}.
		\label{eq:U3Sol3d}
		\end{align}
		
		\item [5:] 
		Compute $\mathbf{\hat{F}}^{x}$, $\mathbf{\hat{F}}^{y}$ and $\mathbf{\hat{F}}^{z}$ by
		\begin{equation}
		\label{eq:MatFchihat}
		\mathbf{\hat{F}}^{\chi} = {\rm ivec}(\bs{\vec{\hat{F}}^{\chi}}), \quad 
		\bs{\vec{\hat{F}}^{\chi}} = ( \mathbf{Q}^{\intercal}\otimes\mathbf{Q}^{\intercal} ) \bs{\vec{F}^{\chi}},\quad  
		\chi = x,y,z.
		\end{equation}
		
		\item [6:]
		Compute $\bs{\vec{U}^{x}}$, $\bs{\vec{U}^{y}}$ and $\bs{\vec{U}^{z}}$ by 
		\begin{align}
		& 
		\bs{\vec{U}^{x}} = {\rm vec}(\mathbf{U}^{x}), \quad 
		\mathbf{U}^{x} = \mathbf{Q}\mathbf{\hat{U}}^{x}\mathbf{Q}^{\intercal}, \quad  
		\mathbf{\hat{U}}^{x} = (\hat{U}_{jk}^{x}), \quad
		\hat{U}_{jk}^{x} = \frac{\hat{F}_{jk}^{x}}{ d_{j} + d_{k} + {\kappa}d_{j}d_{k}}, \quad
		\label{eq:UxSol3d}
		\\
		& 
		\bs{\vec{U}^{y}} = {\rm vec}(\mathbf{U}^{y}), \quad 
		\mathbf{U}^{y} = \mathbf{Q}\mathbf{\hat{U}}^{y}\mathbf{Q}^{\intercal}, \quad  
		\mathbf{\hat{U}}^{y} = (\hat{U}_{ik}^{y}), \quad
		\hat{U}_{ik}^{y} = \frac{\hat{F}_{ik}^{y}}{ d_{i} + d_{k} + {\kappa}d_{i}d_{k}}, \quad
		\label{eq:UySol3d}
		\\
		&
		\bs{\vec{U}^{z}} = {\rm vec}(\mathbf{U}^{z}), \quad 
		\mathbf{U}^{z} = \mathbf{Q}\mathbf{\hat{U}}^{z}\mathbf{Q}^{\intercal}, \quad   
		\mathbf{\hat{U}}^{z} = (\hat{U}_{ij}^{z}), \quad
		\hat{U}_{ij}^{z} = \frac{\hat{F}_{ij}^{z}}{ d_{i} + d_{j} + {\kappa}d_{i}d_{j}}.
		\label{eq:UzSol3d}
		\end{align}
		
		\item [7:] 
		Return $\bs{\vec{u}}$, $\bs{\vec{v}}$ and $\bs{\vec{w}}$ by 
		\begin{align}
		& \bs{\vec{u}} = ( \mathbf{I}_{N-1}\otimes\mathbf{I}_{N-1}\otimes\mathbf{\ushort{E}}^{\intercal} ) \bs{\vec{U}^{1}} 
		+ ( \mathbf{I}_{N-1}\otimes\mathbf{I}_{N-1}\otimes\bs{\ushort{e}}_{1} ) \bs{\vec{U}^{x}},
		\label{eq:uVec3d}
		\\
		& \bs{\vec{v}} = ( \mathbf{I}_{N-1}\otimes\mathbf{\ushort{E}}^{\intercal}\otimes\mathbf{I}_{N-1} ) \bs{\vec{U}^{2}}
		+ ( \mathbf{I}_{N-1}\otimes\bs{\ushort{e}}_{1}\otimes\mathbf{I}_{N-1} ) \bs{\vec{U}^{y}},
		\label{eq:vVec3d}
		\\
		& \bs{\vec{w}} =  ( \mathbf{\ushort{E}}^{\intercal}\otimes\mathbf{I}_{N-1}\otimes\mathbf{I}_{N-1} ) \bs{\vec{U}^{3}}
		+ ( \bs{\ushort{e}}_{1}\otimes\mathbf{I}_{N-1}\otimes\mathbf{I}_{N-1} ) \bs{\vec{U}^{z}}.
		\label{eq:wVec3d}
		\end{align}
		
	\end{algorithmic}
\end{algorithm}

\begin{rem}
\label{Remk3-Algo2}
The bottleneck of the computational cost of Algorithm \ref{Algo2} resides in steps 3-4, where the Kronecker product of the form 
$\bs{\vec{v}}=(\mathbf{Q}\otimes\mathbf{Q}\otimes\mathbf{Q})\bs{\vec{u}}$ has been invoked four times.  
With a fully assembled manner, the computational complexity for these Kronecker products requires $\mathcal{O}(N^6)$ operations. 
However, here we utilize the tensor structure and adopt the partial assemble technique to calculate the Kronecker product by
\begin{equation*}
\bs{\vec{v}}={\rm vec}(\mathbf{V}),\;\; \mathbf{V}=(V_{ijk}), \;\; 
V_{ijk}=\sum_{\hat{i},\hat{j},\hat{k}} Q_{i\hat{i}} Q_{j\hat{j}} Q_{k\hat{k}} U_{\hat{i}\hat{j}\hat{k}}, \;\; \mathbf{U}=(U_{ijk}), \;\; \mathbf{U}={\rm ivec}(\bs{\vec{u}}). 
\end{equation*} 
Moreover, when compared with direct solution method  based on LU factorization, 
the total computational complexity for solving system \eqref{eq:MaxwAlgSys3d} reduces from $\mathcal{O}(N^9)$ to $\mathcal{O}(N^4)$, 
and can further be reduced to $\mathcal{O}(N^{1+\log_{2}7})$ with the help of Strassen's matrix multiplication algorithm \cite{Strassen-1969}, 
as verified by Example \ref{Exmp5-4} in Section \ref{Sec5}.
\end{rem}

\subsection{Double-curl eigenvalue problems in 3D}
\label{Sec3-3}
Likewise the 2D eigenvalue problem, the approximation scheme for the weak formulation \eqref{eq:EigWeak} in three dimensions takes the form: 
find ${\lambda}^{N}\in\mathbb{R}$ and $\bs{u}_{N}\in\bs{H}_{N,0}({\rm curl};{\Lambda}^3)$ with $\bs{u}_N\neq\bs{0}$ such that
\begin{equation}
\label{eq:EigDis3d}
\begin{aligned}
& (\nabla\times\bs{u}_N,\nabla\times\bs{v}) = {\lambda}^N (\bs{u}_N,\bs{v}) ,\quad \forall \bs{v}\in\bs{H}_{N,0}({\rm curl};{\Lambda}^3).
\end{aligned}
\end{equation}
With a slight abuse of notation, we arrive at the generalised matrix eigenvalue problem as follows:
\begin{equation}
\label{eq:EigAlgSys3d}
\begin{pmatrix}
\mathbf{{A}}_{11}  &  \mathbf{{A}}_{12}  &  \mathbf{{A}}_{13}  
\\
\mathbf{{A}}_{21}  &  \mathbf{{A}}_{22}  &  \mathbf{{A}}_{23}  
\\  
\mathbf{{A}}_{31}  &  \mathbf{{A}}_{32}  &  \mathbf{{A}}_{33}  
\end{pmatrix}
\begin{pmatrix}
\bs{\vec{u}} 
\\
\bs{\vec{v}} 
\\
\bs{\vec{w}} 
\end{pmatrix}
={\lambda}^{N}
\begin{pmatrix}
\mathbf{{A}_1} & & 
\\
& \mathbf{{A}_2}  & 
\\
& & \mathbf{{A}_3}   
\end{pmatrix}
\begin{pmatrix}
\bs{\vec{u}} 
\\
\bs{\vec{v}} 
\\
\bs{\vec{w}} 
\end{pmatrix},
\end{equation}
where
\begin{align*}
& \mathbf{{A}}_{11} = \mathbf{I}_{N-1}\otimes\mathbf{M}\otimes\mathbf{I}_{N} + \mathbf{M}\otimes\mathbf{I}_{N-1}\otimes\mathbf{I}_{N},
\quad
\mathbf{{A}}_{22} = \mathbf{I}_{N-1}\otimes\mathbf{I}_{N}\otimes\mathbf{M} + \mathbf{M}\otimes\mathbf{I}_{N}\otimes\mathbf{I}_{N-1},
\\
&
\mathbf{{A}}_{33} = \mathbf{I}_{N}\otimes\mathbf{I}_{N-1}\otimes\mathbf{M} + \mathbf{I}_{N}\otimes\mathbf{M}\otimes\mathbf{I}_{N-1}, 
\quad
\mathbf{{A}}_{12} =\mathbf{{A}}_{21}^{\intercal} = - \mathbf{M}\otimes\mathbf{\ushort{E}}^{\intercal}\otimes\mathbf{\ushort{E}}, 
\\
&
\mathbf{{A}}_{13} = \mathbf{{A}}_{31}^{\intercal} = - \mathbf{\ushort{E}}^{\intercal}\otimes\mathbf{M}\otimes\mathbf{\ushort{E}}, 
\quad \qquad \qquad \qquad
\mathbf{{A}}_{23} = \mathbf{{A}}_{32}^{\intercal} = - \mathbf{\ushort{E}}^{\intercal}\otimes\mathbf{\ushort{E}}\otimes\mathbf{M},
\\
&
\mathbf{{A}}_{1} = \mathbf{M}\otimes\mathbf{M}\otimes\mathbf{I}_{N},
\quad
\mathbf{{A}}_{2} = \mathbf{M}\otimes\mathbf{I}_{N}\otimes\mathbf{M},
\quad
\mathbf{{A}}_{3} = \mathbf{I}_{N}\otimes\mathbf{M}\otimes\mathbf{M}.
\end{align*}

The solution algorithm follows the same derivation procedure. For the sake of conciseness, we postpone the detailed proof in Appendix \ref{Appx2-Thm2} 
and present the matrix-free, semi-analytic expression of the egien-solution for system \eqref{eq:EigDis3d} in Theorem \ref{Thm2-LamU3dSol}.
\begin{thm}\label{Thm2-LamU3dSol}
	The exact eigen-solutions $({\lambda}^{N},\bs{\vec{u}},\bs{\vec{v}},\bs{\vec{w}})$ of system \eqref{eq:EigAlgSys3d} shall be divided into two categories. 
	\begin{itemize}
	\item The eigen-solution $({\lambda}_{ijk}^{N},\mathbf{U}_{ijk}^{1},\mathbf{U}_{ijk}^{2},\mathbf{U}_{ijk}^{3})$  with non-zero eigenvalues includes:     
	\begin{itemize}
		\item[---] 
		For $1\leq{i,j,k}\leq{N-1}$, 
		the interior $(i,j,k)$-th eigenvalue ${\lambda}_{ijk}^{N}$ with a multiplicity of $2$:
		\begin{equation}
		\label{eq:Lam3dSolIn}
		{\lambda}_{ijk}^{N} = \frac{1}{d_i} +  \frac{1}{d_j} +  \frac{1}{d_k}, 
		\end{equation}
		and the associated two independent eigenfunctions $\mathbf{U}_{ijk}^{s} = {\rm ivec}(\bs{\vec{U}_{ijk}^s})$, $s=1,2,3$ with
		\begin{align}
			\label{eq:U3dSolIn1}
			\bs{\vec{U}_{ijk}^1} = (\mathbf{Q}\otimes\mathbf{Q}\otimes\mathbf{Q}) \bs{\vec{I}}_{ijk}, \quad 
			\bs{\vec{U}_{ijk}^2} = \bs{\vec{0}}, \quad 
			\bs{\vec{U}_{ijk}^3} = -(\mathbf{Q}\otimes\mathbf{Q}\otimes\mathbf{Q}) \big( {d_k}/{d_i}\bs{\vec{I}}_{ijk} \big), 
		\end{align}
		and
		\begin{align}
			\label{eq:U3dSolIn2}
			\bs{\vec{U}_{ijk}^1} = \bs{\vec{0}}, \quad 
			\bs{\vec{U}_{ijk}^2} = (\mathbf{Q}\otimes\mathbf{Q}\otimes\mathbf{Q})  \bs{\vec{I}}_{ijk}, \quad 
			\bs{\vec{U}_{ijk}^3} = -(\mathbf{Q}\otimes\mathbf{Q}\otimes\mathbf{Q})  \big( {d_k}/{d_j}\bs{\vec{I}}_{ijk} \big), 
		\end{align}
		where $\bs{\vec{I}}_{ijk} = {\rm vec}(\mathbf{I}_{ijk})$, the $(N-1)\times(N-1)\times(N-1)$ matrix $\mathbf{I}_{ijk}$ 
		whose $(i,j,k)$-th element is $1$, while the rest elements are $0$.
		Then $({\lambda}^{N},\bs{\vec{u}},\bs{\vec{v}},\bs{\vec{w}})$ are given by ${\lambda}^{N} = {\lambda}_{ijk}^{N}$ and 
		\begin{equation}
		\label{eq:LamU3dSolIn}
		\begin{aligned}
		&	\bs{\vec{u}} = ( \mathbf{I}_{N-1}\otimes\mathbf{I}_{N-1}\otimes\mathbf{\ushort{E}}^{\intercal} ) \bs{\vec{U}_{ijk}^{1}} 
		+ ( \mathbf{I}_{N-1}\otimes\mathbf{I}_{N-1}\otimes\bs{\ushort{e}}_{1} ) \bs{\vec{U}^{x}}, \quad \; 
		\bs{\vec{U}^{x}}={\rm vec}(\mathbf{0}_{N-1}),
		\\
		&	\bs{\vec{v}} = ( \mathbf{I}_{N-1}\otimes\mathbf{\ushort{E}}^{\intercal}\otimes\mathbf{I}_{N-1} ) \bs{\vec{U}_{ijk}^{2}}
		+ ( \mathbf{I}_{N-1}\otimes\bs{\ushort{e}}_{1}\otimes\mathbf{I}_{N-1} ) \bs{\vec{U}^{y}}, \quad \; 
		\bs{\vec{U}^{y}}={\rm vec}(\mathbf{0}_{N-1}),
		\\
		& \bs{\vec{w}} =  ( \mathbf{\ushort{E}}^{\intercal}\otimes\mathbf{I}_{N-1}\otimes\mathbf{I}_{N-1} ) \bs{\vec{U}_{ijk}^{3}}
		+ ( \bs{\ushort{e}}_{1}\otimes\mathbf{I}_{N-1}\otimes\mathbf{I}_{N-1} ) \bs{\vec{U}^{z}}, \quad \; 
		\bs{\vec{U}^{z}}={\rm vec}(\mathbf{0}_{N-1}).
		\end{aligned}
		\end{equation}
		
		\item [---]
		For $i=0$ and $1\leq{j,k}\leq{N-1}$, the $(j,k)$-th eigen-solution $({\lambda}_{0jk}^{N},\bs{\vec{U}_{0jk}})$ on the face takes the form
		\begin{equation}
		\label{eq:LamUyz3dFa1}
		{\lambda}_{0jk}^{N} = \frac{1}{d_j} + \frac{1}{d_k},
		\quad \; 
		\bs{\vec{U}_{0jk}} = {\rm vec}(\mathbf{U}_{0jk}), 
		\quad \;
		\mathbf{U}_{0jk} = \mathbf{Q}\mathbf{I}_{jk}\mathbf{Q}^{\intercal}.
		\end{equation}
		Then $({\lambda}^{N},\bs{\vec{u}},\bs{\vec{v}},\bs{\vec{w}})$ are given by ${\lambda}^{N} = {\lambda}_{0jk}^{N}$ and 
		\begin{equation}
		\label{eq:LamU3dSolFa1}
		\begin{aligned}
		&	\bs{\vec{u}} = ( \mathbf{I}_{N-1}\otimes\mathbf{I}_{N-1}\otimes\mathbf{\ushort{E}}^{\intercal} ) \bs{\vec{U}^{1}} 
		+ ( \mathbf{I}_{N-1}\otimes\mathbf{I}_{N-1}\otimes\bs{\ushort{e}}_{1} ) \bs{\vec{U}_{0jk}}, \quad \; 
		\bs{\vec{U}^{1}} = {\rm vec}(\mathbf{0}),  
		\\
		& 	\bs{\vec{v}} = ( \mathbf{I}_{N-1}\otimes\mathbf{\ushort{E}}^{\intercal}\otimes\mathbf{I}_{N-1} ) \bs{\vec{U}^{2}}
		+ ( \mathbf{I}_{N-1}\otimes\bs{\ushort{e}}_{1}\otimes\mathbf{I}_{N-1} ) \bs{\vec{0}}, \quad \; \quad \; \;
		\bs{\vec{U}^{2}} = {\rm vec}(\mathbf{0}),
		\\
		& \bs{\vec{w}} =  ( \mathbf{\ushort{E}}^{\intercal}\otimes\mathbf{I}_{N-1}\otimes\mathbf{I}_{N-1} ) \bs{\vec{U}^{3}}
		+ ( \bs{\ushort{e}}_{1}\otimes\mathbf{I}_{N-1}\otimes\mathbf{I}_{N-1} ) \bs{\vec{0}}, \quad \; \quad \;
		\bs{\vec{U}^{3}}={\rm vec}(\mathbf{0}),
		\end{aligned}
		\end{equation}
		where $\mathbf{0}$ is short for the zero matrix in $\mathbb{R}^{(N-1)\times{(N-1)}\times{(N-1)}}$.
		
		\item [---]
		For $j=0$ and $1\leq{i,k}\leq{N-1}$, the $(i,k)$-th eigen-solution $({\lambda}_{i0k}^{N},\bs{\vec{U}_{i0k}})$ on the face takes the form
		\begin{equation}
		\label{eq:LamUxz3dFa2}
		{\lambda}_{i0k}^{N} =  \frac{1}{d_i} + \frac{1}{d_k},
		\quad \; 
		\bs{\vec{U}_{i0k}} = {\rm vec}(\mathbf{U}_{i0k}),  
		\quad \;
		\mathbf{U}_{i0k} = \mathbf{Q}\mathbf{I}_{ik}\mathbf{Q}^{\intercal}. 
		\end{equation}
		Then $({\lambda}^{N},\bs{\vec{u}},\bs{\vec{v}},\bs{\vec{w}})$ are given by ${\lambda}^{N} = {\lambda}_{i0k}^{N}$ and 
		\begin{equation}
		\label{eq:LamU3dSolFa2}
		\begin{aligned}
		&	\bs{\vec{u}} = ( \mathbf{I}_{N-1}\otimes\mathbf{I}_{N-1}\otimes\mathbf{\ushort{E}}^{\intercal} ) \bs{\vec{U}^{1}} 
		+ ( \mathbf{I}_{N-1}\otimes\mathbf{I}_{N-1}\otimes\bs{\ushort{e}}_{1} ) \bs{\vec{0}}, \quad \;  \quad \;
		\bs{\vec{U}^{1}}={\rm vec}(\mathbf{0}),
		\\
		& 	\bs{\vec{v}} = ( \mathbf{I}_{N-1}\otimes\mathbf{\ushort{E}}^{\intercal}\otimes\mathbf{I}_{N-1} ) \bs{\vec{U}^{2}}
		+ ( \mathbf{I}_{N-1}\otimes\bs{\ushort{e}}_{1}\otimes\mathbf{I}_{N-1} ) \bs{\vec{U}_{i0k}}, \quad \; 
		\bs{\vec{U}^{2}}={\rm vec}(\mathbf{0}),
		\\
		& \bs{\vec{w}} =  ( \mathbf{\ushort{E}}^{\intercal}\otimes\mathbf{I}_{N-1}\otimes\mathbf{I}_{N-1} ) \bs{\vec{U}^{3}}
		+ ( \bs{\ushort{e}}_{1}\otimes\mathbf{I}_{N-1}\otimes\mathbf{I}_{N-1} ) \bs{\vec{0}}, \quad \;  \quad \;
		\bs{\vec{U}^{3}}={\rm vec}(\mathbf{0}).
		\end{aligned}
		\end{equation}
		
		\item [---]
		For $k=0$ and $1\leq{i,j}\leq{N-1}$, the $(i,j)$-th eigen-solution $({\lambda}_{ij0}^{N},\bs{\vec{U}_{ij0}})$ on the face takes the form
		\begin{equation}
		\label{eq:LamUxy3dFa3}
		{\lambda}_{ij0}^{N} = \frac{1}{d_i} + \frac{1}{d_j}, 
		\quad \; 
		\bs{\vec{U}_{ij0}} = {\rm vec}(\mathbf{U}_{ij0}), 
		\quad \;
		\mathbf{U}_{ij0} = \mathbf{Q}\mathbf{I}_{ij}\mathbf{Q}^{\intercal}. 
		\end{equation}
		Then $({\lambda}^{N},\bs{\vec{u}},\bs{\vec{v}},\bs{\vec{w}})$ are given by $\lambda^{N} = \lambda_{ij0}^{N}$ and 
		\begin{equation}
		\label{eq:LamU3dSolFa3}
		\begin{aligned}
		&	\bs{\vec{u}} = ( \mathbf{I}_{N-1}\otimes\mathbf{I}_{N-1}\otimes\mathbf{\ushort{E}}^{\intercal} ) \bs{\vec{U}^{1}} 
		+ ( \mathbf{I}_{N-1}\otimes\mathbf{I}_{N-1}\otimes\bs{\ushort{e}}_{1} ) \bs{\vec{0}}, \quad \;  \quad \;
		\bs{\vec{U}^{1}} = {\rm vec}(\mathbf{0}),
		\\
		& 	\bs{\vec{v}} = ( \mathbf{I}_{N-1}\otimes\mathbf{\ushort{E}}^{\intercal}\otimes\mathbf{I}_{N-1} ) \bs{\vec{U}^{2}}
		+ ( \mathbf{I}_{N-1}\otimes\bs{\ushort{e}}_{1}\otimes\mathbf{I}_{N-1} ) \bs{\vec{0}}, \quad \; \quad \;
		\bs{\vec{U}^{2}} = {\rm vec}(\mathbf{0}),
		\\
		& \bs{\vec{w}} =  ( \mathbf{\ushort{E}}^{\intercal}\otimes\mathbf{I}_{N-1}\otimes\mathbf{I}_{N-1} ) \bs{\vec{U}^{3}}
		+ ( \bs{\ushort{e}}_{1}\otimes\mathbf{I}_{N-1}\otimes\mathbf{I}_{N-1} )  \bs{\vec{U}_{ij0}}, \quad  
		\bs{\vec{U}^{3}} = {\rm vec}(\mathbf{0}).
		\end{aligned}
		\end{equation}
		
		\end{itemize}
		
		\item 
		The eigen-solution $({\lambda}_{ijk}^{N},\mathbf{U}_{ijk}^{1},\mathbf{U}_{ijk}^{2},\mathbf{U}_{ijk}^{3})$ 
		associated with zero eigenvalue for $1\leq{i,j,k}\leq{N-1}$ takes the form 
		\begin{align}
		\label{eq:Lam0Uijk3d}
		& {\lambda}_{ijk}^{N} = 0, \quad \;
		\mathbf{U}_{ijk}^{s} = {\rm ivec}(\bs{\vec{U}_{ijk}^s}), \quad \;
		\bs{\vec{U}_{ijk}^{s}} = (\mathbf{Q}\otimes\mathbf{Q}\otimes\mathbf{Q})\bs{\vec{I}}_{ijk}, \quad \; s=1,2,3. 
		\end{align}
		Then $({\lambda}^{N},\bs{\vec{u}},\bs{\vec{v}},\bs{\vec{w}})$ are given by ${\lambda}^{N} = 0$ and 
		\begin{equation}
		\label{eq:Lam0U3dSol}
		\begin{aligned}
		&	\bs{\vec{u}} = ( \mathbf{I}_{N-1}\otimes\mathbf{I}_{N-1}\otimes\mathbf{\ushort{E}}^{\intercal} ) \bs{\vec{U}_{ijk}^{1}} 
		+ ( \mathbf{I}_{N-1}\otimes\mathbf{I}_{N-1}\otimes\bs{\ushort{e}}_{1} ) \bs{\vec{U}^{x}}, \quad \; 
		\bs{\vec{U}^{x}}={\rm vec}(\mathbf{0}_{N-1}),
		\\
		&	\bs{\vec{v}} = ( \mathbf{I}_{N-1}\otimes\mathbf{\ushort{E}}^{\intercal}\otimes\mathbf{I}_{N-1} ) \bs{\vec{U}_{ijk}^{2}}
		+ ( \mathbf{I}_{N-1}\otimes\bs{\ushort{e}}_{1}\otimes\mathbf{I}_{N-1} ) \bs{\vec{U}^{y}}, \quad \; 
		\bs{\vec{U}^{y}}={\rm vec}(\mathbf{0}_{N-1}),
		\\
		& \bs{\vec{w}} =  ( \mathbf{\ushort{E}}^{\intercal}\otimes\mathbf{I}_{N-1}\otimes\mathbf{I}_{N-1} ) \bs{\vec{U}_{ijk}^{3}}
		+ ( \bs{\ushort{e}}_{1}\otimes\mathbf{I}_{N-1}\otimes\mathbf{I}_{N-1} ) \bs{\vec{U}^{z}}, \quad \;
		\bs{\vec{U}^{z}}={\rm vec}(\mathbf{0}_{N-1}).
		\end{aligned}
		\end{equation}
	\end{itemize}
\end{thm}

\begin{rem}
	\label{Remk4-Thm2}
	Similar to Remark \ref{Remk2-Thm1}, the above proposed algorithm also strictly obeys the Helmholtz-Hodge decomposition, 
	i.e., for any $\bs{u}\in \bs{H}_{0}({\rm curl};\Omega)$, there exists a unique scalar potential $\phi\in H_0^1(\Omega)$ 
	and a unique vector potential $\bs{\psi}\in\bs{H}({\rm curl};\Omega)\cap\bs{H}_0({\rm div};\Omega)$ with $\nabla\cdot\bs\psi=0$ such that
	$\bs{u}=\nabla\phi + \nabla\times\bs{\psi}$.
	This conclusion can be observed as follows: 
	\begin{itemize}
		\item
		When ${\lambda}^{N}\neq{0}$, 
		following the same derivation procedure of sub-systems \eqref{eq:3dSub1MatEq}-\eqref{eq:3dSub2MatEq},
		we decompose linear system \eqref{eq:EigAlgSys3d} into the decoupled sub-systems by left multiplying the first 6-by-3 block submatrix of the permutation matrix 
		given in \eqref{eq:MaxwPemMat3d}, 
		i.e., 
		\begin{equation*}
		\begin{pmatrix}
		\mathbf{I}_{N-1}\otimes\mathbf{I}_{N-1}\otimes\bs{\ushort{e}}_{1}^{\intercal}  &  &  
		\\
		\mathbf{I}_{N-1}\otimes\mathbf{I}_{N-1}\otimes\mathbf{\ushort{E}} &  & 
		\\
		&   \mathbf{I}_{N-1}\otimes\bs{\ushort{e}}_{1}^{\intercal}\otimes\mathbf{I}_{N-1} & 
		\\
		&   \mathbf{I}_{N-1}\otimes\mathbf{\ushort{E}}\otimes\mathbf{I}_{N-1}    & 
		\\
		&  & \bs{\ushort{e}}_{1}^{\intercal}\otimes\mathbf{I}_{N-1}\otimes\mathbf{I}_{N-1} 
		\\
		&   & \mathbf{\ushort{E}}\otimes\mathbf{I}_{N-1}\otimes\mathbf{I}_{N-1}   
		\end{pmatrix},
		\end{equation*}
		and then obtain the constraint 
		\[
			( \mathbf{M}\otimes\mathbf{M}\otimes\mathbf{I}_{N-1} ) \bs{\vec{U}^{1}} + ( \mathbf{M}\otimes\mathbf{I}_{N-1}\otimes\mathbf{M} ) \bs{\vec{U}^{2}}
			+ ( \mathbf{I}_{N-1}\otimes\mathbf{M}\otimes\mathbf{M} ) \bs{\vec{U}^{3}} = \bs{\vec{0}}
		\]
		by summing the first three equations of sub-system involving $\bs{\vec{U}^{1}}$-$\bs{\vec{U}^{3}}$, 
		which is exactly the divergence-free condition satisfied by the numerical eigen-functions.
		Thus, all eigen-solutions associated with non-zero eigenvalues in Theorem \ref{Thm2-LamU3dSol} are exactly divergence-free solutions at the discrete level.

		\item 
		When ${\lambda}^{N}=0$, the eigenfunctions in \eqref{eq:Lam0U3dSol} correspond to  
		\begin{align*}
		\bs{u}
		& = \bs{\Phi}_{m,n,l}^1(\bs{x}) + \bs{\Phi}_{m,n,l}^2(\bs{x})  + \bs{\Phi}_{m,n,l}^3(\bs{x})  
		\\
		& = \big( \phi_{m}(x_1)\psi_{n+1}(x_2)\psi_{l+1}(x_3), \psi_{m+1}(x_1)\phi_{n}(x_2)\psi_{l+1}(x_3), \psi_{m+1}(x_1)\psi_{n+1}(x_2)\phi_{l}(x_3) \big)^{\intercal}
		\\
		& = \nabla\big( \psi_{m+1}(x_1)\psi_{n+1}(x_2)\psi_{l+1}(x_3) \big),
		\end{align*}
		where $\big\{{\Psi}_{m,n,l}(\bs{x})=\psi_{m+1}(x_1)\psi_{n+1}(x_2)\psi_{l+1}(x_3)\big\}_{m,n,l=1}^{+\infty}$ is the complete basis for $H_0^1(\Omega)$.
	\end{itemize}
\end{rem}

\section{Miscellaneous extensions}
\label{Sec4}
The proposed Gauss's law-preserving spectral methods can readily be applied to more general cases. 
Here, we discuss three direct extensions, i.e., problems (i) with non-homogeneous boundary condition, (ii) on irregular geometry, and (iii) with variable coefficients.

\subsection{Extension to non-homogeneous boundary condition}
\label{Sec4-1}
Now, we consider double-curl source problems with non-homogeneous boundary condition as follows:
\begin{equation}
\label{eq:MaxwSysNBC}
\nabla\times\nabla\times\bs{u} + \kappa\,\bs{u}=\bs{f},  \quad \nabla\cdot\bs{u} = \rho \;\; \text{in} \;\; \Omega,\quad \bs{n}\times\bs{u} = \bs{b}\;\; \text{on} \;\; \partial\Omega.
\end{equation}
It suffices to construct $\bs{u}_{\rm{b}}$ such that it satisfies $\bs{n}\times\bs{u}_{\rm{b}} =\bs{b}$ on $\partial\Omega$.
For instance, $\bs{u}_b$ in 2D can be constructed by 
\begin{equation}
\label{eq:ub2D}
\bs{u}_{\rm{b}}(\bs{x})=\sum_{m=0}^{N-1}
 \big( \hat{b}_{m}^{3}\phi_{m}(x_1) \psi_{0}(x_2) +\hat{b}_{m}^{4}\phi_{m}(x_1) \psi_{1}(x_2)    \big) \bs{i}
 + \sum_{n=0}^{N-1} \big( \hat{b}_{n}^{1}\psi_{0}(x_1) \phi_{n}(x_2)+\hat{b}_{n}^{2}\psi_{1}(x_1) \phi_{n}(x_2)  \big) \bs {j}
\end{equation}
with $\bs{i}=(1,0)^{\intercal}$, $\bs{j}=(0,1)^{\intercal}$, $\psi_{0}(\xi)=(1-\xi)/2$, $\psi_{1}(\xi)=(1+\xi)/2$, and
\begin{equation*}
\begin{aligned}
& \hat{b}_{n}^{1} = -\int_{\Lambda} b(x_1,x_2)|_{x_1=-1} \phi_{n}(x_2) \,\mathrm{d}x_2, \quad 
\hat{b}_{n}^{2} = \int_{\Lambda} b(x_1,x_2)|_{x_1=1} \phi_{n}(x_2) \,\mathrm{d}x_2, 
\\
& \hat{b}_{m}^{3} = \int_{\Lambda} b(x_1,x_2)|_{x_2=-1} \phi_{m}(x_1) \,\mathrm{d}x_1, \quad 
\hat{b}_{m}^{4} = -\int_{\Lambda} b(x_1,x_2)|_{x_2=1} \phi_{m}(x_1) \,\mathrm{d}x_1.
\end{aligned}
\end{equation*}
Next, with the decomposition $\bs{u}=\bs{u}_{\rm{in}}+\bs{u}_{\rm{b}}$, one can obtain $\bs{u}_{\rm{in}}$ via 
\begin{equation}
\label{eq:MaxwSys_uin}
\begin{aligned}
& \nabla\times\nabla\times\bs{u}_{\rm{in}} + \kappa\,\bs{u}_{\rm{in}}=\bs{\tilde{f}},  \quad 
\nabla\cdot\bs{u}_{\rm in} = \tilde{\rho} \;\; \text{in} \;\; \Omega,\quad \bs{n}\times\bs{u}_{\rm in} = \bs{0}\;\; \text{on} \;\; \partial\Omega,\\
& \bs{\tilde{f}}=\bs{f}-\nabla\times\nabla\times\bs{u}_{\rm{b}} - \kappa\bs{u}_{\rm{b}},\quad \tilde{\rho}=\rho-\nabla\cdot\bs{u}_{\rm{b}} \;\; \text{in} \;\; \Omega.
\end{aligned}
\end{equation} 
Thus, it is direct to apply the proposed method to compute $\bs{u}$. This procedure is applicable to the 3D case without difficulty.

\subsection{Extension to complex geometry}
\label{Sec4-2}
Next, we discuss how to extend the proposed method to irregular geometry diffeomorphic to square or cube ${\Lambda}^d$, $d=2,3$. 
The following covariant Piola transformation (see \cite[p.~39]{Ciarlet1988} more details) finds useful for constructing $\bs{H}({\rm curl})$-conforming basis on such domains.
\begin{lm}
	\label{Lem3-Covar-Piola}
	Let $\mathcal{F}$ be a diffeomorphism from ${\Lambda}^d$ onto $\mathcal{F}({\Lambda}^d)=\Omega\subset\mathbb{R}^{d}$, 
	the covariant Piola transformation is given by
	\begin{equation*}
		\mathcal{F}^{\rm curl}[\bs{u}](\bs{x}) = ( \partial_{\bs{\xi}}{\bs{x}} )^{-\intercal}\bs{u}\circ\mathcal{F}^{-1}(\bs{x}), \quad \; \bs{x}=\mathcal{F}(\bs{\xi}),
	\end{equation*}
	which satisfies the identity
	\begin{equation*}
		\nabla_{\bs{x}}\times \mathcal{F}^{\rm curl}[\bs{u}] = \nabla_{\bs{\xi}}\times\bs{u}.
	\end{equation*}
\end{lm}
Note that for $\Omega$ diffeomorphic to $\Lambda^d$, $d=2,3$, $\mathcal{F}$ can be chosen as the commonly-used Gordon-Hall transformation (see \cite{Gordon-Hall-1973}).  
Then equipped with the covariant Piola transformation, one can define the $\bs{H}({\rm curl})$-conforming spectral approximation basis on domain $\Omega$ 
by $\bs{\tilde{\Phi}}(\bs{x})=\mathcal{F}^{\rm curl}[\bs{\Phi}](\bs{x})$, 
where $\bs{\Phi}$ represents the proposed basis on the reference square or cube in Propositions \ref{Prop1-DiscSpace2d} and \ref{Prop2-DiscSpace3d}. 
While for more complex domains, it is possible to combine the proposed method with spectral element methods, 
where study along this line will be reported in a separate work.

\subsection{Extension to problems with variable coefficients}
\label{Sec4-3}
Finally, we consider the double-curl problem with variable coefficients, which arises naturally when the permittivity or permeability in Maxwell's equations vary with space 
or through coordinate transformation of the original double-curl problem with constant coefficients in general domain to reference domain, as described in the previous subsection. 
To be specific, let us consider
\begin{align}
\label{eq:MaxwVarProb}
&  \nabla\times\big({\bs \alpha}(\bs{x})\nabla\times{\bs{u}}\big) + \kappa\,{\bs u}={\bs f}\;\; \text{in}\;\;  \Omega, \quad  
 \nabla\cdot\bs{u} = \rho \;\; \text{in} \;\;\Omega,  \quad
\bs{n}\times\bs{u} = \bs{0} \;\; \text{on}\;\;\partial\Omega.
\end{align}
Due to the variable coefficient $\bs{\alpha}(\bs{x})$ involved in the above equation, the proposed mixed spectral approximation scheme leads to dense coefficient matrix and cannot be solved efficiently by a direct method. 
Nevertheless, the fast solution algorithm for problem with constant coefficients described in previous sections could be used as a preconditioner to solve problem \eqref{eq:MaxwVarProb} efficiently 
when combined with matrix-free iterative methods, as verified by Example \ref{Exmp5-3} in Section \ref{Sec5}.

\section{Representative numerical experiments} 
\label{Sec5}
In this section, we conduct several numerical experiments to demonstrate the accuracy and efficiency of the proposed Gauss's law-preserving spectral methods. 
We first employ manufactured solutions to test the convergence rates of the proposed methods for double-curl source problems in both two and three dimensions. 
Numerical results of the classical finite difference Yee scheme (see e.g. \cite{Yee-1966}) are also included for comparison in terms of accuracy, efficiency and scalability.
Then we subject the proposed methods to challenging numerical tests for indefinite systems with highly oscillatory solutions and problems with variable coefficients. 
Finally, we investigate the performance of the proposed spectral eigen-solver for two and three dimensional Maxwell's eigenvalue problems.

\begin{example}[\bf Convergence test for Maxwell's problems in 2D] 
	\label{Exmp5-1}
	Firstly, we consider problem \eqref{eq:MaxwSouProb} with the manufactured solution $\bs{u} = (u_1,u_2)^{\intercal}$ given by
	\begin{equation*}
	\label{eq:ConvTrig_u2d}
	u_1(x_1,x_2)  =  (\cos(\pi{x_1}) + \sin(\pi{x_1}) ) \sin(\pi{x_2}),
	\quad \;
	u_2(x_1,x_2) = \sin(\pi{x_1}) ( -\cos(\pi{x_2})  + \sin(\pi{x_2}) ).
	\end{equation*}
	Accordingly, the source terms $\bs{f}$ and $\rho$ in \eqref{eq:MaxwSouProb} are chosen such that the above analytic solution satisfies \eqref{eq:MaxwSouProb}. 
\end{example}

\begin{figure}[!htb]
	\centering
	\subfigure[Errors versus $N$] 
	{\includegraphics[scale=0.40]{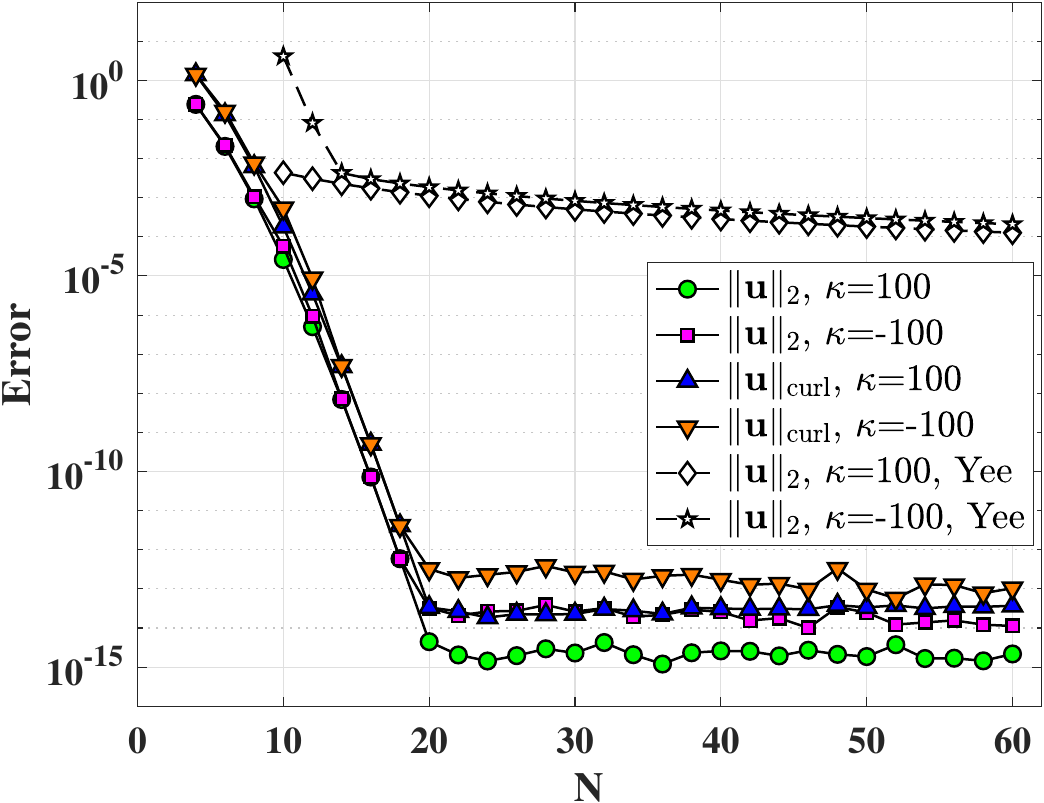}} 
		\subfigure[CPU time versus DoFs] 
	{\includegraphics[scale=0.40]{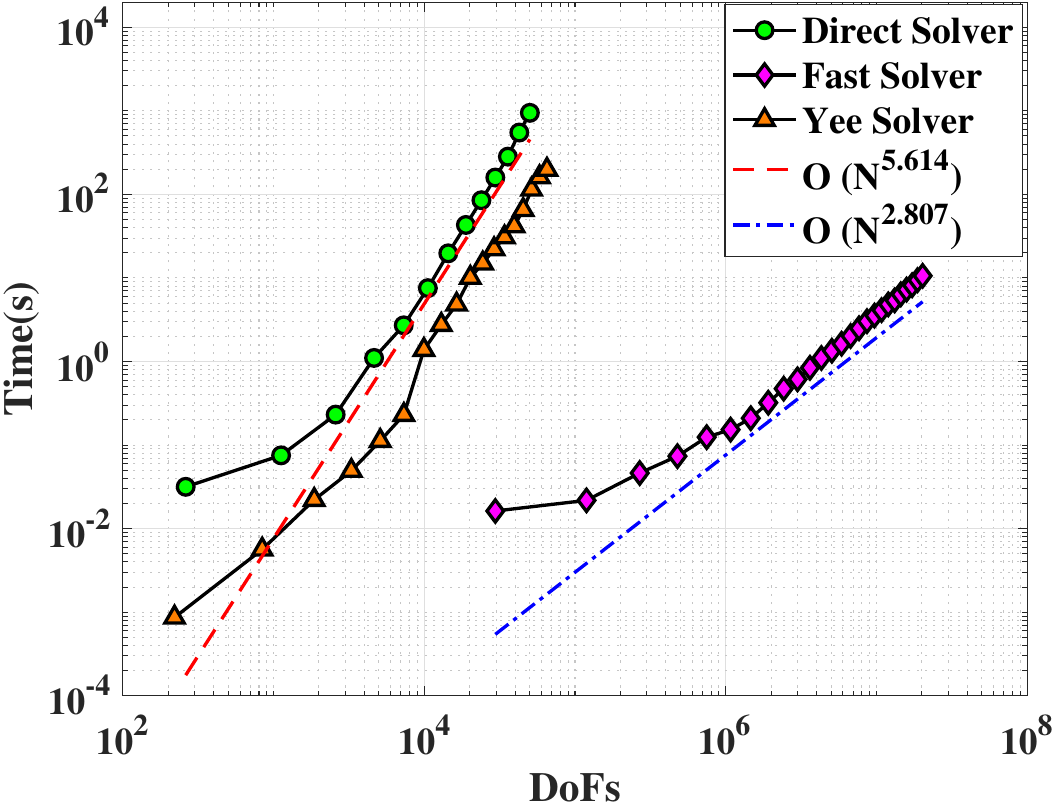}}
	\caption{Example \ref{Exmp5-1}: convergence and computational cost tests in 2D. 
		(a) $L^2$- and $H(\rm curl)$-errors of $\bs{u}$ versus polynomial order $N$ for $\kappa = \pm{100}$;
		(b) comparison on the computational time versus DoFs for $\kappa=100$ of the direct, Yee and the proposed method.
		}
	\label{Fig1_Exmp1_ConvComp}
\end{figure}

In Figure~\ref{Fig1_Exmp1_ConvComp}(a), we plot the discrete $L^2$- and $H(\rm curl)$-errors on the semi-log scale against various $N$ by using the proposed Gauss's law preserving spectral method. 
Comparisons of errors with the widely-used finite difference Yee scheme is also included. 
We observe the expected exponential convergence of the proposed method for both positive and negative parameters $\kappa$ as $N$ increases, 
while it converges quadratically for the Yee scheme. 
Moreover, it can be observed that the proposed method reaches machine precision when $N\geq{20}$, 
while the discrete $L^2$-errors of the Yee method with the same size of computational grid is around $10^{-3}$.

In Figure~\ref{Fig1_Exmp1_ConvComp}(b), we depict the computational time consumption on log-log scale versus DoFs for the direct method based on LU factorization, 
the Yee method and our proposed method with fast solution algorithm for fixed $\kappa=100$.
Hereinafter, DoFs denotes the total degrees of freedom and DoFs$=3(N-1)^2+2(N-1)$ for the proposed method. 
Note that our proposed algorithm is matrix-free, making it possible to deal with large-scale computation. 
While direct algorithm requires the storage of coefficient matrix, thus the problem size is greatly affected by memory limit. 
We list in below the maximum DoFs that we computed for these three methods and their corresponding CPU times:
\begin{itemize}
\item[(1)] the Yee scheme, Max DoFs$=65160$, CPU time$=197.7$s;
\item[(2)] the proposed spectral method with direct method,  Max DoFs$=50181$, CPU time$=945.9$s;
\item[(3)] the proposed spectral method with fast algorithm,  Max DoFs$=20269601$, CPU time$=10.61$s.
\end{itemize}
It can be seen that the curve of computational complexities for the Yee scheme and direct method shows an asymptotic trend of $\mathcal{O}(N^{5.614})$, 
while that for our fast algorithm exhibits complexity of $\mathcal{O}(N^{\log_{2}{7}})$ ($\log_2{7}\approx{2.807}$), 
agreeing well with the complexity analysis in Remark \ref{Remk1-Algo1}.
These results clearly demonstrate the significant superiority of our proposed algorithm in terms of accuracy, efficiency and scalability.

\begin{example}[\bf Maxwell's problems with sharp point source in 2D]
	\label{Exmp5-2}
	Then we consider problem \eqref{eq:MaxwSouProb} with the Gaussian point source 
	\begin{equation*}
	\bs{f}(x_1,x_2) = \nabla \times \Big(
	\exp\Big( -\frac{(x_1+0.5)^2+x_2^2}{\sigma^2} \Big) + \exp\Big( -\frac{(x_1-0.5)^2+x_2^2}{\sigma^2} \Big)
	\Big),
	\end{equation*}
	where $\sigma$=0.01.
\end{example}

\begin{figure}[!htb]
	\centering
	\subfigure[$\kappa=-100$] 
	{\includegraphics[scale=0.38]{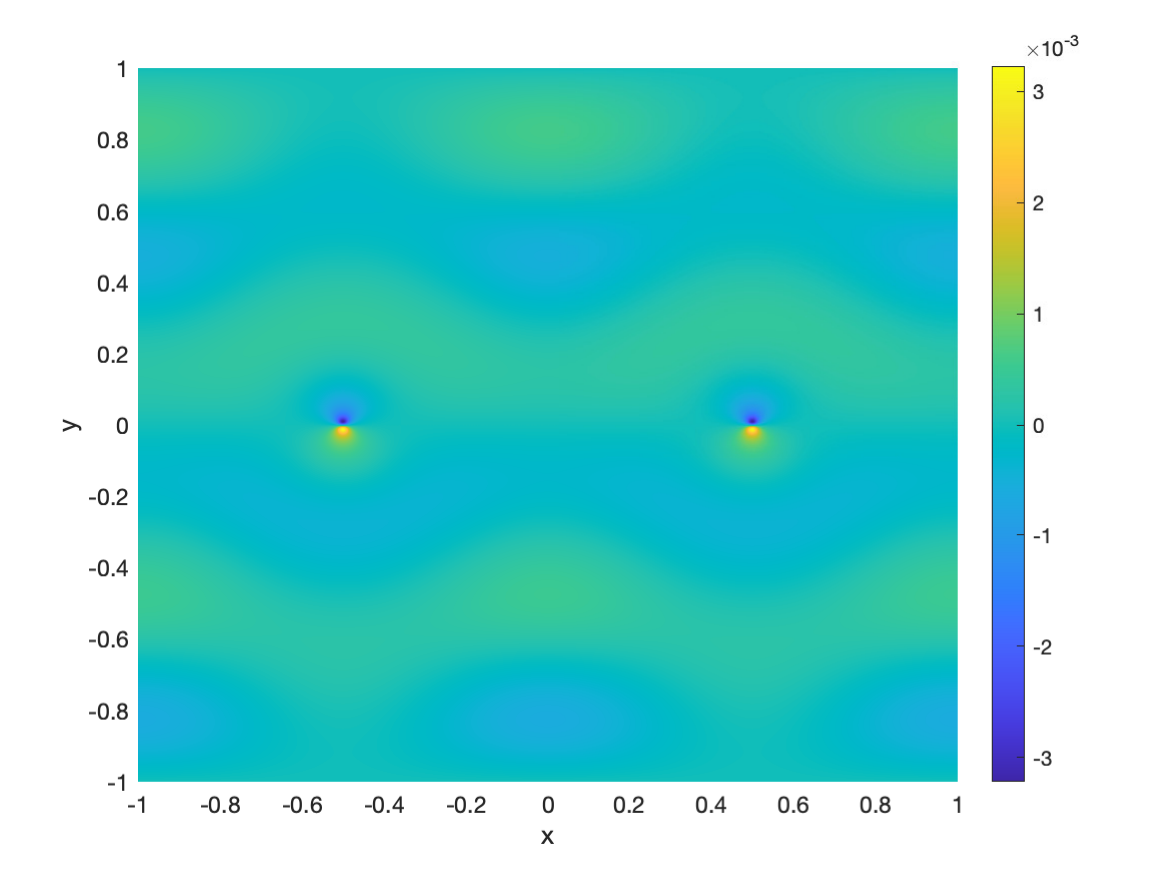}}
	\subfigure[$\kappa=-1600$] 
	{\includegraphics[scale=0.38]{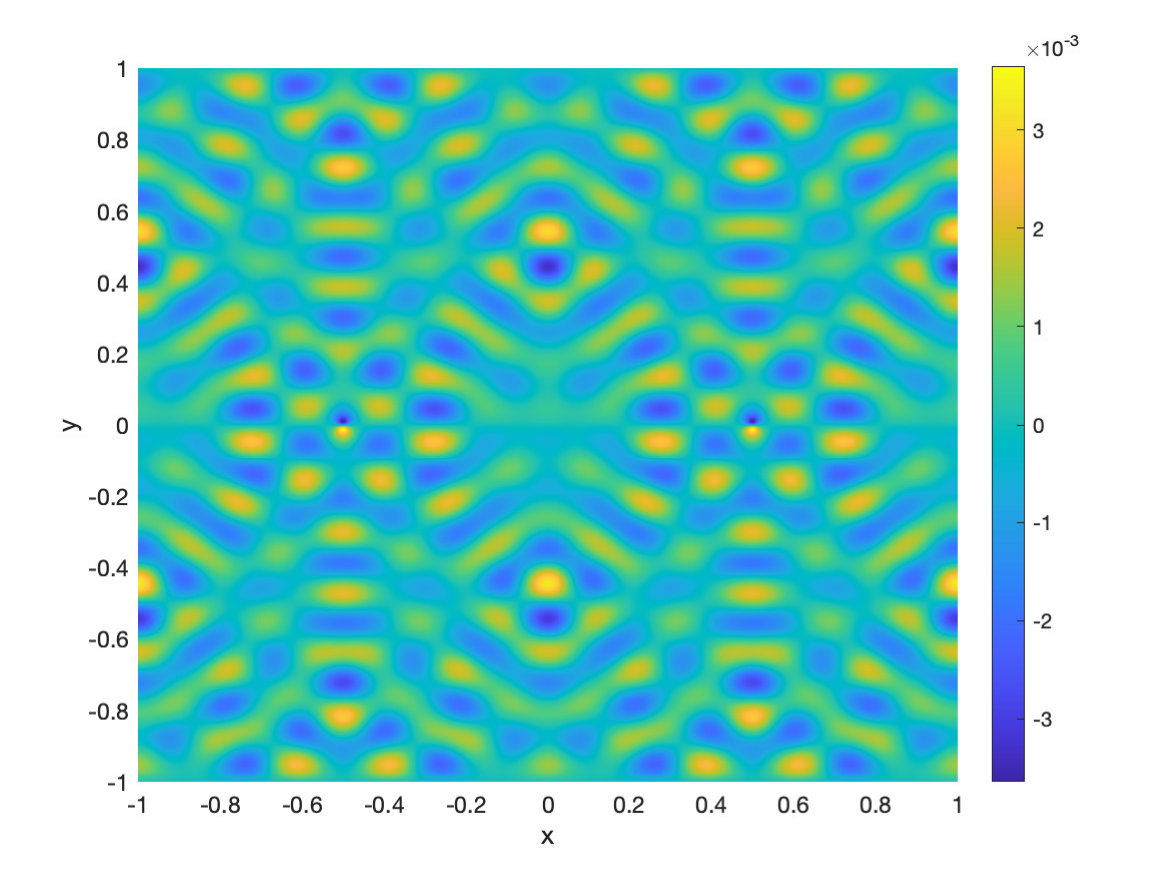}}
	\subfigure[$\kappa=-4900$] 
	{\includegraphics[scale=0.38]{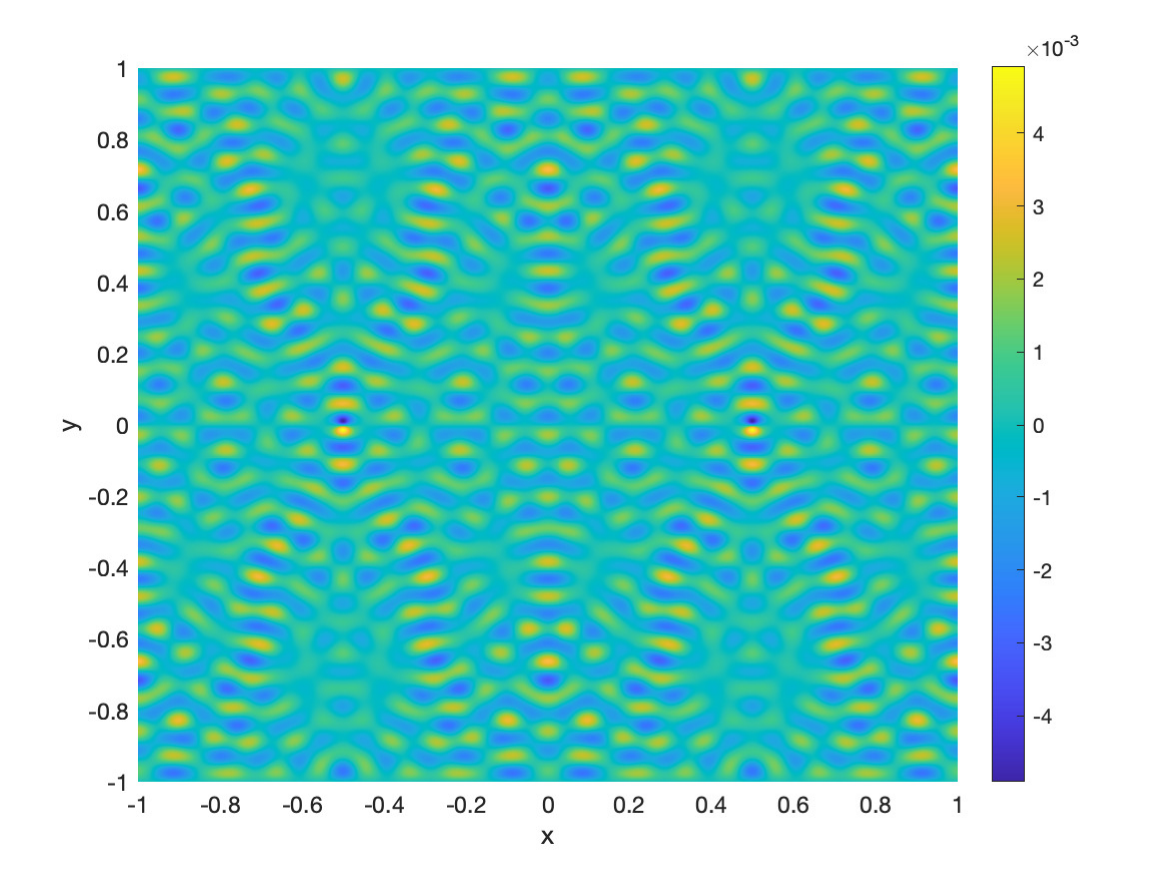}}
	\subfigure[$\kappa=-10000$] 
	{\includegraphics[scale=0.38]{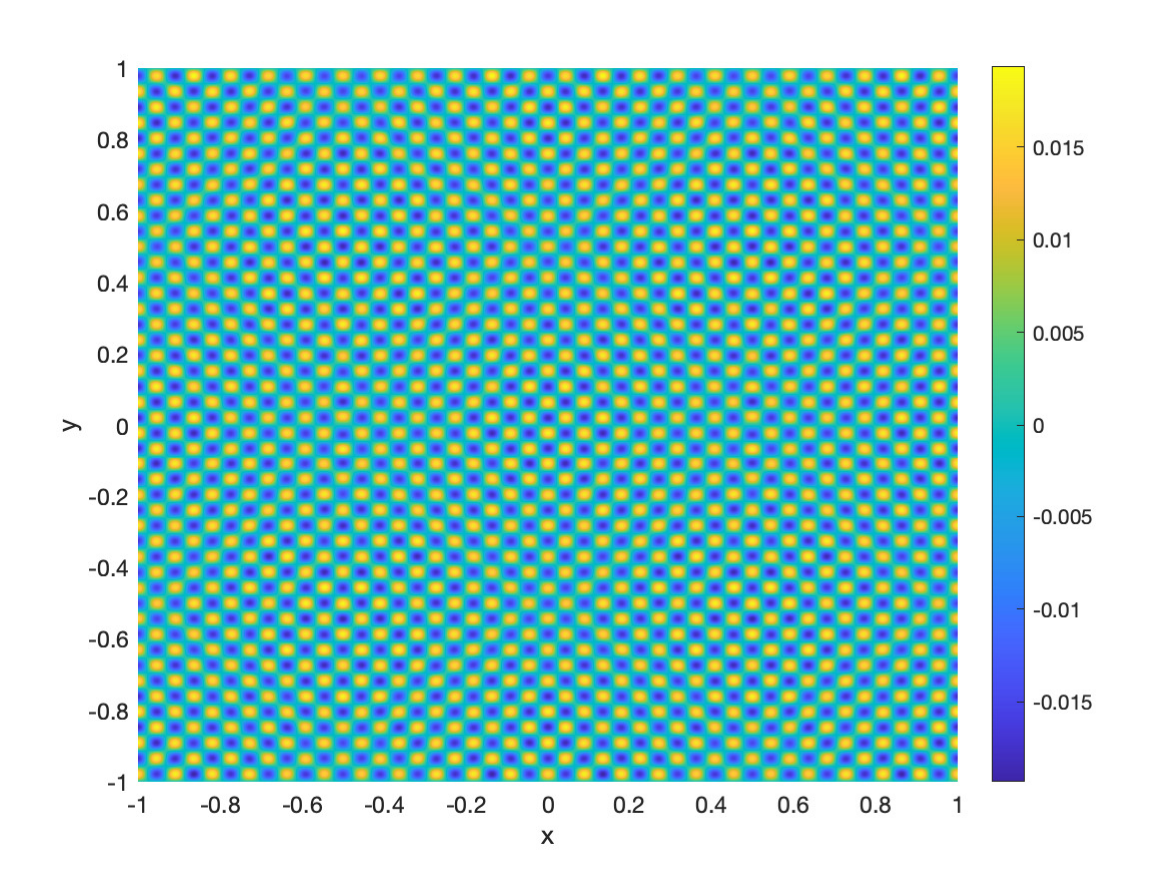}}
	\caption{Example \ref{Exmp5-2}: the profiles of the component $\mathbf{U}$ of the numerical solution versus different $\kappa$. 
		These results are obtained with a fixed polynomial order $N=512$.
	}
	\label{Fig2_Exmp2_kappa_U1Sol}
\end{figure}

Due to the source term here is sharp, thus the accuracy required for solving this problem is very high.
At the same time, the solution itself becomes increasingly oscillatory when $\kappa$ is chosen as large negative values, which poses great difficulties for numerical simulation.
We challenge our spectral algorithm via indefinite system with highly oscillatory solutions by choosing the parameter $\kappa$ as large negative values -100, -1600, -4900 and -10000 successively.
We take the obtained numerical solution with $N=1024$ as the reference solution, and compute the approximate solution as $N$ varies from 4 to 512. 
We find that the discrete $L^2$-errors between the numerical solution and the reference solution versus $N$ exhibit the expenential convergence trend,
and the discrete $L^2$-errors between the numerical solution for $N=512$ and the reference solution for $\kappa$ as values -100, -1600, -4900 and -10000
are 9.153e-9, 9.201e-9, 9.309e-9 and 9.481e-9, respectively.
Figure~\ref{Fig2_Exmp2_kappa_U1Sol} shows the component $\mathbf{U}$ of the numerical solution with $N=512$ for different values of $\kappa$, 
which illsutrates that our spectral algorithm provides good approximations, even for solutions with high oscillation.

\begin{example}[\bf Maxwell's problems with variable coefficients in 2D] 
	\label{Exmp5-3}
	Now we consider the problem \eqref{eq:MaxwVarProb} on $\Omega=\Lambda^2$, 
	where the coefficient $\alpha(\bs{x})$ is given by
	\begin{equation*}
		\alpha(x_1,x_2) = 5.5 + 0.2\sum_{i=1}^{5}\sum_{j=1}^{5} \tanh\Big(\frac{\sqrt{(x_1-a_i)^2+(x_2-b_j)^2} - 0.16}{\sqrt{2}\gamma}\Big),
	\end{equation*}
	and source term $\bs{f}$ is 
	\begin{equation*}
		\bs{f}=(f_1,f_2)^{\intercal}, \quad
		f_1={\rm exp}\Big( \frac{x_1^2+x_2^2}{\sigma^2} \Big), \quad 
		f_2 = 0.
	\end{equation*}
	Here, we choose $\kappa=-400$, $\gamma=0.04$, $\sigma=0.05$ and 
	\begin{equation*}
	a_1=b_1=-0.8, \; a_2=b_2=-0.4, \; a_3=b_3=0, \; a_4=b_4=0.4, \; a_5=b_5=0.8.
	\end{equation*}
\end{example}

\begin{figure}[!htb]
	\centering
	\subfigure[Error versus $N$] 
	{\includegraphics[scale=0.38]{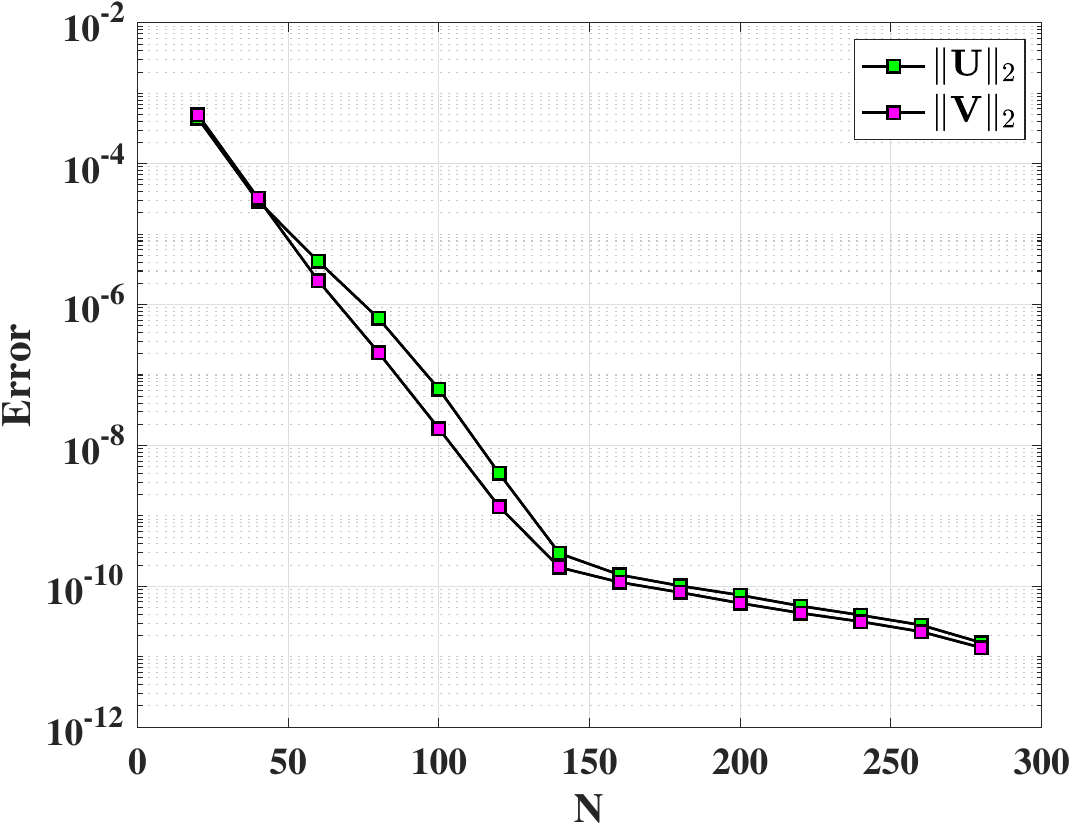}}
	\subfigure[$\alpha(\bs{x})$] 
	{\includegraphics[scale=0.38]{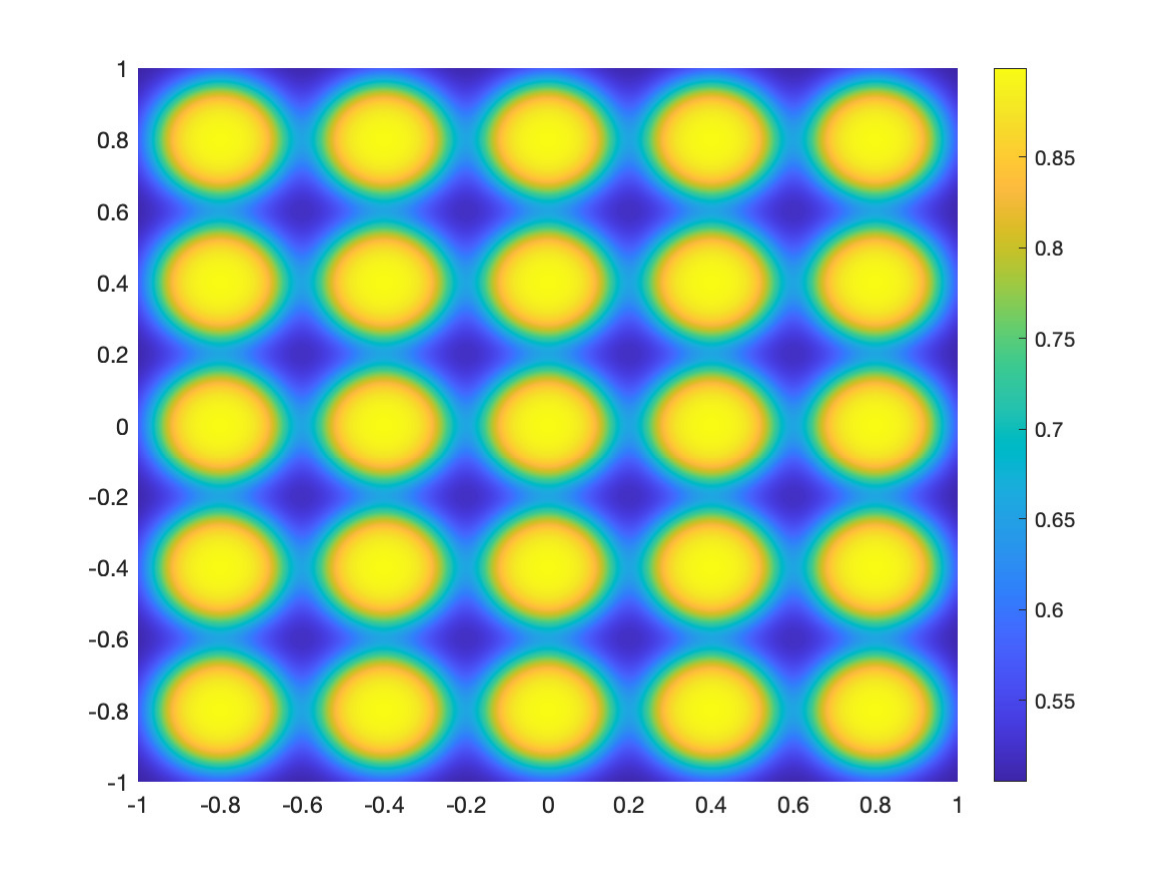}}
	\subfigure[$\mathbf{U}$] 
	{\includegraphics[scale=0.38]{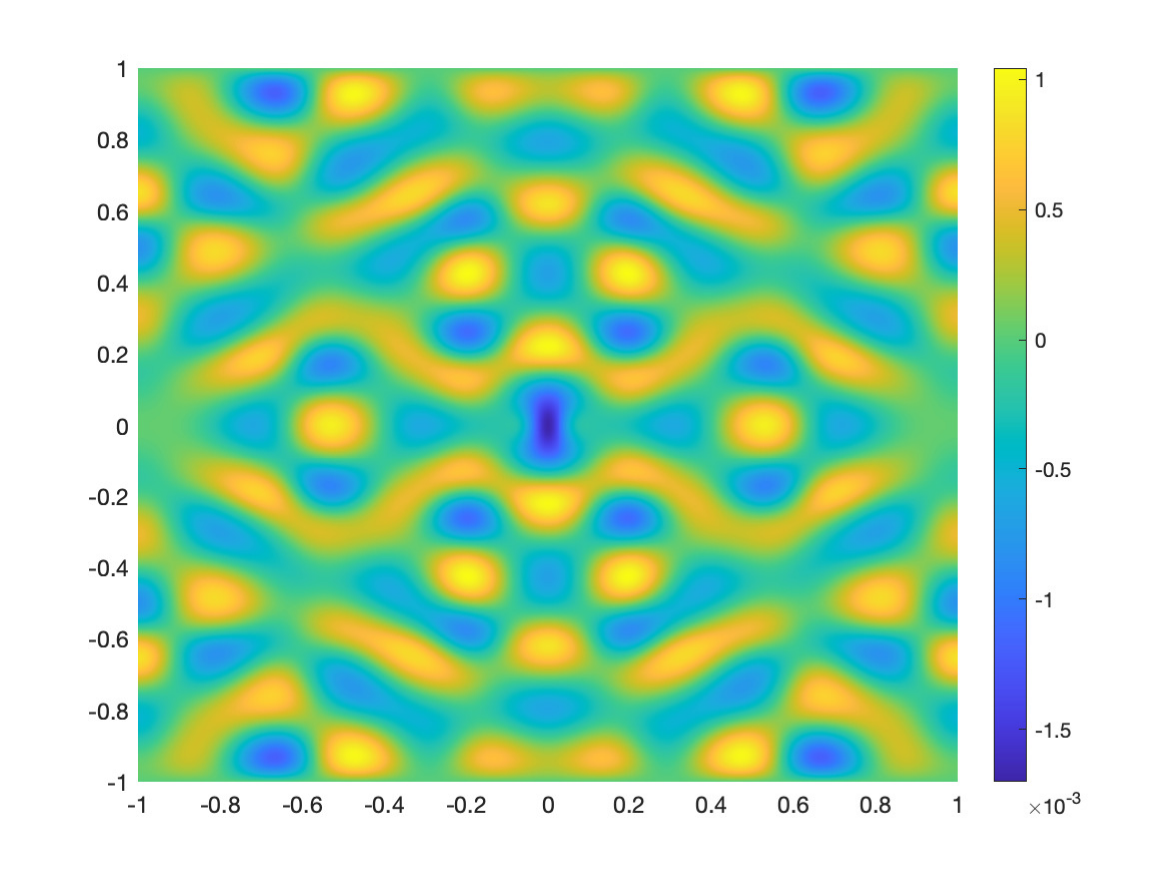}}
	\subfigure[$\mathbf{V}$] 
	{\includegraphics[scale=0.38]{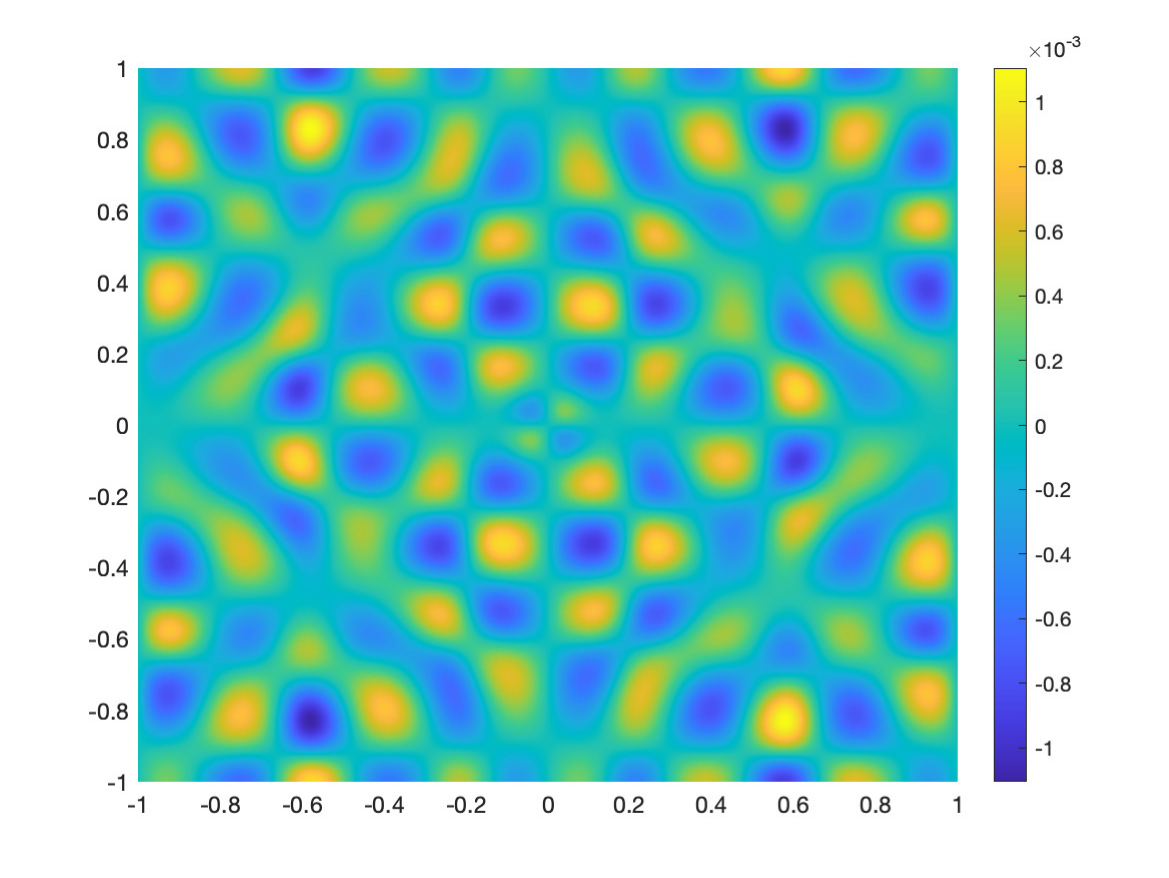}}
	\caption{Example \ref{Exmp5-3}: 
		(a) $L^2$-errors of $\mathbf{U}$ and $\mathbf{V}$ versus $N$; 
		(b-d) the profiles of $\alpha(\bs{x})$ and two components $\mathbf{U}$ and $\mathbf{V}$ of the numerical solution with $\kappa=-400$ for $N=200$. 
	}
	\label{Fig3_Exmp3_UVSolN200}
\end{figure}

We then employ the proposed algorithm for problem with constant coefficients as a preconditioner to solve the above problem using matrix-free preconditioned GMRES iterative method.  
We take the obtained numerical solution with $N=300$ as the reference solution, and compute the approximate solution as $N$ varies from $20$ to $280$. 
We observe that the discrete $L^2$-errors of two components $\mathbf{U}$ and $\mathbf{V}$  between the numerical solution and the reference solution versus $N$ 
exhibit the expenential convergence trend in Figure~\ref{Fig3_Exmp3_UVSolN200}(a).
Besides, we depict the profiles of the coefficient $\alpha(\bs{x})$ and the two components $\mathbf{U}$ and $\mathbf{V}$ of the numerical solution with $N=200$ in Figure~\ref{Fig3_Exmp3_UVSolN200}(b-d).

\begin{example}[\bf Convergence test for Maxwell's problems in 3D] 
	\label{Exmp5-4}
	Next, we consider problem \eqref{eq:MaxwSouProb} in three dimensions with the manufactured solution $\bs{u}=(u_1,u_2,u_3)^{\intercal}$ given by
	\begin{align*}
	u_1(x_1,x_2,x_3) &= 2\cos(\pi(x_1+1)/2) \sin(\pi(x_2+1)/2) \sin(\pi(x_3+1)/2) + (x_1^2-1)(x_2^2-1)(x_3^2-1),
	\\
	u_2(x_1,x_2,x_3) &= - \sin(\pi(x_1+1)/2) \sin(\pi(x_2+1)/2) \cos(\pi(x_3+1)/2) + (x_1^2-1)(x_2^2-1)(x_3^2-1),
	\\
	u_3(x_1,x_2,x_3) &= - \sin(\pi(x_1+1)/2) \sin(\pi(x_2+1)/2) \cos(\pi(x_3+1)/2) + (x_1^2-1)(x_2^2-1)(x_3^2-1).
	\end{align*}
\end{example}

\begin{figure}[!htb]
	\centering
	\subfigure[Errors versus $N$] 
	{\includegraphics[scale=0.40]{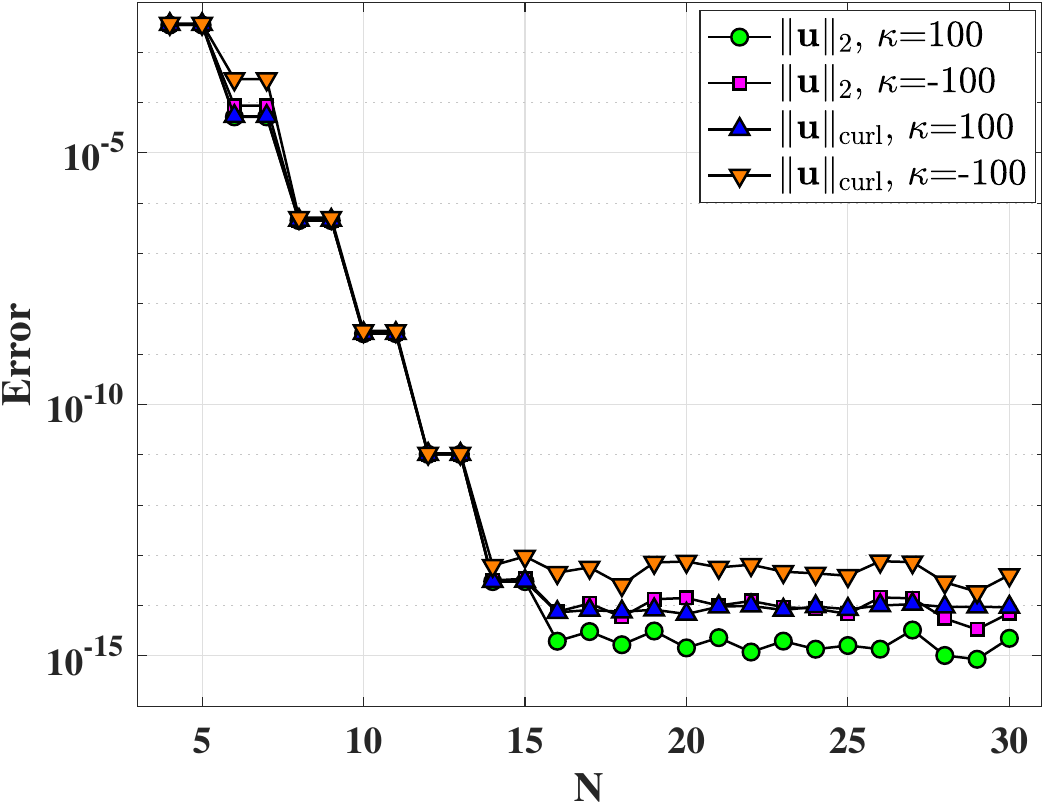}}
	\subfigure[CPU time versus DoFs] 
	{\includegraphics[scale=0.40]{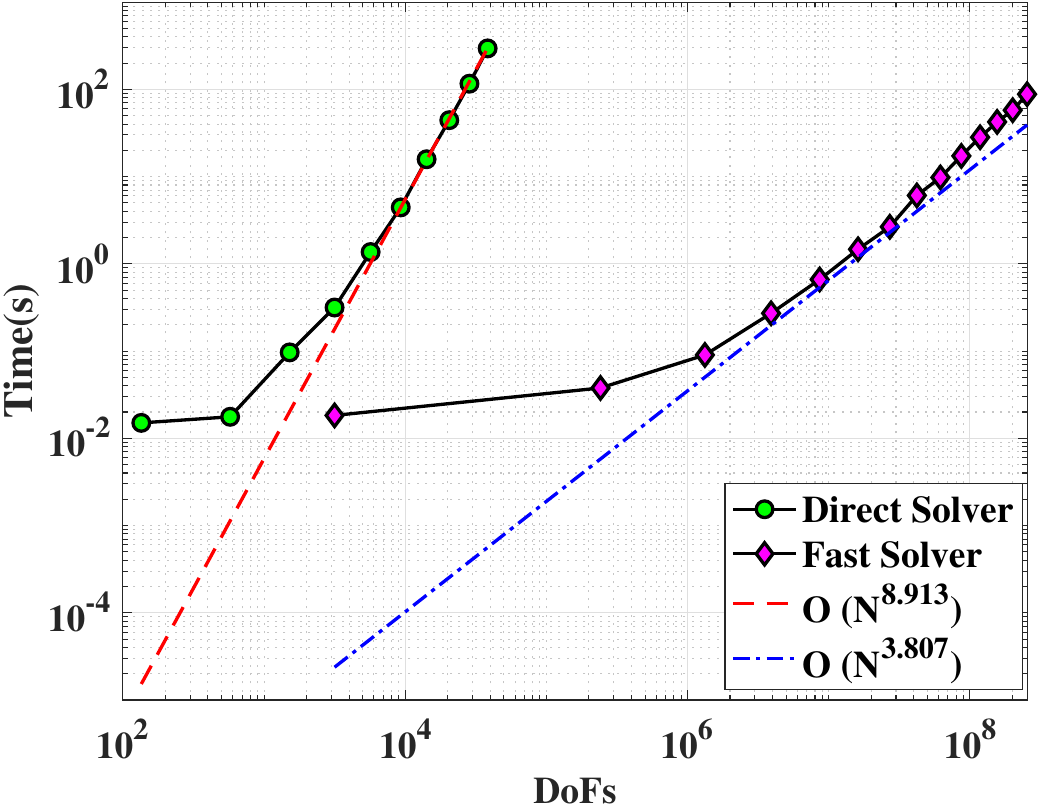}}
	\caption{Example \ref{Exmp5-4}: convergence and comparison tests of curl-curl problem in 3D. 
		(a) $L^2$- and $H(\rm curl)$-errors of $\bs{u}$ versus polynomial order $N$ for $\kappa = \pm{100}$;
		(b) comparison on the computational time versus DoFs for $\kappa=100$ of the direct and our proposed method.
	}
	\label{Fig4_Exmp4_ConvComp}
\end{figure}

In Figure~\ref{Fig4_Exmp4_ConvComp}(a),  we plot the discrete $L^2$- and $H(\rm curl)$-errors on the semi-log scale against various $N$ 
by using the proposed spectral algorithm. 
We observe the expected exponential convergence of the proposed method for both positive and negative parameters $\kappa$ as $N$ increases.

In Figure~\ref{Fig4_Exmp4_ConvComp}(b), we compare the computational time consumption on log-log scale versus DoFs for the direct method and our proposed method.
Here, DoFs$=4(N-1)^3+3(N-1)^2$ in the 3D case. The maximum DoFs by using the direct method is $38367$ with CPU time$=294.4$s, 
while the performance of our proposed algorithm far exceeds the former with the maximum DoFs$=254562399$ with CPU time$=87.97$s, 
due to the matrix-free and semi-analytic nature of the proposed method. 
Finally, the curve of computational time cost for the direct method shows a convergence trend of $\mathcal{O}(N^{8.913})$, 
and that for our method exhibits the optimal convergence rate $\mathcal{O}(N^{\log_{2}{14}})$ ($\log_2{14}\approx{3.807}$),
agreeing well with the complexity analysis in Remark \ref{Remk3-Algo2}.

\begin{example}[\bf Maxwell's problems with sharp point source in 3D]
\label{Exmp5-5}
Now we consider problem \eqref{eq:MaxwSouProb} in three dimensions with the Gaussian point source 
$\bs{f}=(f_1,f_2,f_3)^{\intercal}$ given by 
\begin{align*}
& f_1(x_1,x_2,x_3) = \exp\Big( -\frac{(x_1+0.9)^2+x_2^2+x_3^2}{\sigma^2} \Big) + \exp\Big( -\frac{(x_1-0.9)^2+x_2^2+x_3^2}{\sigma^2} \Big), 
\\ 
& f_2(x_1,x_2,x_3) = \exp\Big( -\frac{x_1^2+(x_2+0.95)^2+x_3^2}{\sigma^2} \Big) + \exp\Big( -\frac{x_1^2+(x_2-0.95)^2+x_3^2}{\sigma^2} \Big), 
\\
& f_3(x_1,x_2,x_3) = \exp\Big( -\frac{x_1^2+x_2^2+(x_3+0.98)^2}{\sigma^2} \Big) + \exp\Big( -\frac{x_1^2+x_2^2+(x_3-0.98)^2}{\sigma^2} \Big).
\end{align*}
Here, we choose $\sigma=0.05$ and $\kappa=-500$. 
\end{example}

\begin{figure}[!htb]
	\centering
	{\includegraphics[scale=0.40]{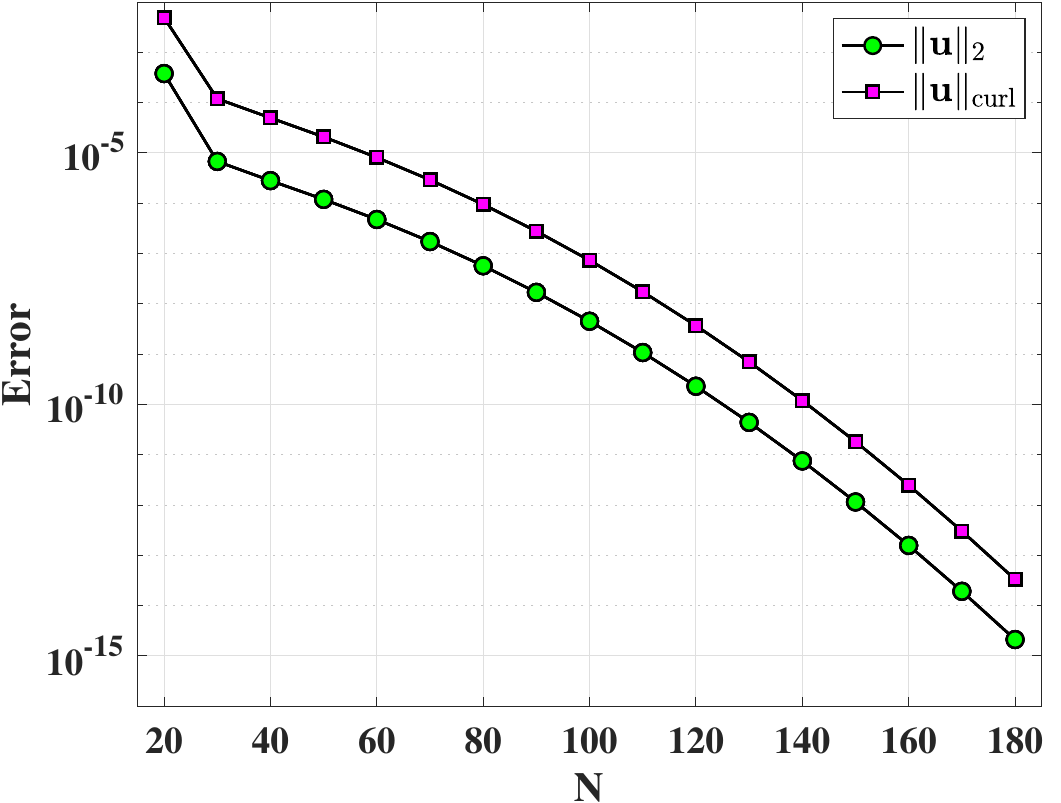}}
	\caption{Example \ref{Exmp5-5}: 
		$L^2$- and $H(\rm curl)$-errors of $\bs{u}$ versus polynomial order $N$ for $\kappa=-500$.
	}
	\label{Fig5_Exmp5_ConvN180}
\end{figure}

In this case, the solution is high oscillatory due to the indefiniteness of the system. 
Besides, the three components of the source term have a certain sharp phenomenon when approaching the six faces of the domain, 
which may lead to difficulties in calculating the numerical solution.
Similar to Example \ref{Exmp5-2}, we take the obtained numerical solution as the reference solution with the fixed $N=200$, 
and compute the numerical solution as $N$ varies from $20$ to $180$. 
In Figure~\ref{Fig5_Exmp5_ConvN180}, we plot the discrete $L^2$- and $H(\rm curl)$-errors on the semi-log scale against various $N$ from $20$ to $180$. 
The discrete $L^2$-error between the numerical solution with $N$=180 and the reference solution for $\kappa=-500$ is 2.122e-15.
We find that the discrete $L^2$- and $H(\rm curl)$-errors between the numerical solution and the reference solution versus $N$ exhibit the expenential convergence trend.

\begin{figure}[!htb]
	\centering
	\includegraphics[scale=0.40]{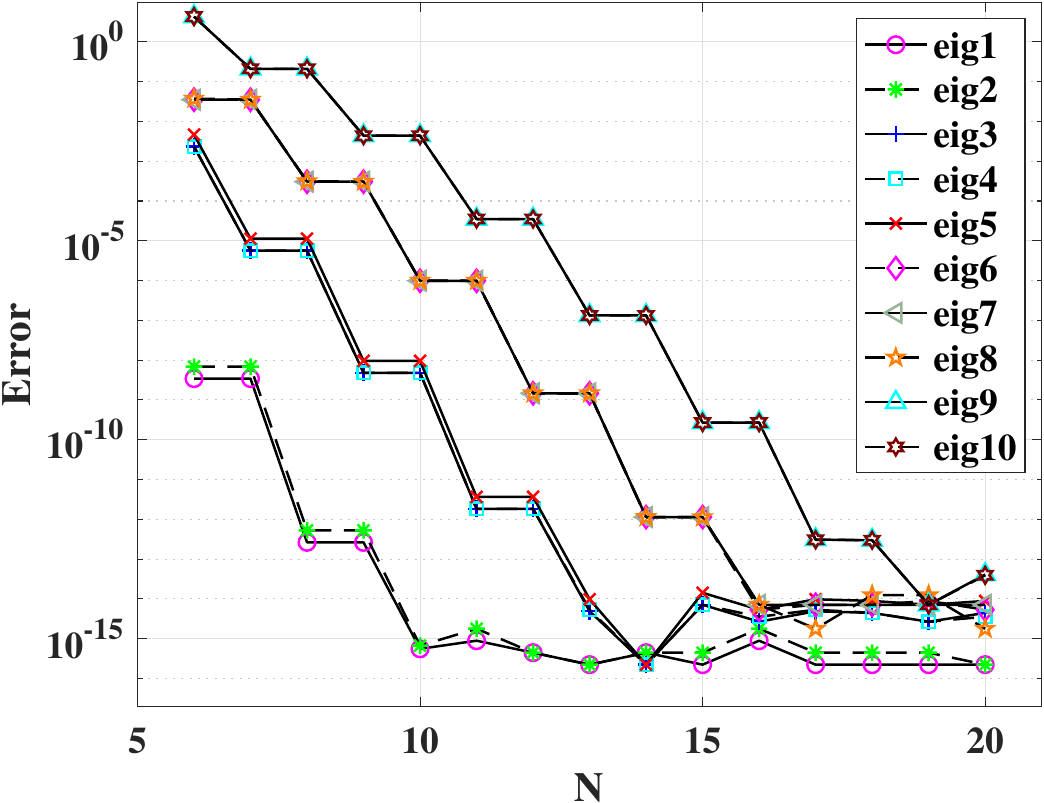} 
	\caption{Example \ref{Exmp5-6} in 2D: the $L^{\infty}$-error of the first 10 egienvalues versus $N$.}
	\label{Fig6_Exmp6_Lambda2DConv}
\end{figure}

\begin{example}[\bf Double-curl eigenvalue problems]
\label{Exmp5-6}
Finally, we focus on the double-curl eigenvalue problem \eqref{eq:MaxwEigProb} on $\Omega={\Lambda}^d$, $d=2,3$. 
\end{example}
The eigenvalues of the Maxwell's double-curl eigenvalue problem \eqref{eq:MaxwEigProb} on the reference square $\Omega=(-1,1)^2$ can be computed exactly, 
we refer to the literature \cite{Brenner-Li-Sun-2008,Fernandes-Raffetto-2002} and the reference therein for more details. 
The eigenvalues are of form $\{{4k}/{\pi}^2\}$, where $k=r^2+s^2$ with $r$ and $s$ are non-negative integers and $r^2+s^2\neq{0}$. 
We first investigate the convergence of some first numerical eigenvalues against various $N$ by our spectral method. 
In Table~\ref{Tab1_LambdaExt2DSol16}, we list the first 16 smallest exact eigenvalues with their multiplicities, and the coresponding numerical eigenvalues for fixed $N=20$. 
It can be observed that the $L^{\infty}$-errors of these eigenvalues are at the level of $10^{-14}$. 
In addition, Figure~\ref{Fig6_Exmp6_Lambda2DConv} presents the $L^{\infty}$-errors between the first $10$ smallest exact and numerical eigenvalues 
on the semi-log scale against various $N$, in which we observe the expected exponential convergence. 
When $N$ increases to $18$, the $L^{\infty}$-errors between these exact and numerical eigenvalues have reached machine precision.

\begin{table}[!ht]
	\normalsize
	\setlength{\belowcaptionskip}{6pt}
	\renewcommand{\arraystretch}{1.0}
	\centering	
	\caption{Example \ref{Exmp5-6} in 2D: the first 16 smallest eigenvalues with their multiplicities.}
	\label{Tab1_LambdaExt2DSol16}
	\setlength\tabcolsep{1mm}{
		\begin{tabular}{c|c|c|c|c|c|c|c}
			\hline \hline 
			\multicolumn{8}{c}{Exact solution: ${4\lambda}/{\pi^2}$ (multiplicity)} \\ \hline \hline
			1 (2)	& 4 (2)		& 8 (1)  & 10 (2) & 16 (2) & 18 (1)  & 25 (4)  & 29 (2) 
			\\ 
			2 (1)	& 5 (2)		& 9 (2)  & 13 (2) & 17 (2)  & 20 (2)  & 26 (2) & 32 (1)
			\\ \hline  \hline
			\multicolumn{8}{c}{Numerical solution: $4{\lambda}^N/{\pi^2}$} \\ \hline \hline
			\multicolumn{2}{c}{1.000000000000000}	&\multicolumn{2}{c}{9.000000000000005}	&\multicolumn{2}{c}{16.999999999999960}	&\multicolumn{2}{c}{25.000000000000046} 
			\\ \hline
			\multicolumn{2}{c}{1.000000000000000}	&\multicolumn{2}{c}{9.000000000000005}	&\multicolumn{2}{c}{16.999999999999960} &\multicolumn{2}{c}{26.000000000000043}
			\\ \hline
			\multicolumn{2}{c}{2.000000000000000}	&\multicolumn{2}{c}{10.000000000000007}	&\multicolumn{2}{c}{18.000000000000014} &\multicolumn{2}{c}{26.000000000000043}
			\\ \hline
			\multicolumn{2}{c}{3.999999999999996}	&\multicolumn{2}{c}{10.000000000000007}	&\multicolumn{2}{c}{19.999999999999954} &\multicolumn{2}{c}{29.000000000000046}
			\\ \hline
			\multicolumn{2}{c}{3.999999999999996}	&\multicolumn{2}{c}{13.000000000000002}	&\multicolumn{2}{c}{19.999999999999954} &\multicolumn{2}{c}{29.000000000000046}
			\\ \hline
			\multicolumn{2}{c}{4.999999999999996}	&\multicolumn{2}{c}{13.000000000000002}	&\multicolumn{2}{c}{24.999999999999964} &\multicolumn{2}{c}{31.999999999999915}
			\\ \hline
			\multicolumn{2}{c}{4.999999999999996}	&\multicolumn{2}{c}{15.999999999999960}	&\multicolumn{2}{c}{24.999999999999964} 
			\\ \hline
			\multicolumn{2}{c}{7.999999999999991}	&\multicolumn{2}{c}{15.999999999999960}	&\multicolumn{2}{c}{25.000000000000046} 
			\\ \hline
		\end{tabular}
	}
\end{table}

\begin{figure}[!htb]
	\centering
	\subfigure[Relative error versus index] 
	{\includegraphics[scale=0.39]{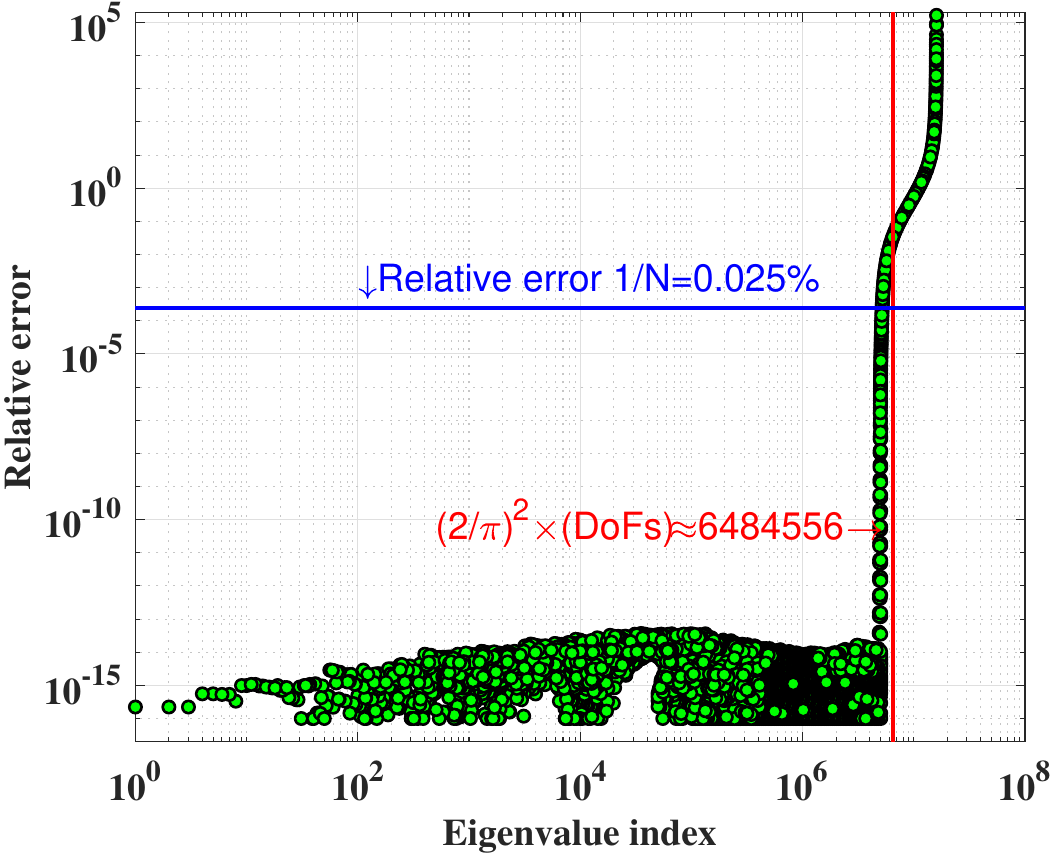}}
	\subfigure[Relative error region] 
	{\includegraphics[scale=0.40]{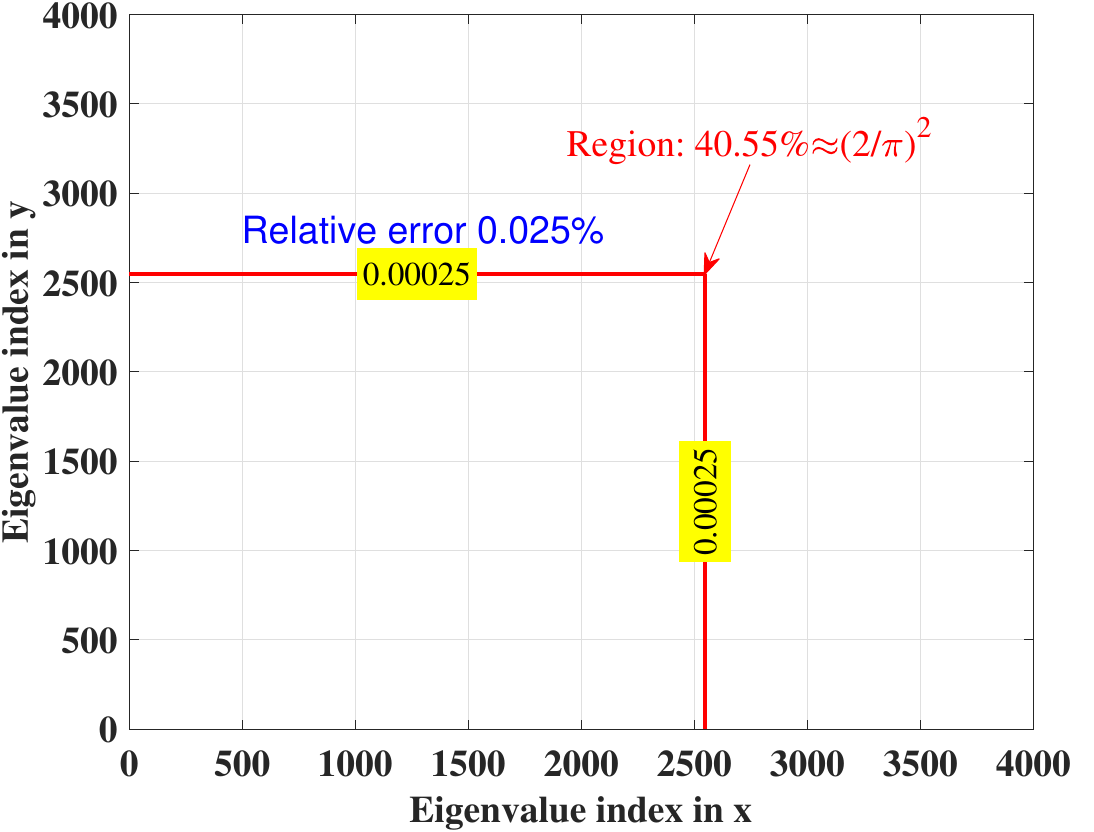}}
	\caption{Example \ref{Exmp5-6} in 2D:
		(a) relative error versus eigenvalue index for fixed $N=4000$;
		(b) relative error region for fixed $N=4000$.
	}
	\label{Fig7_Exmp6_RelErrTrustRegN4000}
\end{figure}

In \cite{Zhang-2015}, it is numerically verified that the percentage of ``trusted'' eigenvalues for solving the Laplace eigenvalue problem using Legendre spectral-Garlerkin method is around $(2/\pi)^d$, $d=1,2$. 
Here, the criterion for ``trusted'' eigenvalues is to have at least $\mathcal{O}(N^{-1})$ accuracy with polynomial order $N$. 
Motivated by this numerical result, we conduct numerical experiments in both two and three dimensions to explore the percentage of ``trusted'' eigenvalues for Maxwell's eigenvalue problems by our proposed method.  
In Figure~\ref{Fig7_Exmp6_RelErrTrustRegN4000}(a), we depict the relative error of exact and numerical eigenvalues on the log-log scale versus the eigenvalue index for fixed $N=4000$. 
Specifically, we draw a horizontal line  $y=1/N=0.025\%$ representing the error level of ``trusted'' eigenvalues and a vertical line $x=(2/\pi)^{2}\times({\rm DoFs})\approx{6484556}$, 
and find that the intersection point of these two lines lies on the error curve of the proposed method, 
which shows that there are about $(2/\pi)^{2}$ numerical eigenvalues for which the relative error converges at rate $\mathcal{O}(N^{-1})$.
In Figure~\ref{Fig7_Exmp6_RelErrTrustRegN4000}(b), we draw $0.025\%=1/4000$ relative error region of our proposed method.
Inside the region, the relative error is less than $0.025\%$. 
We see that the percentage of ``trusted'' eigenvalues in 2D case increases to $40.55\%$($\approx(2/\pi)^2$) for our spectral method.

\begin{figure}[!htb]
	\centering
	\includegraphics[scale=0.40]{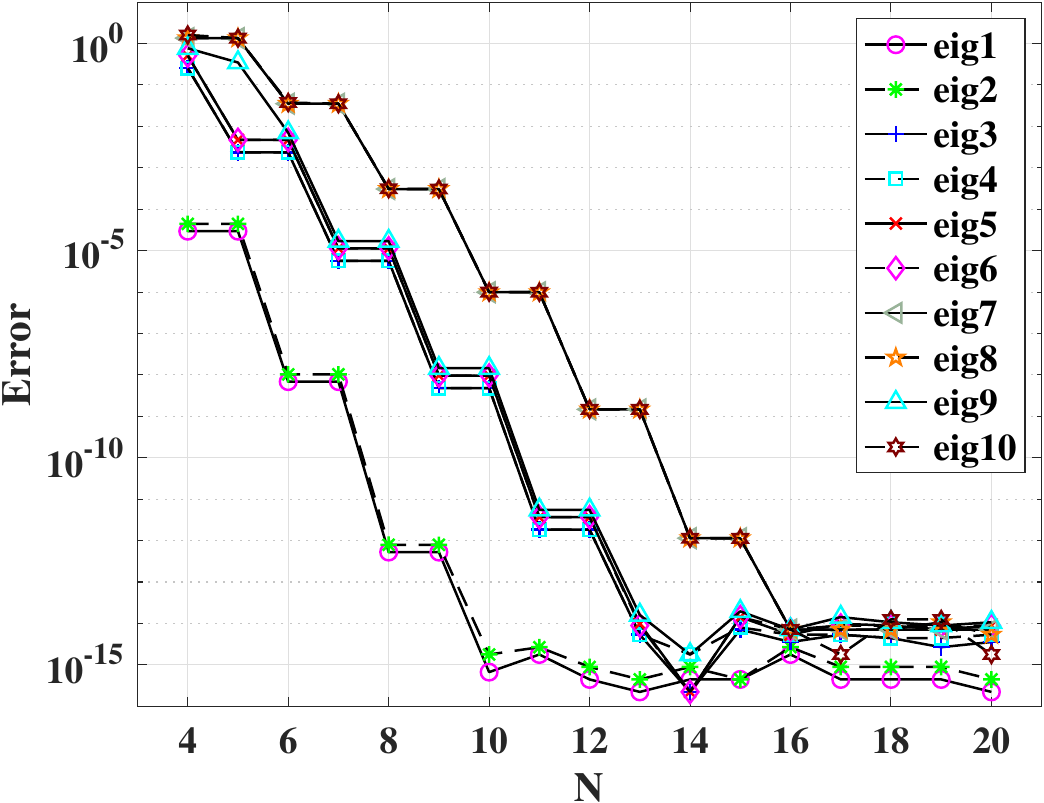} 
	\caption{Example \ref{Exmp5-6} in 3D: the $L^{\infty}$-error of the first 10 eigenvalues versus $N$.}
	\label{Fig8_Exmp7_Lambda3DConv}
\end{figure}

Finally, we focus on the double-curl eigenvalue problem \eqref{eq:MaxwEigProb} on the reference cube $\Omega=(-1,1)^3$. 
Similarly, the exact eigenvalues can be obtained analytically, we refer to the literature \cite{Bramble-Kolev-Pasciak-2005} and the reference therein for more details. 
The exact eigenvalues are of form $\{{4k}/{\pi}^2\}$, where $k=k_1^2+k_2^2+k_3^2$ and $\{k_i\}_{i=1}^3$ are non-negative integers satisfying $k_1k_2+k_2k_3+k_3k_1>0$. 
Triplets with $k_1k_2k_3>0$ generate two linearly independent eigenfunctions.
Figure~\ref{Fig8_Exmp7_Lambda3DConv} presents the $L^{\infty}$-errors between the first $10$ smallest exact and numerical eigenvalues 
on the semi-log scale against various $N$, in which we observe the expected exponential convergence. 
When $N$ increases to $16$, the $L^{\infty}$-errors between these exact and numerical eigenvalues have reached the level of $10^{-14}$.

\begin{figure}[!htb]
	\centering
	\subfigure[Relative error versus index] 
	{\includegraphics[scale=0.40]{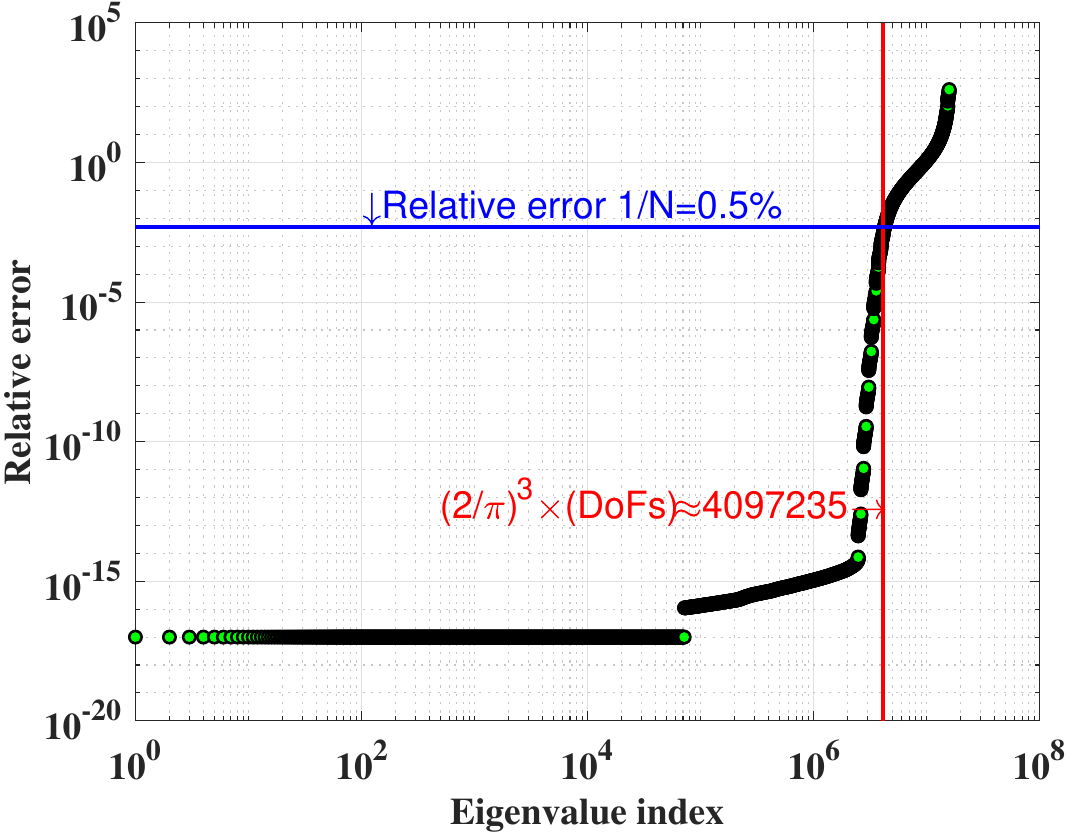}}
	\subfigure[Relative error region] 
	{\includegraphics[scale=0.36]{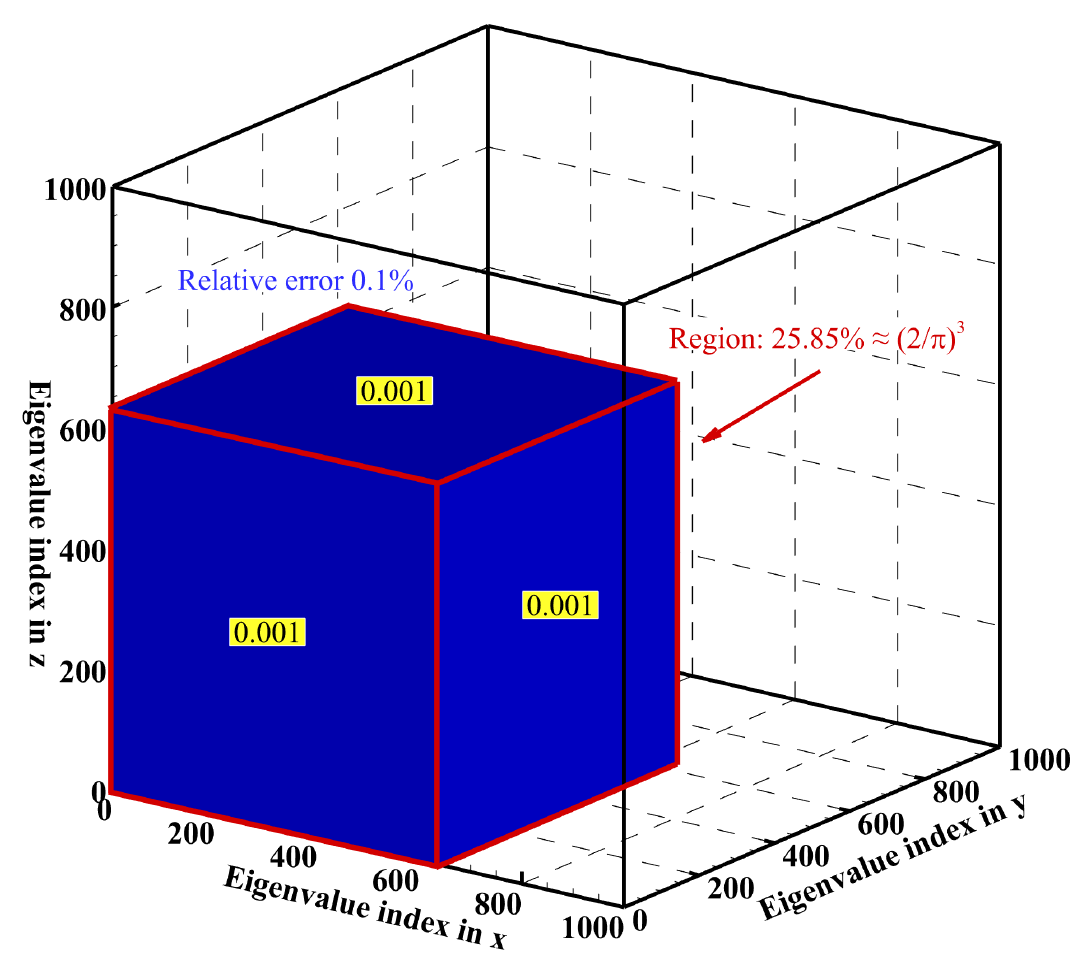}}
	\caption{Example \ref{Exmp5-6} in 3D:
		(a) relative error versus eigenvalue index for fixed $N=200$;
		(b) relative error region for fixed $N=1000$.
	}
	\label{Fig9_Exmp6_3DRelErrN200TrustRegN1000}
\end{figure}

Different from the 2D case, we find that the percentage of  ``trusted'' eigenvalues for 3D case by the proposed method is approximately $(2/{\pi})^3$. 
In Figure~\ref{Fig9_Exmp6_3DRelErrN200TrustRegN1000}(a), we depict the relative error of exact and numerical eigenvalues on the log-log scale versus the eigenvalue index with fixed $N=200$, 
which corresponds to a total DoFs$=15880001$. 
Once again, we draw the horizontal line $y=1/N=0.5\%$ and vertical line $x=(2/{\pi})^3\times({\rm DoFs})\approx{4097235}$, 
and find that they intersect at the error curve of the proposed method. 
Moreover, in Figure~\ref{Fig9_Exmp6_3DRelErrN200TrustRegN1000}(b), we draw $0.1\%=1/1000$ relative error region of our proposed method.
We observe that the percentage of ``trusted'' eigenvalues in 3D case increases to $25.85\%$($\approx(2/{\pi})^3$) for our spectral method.

In conclusion, these numerical results clearly demonstrate that our proposed method is suitable for accurate and large-scale computation for Maxwell's eigenvalue problems.

\section{Concluding Remarks}
\label{Sec6}
In this paper, we present novel spectral methods for double-curl source and eigenvalue problems of Maxwell equations in two and three dimensions that are capable of preserving the Gauss's law at the discrete level 
and exhibit a matrix-free, semi-analytic fast solution algorithm at the same time.  Several ingredients are indispensable for achieving these desirable properties: 
\begin{itemize}
\item an arbitrary order $\bs H({\rm curl})$-conforming vectorial spectral basis functions obtained by compact combination of Legendre polynomials;
\item a Gauss's law preserving approximation scheme based on a combination of the spectral basis with  Kikuchi's mixed formulation of the double-curl problems;  
\item a sophisticated reordering and decoupling procedure to decompose the original coupled algebraic system into decoupled sub-problems; 
\item analytic formulas of solutions in tensor-product form are obtained with the help of matrix diagonalisation of the 1D mass matrix. 
\end{itemize}

Future researches along this line include: 
(i) extend the current method to multi-domains to arrive at Gauss's law preserving spectral-element methods for solving Maxwell's double-curl and eigenvalue problems in more complex geometries; 
(ii) construct efficient matrix-free, scale-invariant iterative algorithms for spectral-element methods with the help of the proposed semi-analytic solution algorithm of single element and domain decomposition technique; 
(iii) develop rigorous error estimates of the proposed spectral-element scheme, especially the derivation of wave-number explicit error bounds for indefinite double-curl problems are of particular interest.  
We believe that this work, along with our subsequent studies, will shed light on the development of efficient and scalable algorithms for high-order vectorial spectral and spectral-element methods.

\vskip 10pt
\noindent\underline{\large\bf Acknowledgments}
\vskip 6pt
S. Lin and Z. Yang acknowledge the support from the National Natural Science Foundation of China (No. 12101399), the Shanghai Sailing Program (No. 21YF1421000), 
and the Strategic Priority Research Program of Chinese Academy of Sciences grant (No. XDA25010402) and the Key Laboratory of Scientific and Engineering Computing (Ministry of Education). 
H. Li acknowledges the support from the National Natural Science Foundation of China (Nos. 11871455 and 11971016).

\vspace{10pt} 
\noindent\textbf{Author Contributions}:
\begin{itemize}
\item 
S. Lin contributed to the software, formal analysis, writing -- original Draft;
\item 
H. Li contributed to the conceptualization, methodology, formal analysis,  writing -- review \& editing, funding acquisition;
\item 
Z. Yang contributed to the conceptualization, methodology, software, formal analysis, supervision, writing – original Draft, writing -- review \& editing, funding acquisition, project administration.
\end{itemize}
\vspace{10pt} 
\noindent\textbf{Data Availability}: Data will be made available on reasonable request. 

\vspace{10pt} 
\noindent\textbf{Competing Interests}: The authors have no relevant financial or non-financial interests to disclose.

\begin{appendix} 
\section{Derivation of Algorithm \texorpdfstring{\ref{Algo2}}{Lg}}
\renewcommand{\theequation}{\thesection.\arabic{equation}}
\label{Appx1-Algo2}
\begin{proof}
	To facilitate the design of the algorithm, we introduce the auxillary vectors as follows:
	\begin{equation}
	\label{eq:UshatVec3d}
	\begin{aligned}
	& \bs{\vec{\hat{U}}^{x}}  = ( \mathbf{Q}^{\intercal}\otimes\mathbf{Q}^{\intercal} ) \bs{\vec{U}^{x}}, \quad
	\bs{\vec{\hat{U}}^{y}}   = ( \mathbf{Q}^{\intercal}\otimes\mathbf{Q}^{\intercal} ) \bs{\vec{U}^{y}}, \quad
	\bs{\vec{\hat{U}}^{z}} = ( \mathbf{Q}^{\intercal}\otimes\mathbf{Q}^{\intercal} ) \bs{\vec{U}^{z}},
	\\
	& \bs{\vec{\hat{F}}^{x}} = ( \mathbf{Q}^{\intercal}\otimes\mathbf{Q}^{\intercal} ) \bs{\vec{F}^{x}}, \quad
	\bs{\vec{\hat{F}}^{y}} = ( \mathbf{Q}^{\intercal}\otimes\mathbf{Q}^{\intercal} ) \bs{\vec{F}^{y}}, \quad
	\bs{\vec{\hat{F}}^{z}} = ( \mathbf{Q}^{\intercal}\otimes\mathbf{Q}^{\intercal} ) \bs{\vec{F}^{z}},
	\\
	& 
	\bs{\vec{\hat{U}}^{s}} = ( \mathbf{Q}^{\intercal}\otimes\mathbf{Q}^{\intercal}\otimes\mathbf{Q}^{\intercal} ) \bs{\vec{U}^{s}}, \quad
	\bs{\vec{\hat{F}}^{s}} = ( \mathbf{Q}^{\intercal}\otimes\mathbf{Q}^{\intercal}\otimes\mathbf{Q}^{\intercal} ) \bs{\vec{F}^{s}}, \; \; s=1,2,3,
	\\
	& \bs{\vec{\hat{p}}} = ( \mathbf{Q}^{\intercal}\otimes\mathbf{Q}^{\intercal}\otimes\mathbf{Q}^{\intercal} ) \bs{\vec{p}}, \quad
	\bs{\vec{\hat{r}}} = ( \mathbf{Q}^{\intercal}\otimes\mathbf{Q}^{\intercal}\otimes\mathbf{Q}^{\intercal} ) \bs{\vec{r}}, 
	\end{aligned}
	\end{equation}
	and adopt the notations
	\begin{equation}
	\label{eq:UshatMat3d}
	\begin{aligned}
	& \mathbf{\hat{U}}^{x}={\rm ivec}(\bs{\vec{\hat{U}}^{x}}), \quad \mathbf{\hat{U}}^{y}={\rm ivec}(\bs{\vec{\hat{U}}^{y}}), \quad \mathbf{\hat{U}}^{z}={\rm ivec}(\bs{\vec{\hat{U}}^{z}}),
	\\
	& \mathbf{\hat{F}}^{x}={\rm ivec}(\bs{\vec{\hat{F}}^{x}}), \quad \mathbf{\hat{F}}^{y}={\rm ivec}(\bs{\vec{\hat{F}}^{y}}), \quad \mathbf{\hat{F}}^{z}={\rm ivec}(\bs{\vec{\hat{F}}^{z}}),
	\\
	& 
	\mathbf{\hat{U}}^{s}={\rm ivec}(\bs{\vec{\hat{U}}^{s}}), \quad \mathbf{\hat{F}}^{s}={\rm ivec}(\bs{\vec{\hat{F}}^{s}}),\;\; s=1,2,3, \quad 
	\mathbf{\hat{P}}={\rm ivec}(\bs{\vec{\hat{p}}}), \quad  \mathbf{\hat{R}}={\rm ivec}(\bs{\vec{\hat{r}}}).
	\end{aligned}
	\end{equation}

	Firstly, one can insert the fourth equation into the summation of the first three equations of system \eqref{eq:3dSub1MatEq} to obtain
	\begin{equation*}
	( \mathbf{M}\otimes\mathbf{M}\otimes\mathbf{I}_{N-1} + \mathbf{M}\otimes\mathbf{I}_{N-1}\otimes\mathbf{M} + \mathbf{I}_{N-1}\otimes\mathbf{M}\otimes\mathbf{M}) 
	\bs{\vec{p}} = \bs{\vec{F}^1} + \bs{\vec{F}^2} + \bs{\vec{F}^3} - {\kappa} \bs{\vec{r}}.
	\end{equation*}
	Then inserting auxiliary vectors in \eqref{eq:UshatVec3d} and matrices in \eqref{eq:UshatMat3d} into the above equation, one arrives at
	\begin{align*}
	& ( \mathbf{D}\otimes\mathbf{D}\otimes\mathbf{I}_{N-1} + \mathbf{D}\otimes\mathbf{I}_{N-1}\otimes\mathbf{D} + \mathbf{I}_{N-1}\otimes\mathbf{D}\otimes\mathbf{D}) 
	\bs{\vec{\hat{p}}} = \bs{\vec{\hat{F}}^1} + \bs{\vec{\hat{F}}^2} + \bs{\vec{\hat{F}}^3} - {\kappa} \bs{\vec{\hat{r}}}
	\\
	& \Longleftrightarrow \;\; 
	( d_{j}d_{k}+ d_{i}d_{k} + d_{i}d_{j} ) \hat{P}_{ijk}= \hat{F}_{ijk}^1 +  \hat{F}_{ijk}^2 +  \hat{F}_{ijk}^3 -  {\kappa}\hat{R}_{ijk}, 
	\quad \; 1\leq{i,j,k}\leq{N-1},
	\end{align*}
	which directly leads to semi-analytic formula for $\mathbf{P}$:
	\begin{equation*}
	\bs{\vec{p}} = (\mathbf{Q}\otimes\mathbf{Q}\otimes\mathbf{Q}) \bs{\vec{\hat{p}}}, \quad
	\bs{\vec{\hat{p}}}={\rm vec}(\mathbf{\hat{P}}), \quad \mathbf{\hat{P}}=(\hat{P}_{ijk}),  \quad
	\hat{P}_{ijk} = \frac{ \hat{F}_{ijk}^1 +  \hat{F}_{ijk}^2 +  \hat{F}_{ijk}^3 -  {\kappa}\hat{R}_{ijk} }{d_{j}d_{k}+ d_{i}d_{k} + d_{i}d_{j}}.
	\end{equation*}
	
	Next, using auxiliary vectors in \eqref{eq:UshatVec3d}, one can transform sub-system \eqref{eq:3dSub1MatEq} into 
	a linear system of $\bs{\vec{\hat{U}}^{1}}$-$\bs{\vec{\hat{U}}^{3}}$ and $\bs{\vec{\hat{P}}}$ as follows:
	\begin{equation*}
	\begin{aligned}
	& \mathbf{\hat{B}}_{11} \bs{\vec{\hat{U}}^{1}} + \mathbf{\hat{B}}_{12} \bs{\vec{\hat{U}}^{2}} + \mathbf{\hat{B}}_{13} \bs{\vec{\hat{U}}^{3}} + \mathbf{\hat{B}}_{14} \bs{\vec{\hat{P}}} = \bs{\vec{\hat{F}}^{1}}, 
	\\
	& \mathbf{\hat{B}}_{21} \bs{\vec{\hat{U}}^{1}} + \mathbf{\hat{B}}_{22}\bs{\vec{\hat{U}}^{2}} + \mathbf{\hat{B}}_{23} \bs{\vec{\hat{U}}^{3}} + \mathbf{\hat{B}}_{24} \bs{\vec{\hat{P}}}= \bs{\vec{\hat{F}}^{2}}, 
	\\
	& \mathbf{\hat{B}}_{31} \bs{\vec{\hat{U}}^{1}} + \mathbf{\hat{B}}_{32} \bs{\vec{\hat{U}}^{2}} + \mathbf{\hat{B}}_{33} \bs{\vec{\hat{U}}^{3}} + \mathbf{\hat{B}}_{34} \bs{\vec{\hat{P}}} = \bs{\vec{\hat{F}}^{3}}, 
	\\
	& \mathbf{\hat{B}}_{41} \bs{\vec{\hat{U}}^{1}} + \mathbf{\hat{B}}_{42} \bs{\vec{\hat{U}}^{2}} + \mathbf{\hat{B}}_{43}\bs{\vec{\hat{U}}^{3}} + \mathbf{\hat{B}}_{44} \bs{\vec{\hat{P}}} = \bs{\vec{\hat{R}}},
	\end{aligned}
	\end{equation*}
	where the components of the coefficient matrix are defined by
	\begin{align*}
	\mathbf{\hat{B}}_{11} &= \mathbf{I}_{N-1}\otimes\mathbf{D}\otimes\mathbf{I}_{N-1} + \mathbf{D}\otimes\mathbf{I}_{N-1}\otimes\mathbf{I}_{N-1} + {\kappa}\mathbf{D}\otimes\mathbf{D}\otimes\mathbf{I}_{N-1},
	\\
	\mathbf{\hat{B}}_{22} &= \mathbf{I}_{N-1}\otimes\mathbf{I}_{N-1}\otimes\mathbf{D} + \mathbf{D}\otimes\mathbf{I}_{N-1}\otimes\mathbf{I}_{N-1} + {\kappa}\mathbf{D}\otimes\mathbf{I}_{N-1}\otimes\mathbf{D},
	\\
	\mathbf{\hat{B}}_{33} &= \mathbf{I}_{N-1}\otimes\mathbf{I}_{N-1}\otimes\mathbf{D} + \mathbf{I}_{N-1}\otimes\mathbf{D}\otimes\mathbf{I}_{N-1} + {\kappa}\mathbf{I}_{N-1}\otimes\mathbf{D}\otimes\mathbf{D},
	\\
	\mathbf{\hat{B}}_{44} &= \mathbf{0}_{N-1}\otimes\mathbf{0}_{N-1}\otimes\mathbf{0}_{N-1},
	\\
	\mathbf{\hat{B}}_{12} &=\mathbf{\hat{B}}_{21}^{\intercal} = - \mathbf{D}\otimes\mathbf{I}_{N-1}\otimes\mathbf{I}_{N-1}, \quad 
	\mathbf{\hat{B}}_{13} = \mathbf{\hat{B}}_{31}^{\intercal} = - \mathbf{I}_{N-1}\otimes\mathbf{D}\otimes\mathbf{I}_{N-1}, 
	\\
	\mathbf{\hat{B}}_{23} &= \mathbf{\hat{B}}_{32}^{\intercal} = - \mathbf{I}_{N-1}\otimes\mathbf{I}_{N-1}\otimes\mathbf{D}, \quad
	\mathbf{\hat{B}}_{14} = \mathbf{\hat{B}}_{41}^{\intercal} = \mathbf{D}\otimes\mathbf{D}\otimes\mathbf{I}_{N-1},
	\\
	\mathbf{\hat{B}}_{24} &= \mathbf{\hat{B}}_{42}^{\intercal} = \mathbf{D}\otimes\mathbf{I}_{N-1}\otimes\mathbf{D}, \quad  \qquad \;
	\mathbf{\hat{B}}_{34} = \mathbf{\hat{B}}_{43}^{\intercal} = \mathbf{I}_{N-1}\otimes\mathbf{D}\otimes\mathbf{D}.
	\end{align*}
	For $1\leq{i,j,k}\leq{N-1}$, one can obtain a four-by-four linear system for each mode as follows:
	\begin{equation*}
	\begin{cases}
	( d_j + d_k+{\kappa}d_jd_k ) \hat{U}_{ijk}^{1}  -  d_k \hat{U}_{ijk}^{2} -  d_j \hat{U}_{ijk}^{3} +  d_jd_k \hat{P}_{ijk} = \hat{F}_{ijk}^1,
	\\
	- d_k \hat{U}_{ijk}^{1} + ( d_i+d_k+{\kappa}d_id_k ) \hat{U}_{ijk}^{2}  -  d_i \hat{U}_{ijk}^{3} +  d_id_k \hat{P}_{ijk} = \hat{F}_{ijk}^2,
	\\
	- d_j \hat{U}_{ijk}^{1}  -  d_i \hat{U}_{ijk}^{2} + ( d_i+d_j+{\kappa}d_id_j ) \hat{U}_{ijk}^{3}  +  d_id_j \hat{P}_{ijk} = \hat{F}_{ijk}^3,
	\\
	d_jd_k \hat{U}_{ijk}^{1}  +  d_id_k \hat{U}_{ijk}^{2} +  d_id_j \hat{U}_{ijk}^{3} = \hat{R}_{ijk},
	\end{cases}
	\end{equation*}
	which directly yields
	\begin{equation*}
	\begin{aligned}
	& 
	\bs{\vec{U}^{1}} = ( \mathbf{Q}\otimes\mathbf{Q}\otimes\mathbf{Q} ) \bs{\vec{\hat{U}}^{1}}, \quad
	\bs{\vec{\hat{U}}^{1}} = {\rm vec}(\mathbf{\hat{U}}^{1}), \quad \mathbf{\hat{U}}^{1} = (\hat{U}_{ijk}^{1}), \quad
	\hat{U}_{ijk}^{1} = \frac{\hat{R}_{ijk} + d_i \hat{F}_{ijk}^1 - d_i d_j d_k \hat{P}_{ijk}}{d_jd_k + d_id_k + d_id_j + {\kappa} d_id_jd_k},
	\\
	& 
	\bs{\vec{U}^{2}} = ( \mathbf{Q}\otimes\mathbf{Q}\otimes\mathbf{Q} ) \bs{\vec{\hat{U}}^{2}}, \quad
	\bs{\vec{\hat{U}}^{2}} = {\rm vec}(\mathbf{\hat{U}}^{2}), \quad \mathbf{\hat{U}}^{2} = (\hat{U}_{ijk}^{2}), \quad
	\hat{U}_{ijk}^{2} = \frac{\hat{R}_{ijk} + d_j \hat{F}_{ijk}^2 - d_i d_j d_k \hat{P}_{ijk}}{d_jd_k + d_id_k + d_id_j + {\kappa} d_id_jd_k},
	\\
	&
	\bs{\vec{U}^{3}} = ( \mathbf{Q}\otimes\mathbf{Q}\otimes\mathbf{Q} ) \bs{\vec{\hat{U}}^{3}}, \quad
	\bs{\vec{\hat{U}}^{3}} = {\rm vec}(\mathbf{\hat{U}}^{3}), \quad \mathbf{\hat{U}}^{3} = (\hat{U}_{ijk}^{3}), \quad
	\hat{U}_{ijk}^{3} = \frac{\hat{R}_{ijk} + d_k \hat{F}_{ijk}^3 - d_i d_j d_k \hat{P}_{ijk}}{d_jd_k + d_id_k + d_id_j + {\kappa} d_id_jd_k}.
	\end{aligned}
	\end{equation*}
	Finally, one also can transform sub-system \eqref{eq:3dSub2MatEq} into a linear and uncoupled system of $\bs{\vec{\hat{U}}^{x}}$-$\bs{\vec{\hat{U}}^{z}}$ as follows:
	\begin{equation*}
	\begin{aligned}
	&
	( \mathbf{I}_{N-1}\otimes\mathbf{D} + \mathbf{D}\otimes\mathbf{I}_{N-1} + {\kappa}\mathbf{D}\otimes\mathbf{D} ) \bs{\vec{\hat{U}}^{x}} = \bs{\vec{\hat{F}}^{x}},
	\\
	&
	( \mathbf{I}_{N-1}\otimes\mathbf{D} + \mathbf{D}\otimes\mathbf{I}_{N-1} + {\kappa}\mathbf{D}\otimes\mathbf{D} ) \bs{\vec{\hat{U}}^{y}} = \bs{\vec{\hat{F}}^{y}},
	\\
	&
	( \mathbf{I}_{N-1}\otimes\mathbf{D} + \mathbf{D}\otimes\mathbf{I}_{N-1} + {\kappa}\mathbf{D}\otimes\mathbf{D} ) \bs{\vec{\hat{U}}^{z}} = \bs{\vec{\hat{F}}^{z}}.
	\end{aligned}
	\end{equation*}
	For $1\leq{i,j,k}\leq{N-1}$, one can obtain three independent equations for each mode as follows:
	\begin{align*}
	&
	 (d_j + d_k + \kappa {d_jd_k}) \hat{U}_{jk}^{x} = \hat{F}_{jk}^{x}, \; \quad 
	(d_i + d_k + \kappa {d_id_k}) \hat{U}_{ik}^{y} = \hat{F}_{ik}^{y}, \; \quad
	(d_i + d_j + \kappa {d_id_j}) \hat{U}_{ij}^{z} = \hat{F}_{ij}^{z},
	\end{align*}
	which directly give
	\begin{equation*}
	\begin{aligned}
	& 
	\bs{\vec{U}^{x}} = {\rm vec}(\mathbf{U}^{x}), \quad 
	\mathbf{U}^{x} = \mathbf{Q}\mathbf{\hat{U}}^{x}\mathbf{Q}^{\intercal}, \quad  
	\mathbf{\hat{U}}^{x}=(\hat{U}_{jk}^{x}), \quad
	\hat{U}_{jk}^{x} = \frac{\hat{F}_{jk}^{x}}{ d_{j} + d_{k} + {\kappa}d_{j}d_{k}}, \quad
	\\
	& 
	\bs{\vec{U}^{y}} = {\rm vec}(\mathbf{U}^{y}), \quad 
	\mathbf{U}^{y} = \mathbf{Q}\mathbf{\hat{U}}^{y}\mathbf{Q}^{\intercal}, \quad  
	\mathbf{\hat{U}}^{y}=(\hat{U}_{ik}^{y}), \quad
	\hat{U}_{ik}^{y} = \frac{\hat{F}_{ik}^{y}}{ d_{i} + d_{k} + {\kappa}d_{i}d_{k}}, \quad
	\\
	&
	\bs{\vec{U}^{z}} = {\rm vec}(\mathbf{U}^{z}), \quad 
	\mathbf{U}^{z} = \mathbf{Q}\mathbf{\hat{U}}^{z}\mathbf{Q}^{\intercal}, \quad 
	\mathbf{\hat{U}}^{z}=(\hat{U}_{ij}^{z}), \quad
	\hat{U}_{ij}^{z} = \frac{\hat{F}_{ij}^{z}}{ d_{i} + d_{j} + {\kappa}d_{i}d_{j}}.
	\end{aligned}
	\end{equation*}
\end{proof}

\section{Proof of Theorem \texorpdfstring{\ref{Thm2-LamU3dSol}}{Lg}}
\renewcommand{\theequation}{\thesection.\arabic{equation}}
\label{Appx2-Thm2}
\begin{proof}
We choose the first 6-by-3 block submatrix of the permutation matrix given in \eqref{eq:MaxwPemMat3d},
adopt the same notations $\mathbf{U}^1$-$\mathbf{U}^3$, $\bs{\vec{U}^{x}}$, $\bs{\vec{U}^{y}}$ and $\bs{\vec{U}^{z}}$ in \eqref{eq:UshatMat3d}, 
and further obtain the equivalent system of \eqref{eq:EigAlgSys3d} in the matrix-vector product form as follows:
\begin{equation}
\label{eq:EigMatVecEq3d}
\begin{aligned}
& \mathbf{{B}}_{11} \bs{\vec{U}^{1}} + \mathbf{{B}}_{12} \bs{\vec{U}^{2}} + \mathbf{{B}}_{13} \bs{\vec{U}^{3}} 
= {\lambda}^{N}( \mathbf{M}\otimes\mathbf{M}\otimes\mathbf{I}_{N-1} ) \bs{\vec{U}^{1}}, 
\\
&
\mathbf{{B}}_{21} \bs{\vec{U}^{1}} + \mathbf{{B}}_{22} \bs{\vec{U}^{2}} + \mathbf{{B}}_{23} \bs{\vec{U}^{3}} 
= {\lambda}^{N}( \mathbf{M}\otimes\mathbf{I}_{N-1}\otimes\mathbf{M} ) \bs{\vec{U}^{2}}, 
\\
& \mathbf{{B}}_{31} \bs{\vec{U}^{1}} + \mathbf{{B}}_{32} \bs{\vec{U}^{2}} + \mathbf{{B}}_{33} \bs{\vec{U}^{3}} 
= {\lambda}^{N}( \mathbf{I}_{N-1}\otimes\mathbf{M}\otimes\mathbf{M} ) \bs{\vec{U}^{3}}, 
\\
& \mathbf{{B}}_{xx} \bs{\vec{U}^{x}} = {\lambda}^{N}( \mathbf{M}\otimes\mathbf{M} ) \bs{\vec{U}^{x}}, \quad
\mathbf{{B}}_{yy} \bs{\vec{U}^{y}} = {\lambda}^{N}( \mathbf{M}\otimes\mathbf{M} ) \bs{\vec{U}^{y}}, \quad
\mathbf{{B}}_{zz}\bs{\vec{U}^{z}} = {\lambda}^{N}( \mathbf{M}\otimes\mathbf{M} ) \bs{\vec{U}^{z}},
\end{aligned}
\end{equation}
where the components of the coefficient matrix are given by 
\begin{align*}
& \mathbf{{B}}_{11} = \mathbf{I}_{N-1}\otimes\mathbf{M}\otimes\mathbf{I}_{N-1} + \mathbf{M}\otimes\mathbf{I}_{N-1}\otimes\mathbf{I}_{N-1}, \quad
\mathbf{{B}}_{22} = \mathbf{I}_{N-1}\otimes\mathbf{I}_{N-1}\otimes\mathbf{M} + \mathbf{M}\otimes\mathbf{I}_{N-1}\otimes\mathbf{I}_{N-1},
\\
& \mathbf{{B}}_{33} = \mathbf{I}_{N-1}\otimes\mathbf{I}_{N-1}\otimes\mathbf{M} + \mathbf{I}_{N-1}\otimes\mathbf{M}\otimes\mathbf{I}_{N-1}, \quad
\mathbf{{B}}_{12} =\mathbf{{B}}_{21}^{\intercal} = - \mathbf{M}\otimes\mathbf{I}_{N-1}\otimes\mathbf{I}_{N-1}, 
\\
&
\mathbf{{B}}_{13} = \mathbf{{B}}_{31}^{\intercal} = - \mathbf{I}_{N-1}\otimes\mathbf{M}\otimes\mathbf{I}_{N-1},  \quad \qquad \qquad \quad  \;\;
\mathbf{{B}}_{23} = \mathbf{{B}}_{32}^{\intercal} = - \mathbf{I}_{N-1}\otimes\mathbf{I}_{N-1}\otimes\mathbf{M},
\\
& \mathbf{{B}}_{xx} = \mathbf{{B}}_{yy} = \mathbf{{B}}_{zz} = \mathbf{I}_{N-1}\otimes\mathbf{M} + \mathbf{M}\otimes\mathbf{I}_{N-1}.
\end{align*}

Introduce auxillary vectors
\begin{equation*}
\bs{\vec{\hat{U}}^{s}}  = ( \mathbf{Q}^{\intercal}\otimes\mathbf{Q}^{\intercal}\otimes\mathbf{Q}^{\intercal}) \bs{\vec{U}^{s}},
\quad \; s=1,2,3.
\end{equation*}
Then the first three equations of \eqref{eq:EigMatVecEq3d} can be transformed into a linear system of $\bs{\vec{\hat{U}}^1}$-$\bs{\vec{\hat{U}}^3}$ as follows:
\begin{equation*}
\begin{aligned}
& \mathbf{\hat{B}}_{11} \bs{\vec{\hat{U}}^1} + \mathbf{\hat{B}}_{12} \bs{\vec{\hat{U}}^2} + \mathbf{\hat{B}}_{13} \bs{\vec{\hat{U}}^3} 
= {\lambda}^{N} ( \mathbf{D}\otimes\mathbf{D}\otimes\mathbf{I}_{N-1} ) \bs{\vec{\hat{U}}^1},
\\
& \mathbf{\hat{B}}_{21} \bs{\vec{\hat{U}}^1} + \mathbf{\hat{B}}_{22} \bs{\vec{\hat{U}}^2} + \mathbf{\hat{B}}_{23} \bs{\vec{\hat{U}}^3} 
= {\lambda}^{N} ( \mathbf{D}\otimes\mathbf{I}_{N-1}\otimes\mathbf{D}) \bs{\vec{\hat{U}}^2},
\\
& \mathbf{\hat{B}}_{31} \bs{\vec{\hat{U}}^1} + \mathbf{\hat{B}}_{32} \bs{\vec{\hat{U}}^2} + \mathbf{\hat{B}}_{33} \bs{\vec{\hat{U}}^3} 
= {\lambda}^{N} ( \mathbf{I}_{N-1}\otimes\mathbf{D}\otimes\mathbf{D}) \bs{\vec{\hat{U}}^3},
\end{aligned}
\end{equation*}
where the components of the coefficient matrix are given by 
\begin{align*}
& \mathbf{\hat{B}}_{11} =  \mathbf{I}_{N-1}\otimes\mathbf{D}\otimes\mathbf{I}_{N-1} + \mathbf{D}\otimes\mathbf{I}_{N-1}\otimes\mathbf{I}_{N-1}, \quad 
\mathbf{\hat{B}}_{22} =  \mathbf{I}_{N-1}\otimes\mathbf{I}_{N-1}\otimes\mathbf{D} + \mathbf{D}\otimes\mathbf{I}_{N-1}\otimes\mathbf{I}_{N-1}, 
\\
& 
\mathbf{\hat{B}}_{33} = \mathbf{I}_{N-1}\otimes\mathbf{I}_{N-1}\otimes\mathbf{D} + \mathbf{I}_{N-1}\otimes\mathbf{D}\otimes\mathbf{I}_{N-1}, \quad 
\mathbf{\hat{B}}_{12} = \mathbf{\hat{B}}_{21} = - \mathbf{D}\otimes\mathbf{I}_{N-1}\otimes\mathbf{I}_{N-1}, 
\\
& \mathbf{\hat{B}}_{13} = \mathbf{\hat{B}}_{31} = - \mathbf{I}_{N-1}\otimes\mathbf{D}\otimes\mathbf{I}_{N-1}, \quad \qquad \qquad \quad \;
\mathbf{\hat{B}}_{23} = \mathbf{\hat{B}}_{32} = - \mathbf{I}_{N-1}\otimes\mathbf{I}_{N-1}\otimes\mathbf{D}.
\end{align*}
One can easily obtain for each mode a four-by-four linear system as follows:
\begin{equation}
\label{eq:EigModeEq3dIn}
\begin{cases}
(d_j+d_k)\hat{U}_{ijk}^{1} - d_k\hat{U}_{ijk}^{2} -  d_j\hat{U}_{ijk}^{3} = {\lambda}_{ijk}^{N}   d_jd_k\hat{U}_{ijk}^{1},
\\
- d_k\hat{U}_{ijk}^{1} + (d_i+d_k)\hat{U}_{ijk}^{2} -  d_i\hat{U}_{ijk}^{3} = {\lambda}_{ijk}^{N}  d_id_k\hat{U}_{ijk}^{2},
\\
- d_j\hat{U}_{ijk}^{1}  -  d_i\hat{U}_{ijk}^{2} +  (d_i+d_j)\hat{U}_{ijk}^{3} = {\lambda}_{ijk}^{N}  d_id_j\hat{U}_{ijk}^{3}.
\end{cases}
\end{equation}
By summing up the above three equations of \eqref{eq:EigModeEq3dIn}, one has 
\begin{align*}
{\lambda}_{ijk}^{N} ( d_jd_k\hat{U}_{ijk}^{1} + d_id_k\hat{U}_{ijk}^{2} + d_id_j\hat{U}_{ijk}^{3} ) = 0
\;\; \Longleftrightarrow  \;\; 
{\lambda}_{ijk}^{N} = 0, \;\; \text{or} \;\;
d_jd_k\hat{U}_{ijk}^{1} + d_id_k\hat{U}_{ijk}^{2} + d_id_j\hat{U}_{ijk}^{3} = 0.
\end{align*}
When ${\lambda}_{ijk}^{N}=0$, one easily obtains
\begin{align*}
\hat{U}_{ijk}^{1} = \hat{U}_{ijk}^{2} = \hat{U}_{ijk}^{3} = 1. 
\end{align*}
When ${\lambda}_{ijk}^{N}\neq0$, it yields 
\begin{align*}
d_jd_k\hat{U}_{ijk}^{1} + d_id_k\hat{U}_{ijk}^{2} + d_id_j\hat{U}_{ijk}^{3} = 0.
\end{align*}
One further has
\begin{align*}
& (d_jd_k + d_id_k + d_id_j) \hat{U}_{ijk}^{1} = {\lambda}_{ijk}^{N}   d_id_jd_k\hat{U}_{ijk}^{1},
\\
& (d_jd_k + d_id_k + d_id_j) \hat{U}_{ijk}^{2} = {\lambda}_{ijk}^{N}   d_id_jd_k\hat{U}_{ijk}^{2},
\\
& (d_jd_k + d_id_k + d_id_j) \hat{U}_{ijk}^{3} = {\lambda}_{ijk}^{N}   d_id_jd_k\hat{U}_{ijk}^{3},
\end{align*}
which can be solved directly by  
\begin{align*}
{\lambda}_{ijk}^{N} = \frac{1}{d_i} + \frac{1}{d_j} + \frac{1}{d_k}. 
\end{align*}
Due to the coefficient matrix of \eqref{eq:EigModeEq3dIn} is rank deficient, we can obtain two linear independent solutions of \eqref{eq:EigModeEq3dIn} as follows:
\begin{align*}
& \hat{U}_{ijk}^{1}  = 1, \quad \hat{U}_{ijk}^{2}  = 0, \quad
\hat{U}_{ijk}^{3}  = -d_k/d_i,
\;\quad  \text{and} \; \quad
\hat{U}_{ijk}^{1}  = 0, \quad  \hat{U}_{ijk}^{2}  = 1, \quad
\hat{U}_{ijk}^{3}  = - d_k/d_j.
\end{align*}

Finally, define 
\begin{equation*}
\bs{\vec{\hat{U}}^x}  = ( \mathbf{Q}^{\intercal}\otimes\mathbf{Q}^{\intercal} ) \bs{\vec{U}^x}, \quad
\bs{\vec{\hat{U}}^y}  = ( \mathbf{Q}^{\intercal}\otimes\mathbf{Q}^{\intercal} ) \bs{\vec{U}^y}, \quad
\bs{\vec{\hat{U}}^z}  = ( \mathbf{Q}^{\intercal}\otimes\mathbf{Q}^{\intercal} ) \bs{\vec{U}^z}.
\end{equation*}
The last three equations of \eqref{eq:EigMatVecEq3d} can also be transformed into linear and decoulped equations of 
$\bs{\vec{\hat{U}}^{x}}$, $\bs{\vec{\hat{U}}^{y}}$ and $\bs{\vec{\hat{U}}^{z}}$ as follows:
\begin{align*}
& ( \mathbf{I}_{N-1}\otimes\mathbf{D} + \mathbf{D}\otimes\mathbf{I}_{N-1} ) \bs{\vec{\hat{U}}^{x}} 
= {\lambda}^{N} ( \mathbf{D}\otimes\mathbf{D}) \bs{\vec{\hat{U}}^{x}}, \;
\\
& ( \mathbf{I}_{N-1}\otimes\mathbf{D} + \mathbf{D}\otimes\mathbf{I}_{N-1 }) \bs{\vec{\hat{U}}^{y}} 
= {\lambda}^{N} ( \mathbf{D}\otimes\mathbf{D}) \bs{\vec{\hat{U}}^{y}}, \;
\\
& ( \mathbf{I}_{N-1}\otimes\mathbf{D} + \mathbf{D}\otimes\mathbf{I}_{N-1} ) \bs{\vec{\hat{U}}^{z}} 
= {\lambda}^{N} ( \mathbf{D}\otimes\mathbf{D} ) \bs{\vec{\hat{U}}^{z}}.
\end{align*}
For $1\leq{i,j}\leq{N-1}$, we have 
\begin{align*}
& (d_i + d_j) \hat{U}_{ij}^{x} = {\lambda}_{ij}^{N} d_id_j \hat{U}_{ij}^{x}, \quad
(d_i + d_j) \hat{U}_{ij}^{y} = {\lambda}_{ij}^{N} d_id_j \hat{U}_{ij}^{y}, \quad
(d_i + d_j) \hat{U}_{ij}^{z} = {\lambda}_{ij}^{N} d_id_j \hat{U}_{ij}^{z}, 
\end{align*}
which can solved directly by 
\begin{align*}
& {\lambda}_{ij}^{N} = \frac{1}{d_i} + \frac{1}{d_j},  
\quad
\hat{U}_{ij}^{x} = 1, \quad \hat{U}_{ij}^{y} = 1, \quad \hat{U}_{ij}^{z} = 1.
\end{align*}
\end{proof}
\end{appendix}


\bibliographystyle{plain}
\bibliography{refpapers}

\begin{thebibliography}{10}

\bibitem{Adam1975}
R.~A. Adams.
\newblock {\em Sobolev Spaces}.
\newblock Acadmic Press, New York, 1975.

\bibitem{Balsara-Kim-2004}
D.~S. Balsara and J.~Kim.
\newblock An intercomparison between divergence-cleaning and staggered mesh
  formulations for numerical magnetohydrodynamics.
\newblock {\em Astrophys J.}, 602:1079--1090, 2004.

\bibitem{BoffiBrezziFortin2013}
D.~Boffi, F.~Brezzi, and M.~Fortin.
\newblock {\em Mixed Finite Element Methods and Applications}, volume~44 of
  {\em Springer Series in Computational Mathematics}.
\newblock Springer, Heidelberg, 2013.

\bibitem{Bramble-Kolev-Pasciak-2005}
J.~H. Bramble, T.~V. Kolev, and J.~E. Pasciak.
\newblock The approximation of the {M}axwell eigenvalue problem using a
  least-squares method.
\newblock {\em Math. Comp.}, 74(252):1575--1598, 2005.

\bibitem{Brenner-Li-Sun-2008}
S.~C. Brenner, F.~Li, and L.-Y. Sung.
\newblock A locally divergence-free interior penalty method for two-dimensional
  curl-curl problems.
\newblock {\em SIAM J. Numer. Anal.}, 46(3):1190--1211, 2008.

\bibitem{Buchholz-1986}
D.~Buchholz.
\newblock Gauss' law and the infraparticle problem.
\newblock {\em Phys. Lett. B}, 174(3):331--334, 1986.

\bibitem{ChenFF1984}
F.~Chen.
\newblock {\em Introduction to Plasma Physics and Controlled Fusion}, volume~1.
\newblock Springer, 1984.

\bibitem{Ciarlet-Wu-Zou-2014}
P.~Ciarlet, H.~Wu, and J.~Zou.
\newblock Edge element methods for {M}axwell's equations with strong
  convergence for {G}auss' laws.
\newblock {\em SIAM J. Numer. Anal.}, 52(2):779--807, 2014.

\bibitem{Ciarlet1988}
P.~G. Ciarlet.
\newblock {\em Mathematical elasticity. {V}ol. {I}. {T}hree-dimensional
  elasticity}, volume~20 of {\em Studies in Mathematics and its Applications}.
\newblock North-Holland Publishing Co., Amsterdam, 1988.

\bibitem{Davidson2001}
P.~A. Davidson.
\newblock {\em An Introduction to Magnetohydrodynamics}.
\newblock Cambridge Texts in Applied Mathematics. Cambridge University Press,
  Cambridge, 2001.

\bibitem{Dedner-Kemm-Kroner-2002}
A.~Dedner, F.~Kemm, D.~Kr\"{o}ner, C.-D. Munz, T.~Schnitzer, and M.~Wesenberg.
\newblock Hyperbolic divergence cleaning for the {MHD} equations.
\newblock {\em J. Comput. Phys.}, 175(2):645--673, 2002.

\bibitem{Duan-Ma-Zou-2021}
H.~Duan, J.~Ma, and J.~Zou.
\newblock Mixed finite element method with {G}auss's law enforced for the
  {M}axwell eigenproblem.
\newblock {\em SIAM J. Sci. Comput.}, 43(6):A3677--A3712, 2021.

\bibitem{Fernandes-Raffetto-2002}
P.~Fernandes and M.~Raffetto.
\newblock Counterexamples to the currently accepted explanation for spurious
  modes and necessary and sufficient conditions to avoid them.
\newblock {\em IEEE transactions on magnetics}, 38(2):653--656, 2002.

\bibitem{Fuentes-Keith-2015}
F.~Fuentes, B.~Keith, L.~Demkowicz, and S.~Nagaraj.
\newblock Orientation embedded high order shape functions for the exact
  sequence elements of all shapes.
\newblock {\em Comput. Math. Appl.}, 70(4):353--458, 2015.

\bibitem{GolubCharles2013}
G.~H. Golub and C.~F. Van~Loan.
\newblock {\em Matrix computations}.
\newblock Johns Hopkins Studies in the Mathematical Sciences. Johns Hopkins
  University Press, Baltimore, MD, fourth edition, 2013.

\bibitem{Gordon-Hall-1973}
W.~J. Gordon and C.~A. Hall.
\newblock Transfinite element methods: blending-function interpolation over
  arbitrary curved element domains.
\newblock {\em Numer. Math.}, 21:109--129, 1973.

\bibitem{Hiptmair-Xu-2007}
R.~Hiptmair and J.~Xu.
\newblock Nodal auxiliary space preconditioning in {${\bf H}({\bf curl})$} and
  {${\bf H}({\rm div})$} spaces.
\newblock {\em SIAM J. Numer. Anal.}, 45(6):2483--2509, 2007.

\bibitem{Hu-Zou-2004}
Q.~Hu and J.~Zou.
\newblock Substructuring preconditioners for saddle-point problems arising from
  {M}axwell's equations in three dimensions.
\newblock {\em Math. Comp.}, 73(245):35--61, 2004.

\bibitem{JoanJohnsonWinn2008}
J.~D. Joannopoulos, S.~G. Johnson, J.~N. Winn, and R.~D. Meade.
\newblock {\em Photonic Crystals: Molding the Flow of Light}.
\newblock Princeton University Press, 2008.

\bibitem{Kijowski-Rudolph-2002}
J.~Kijowski and G.~Rudolph.
\newblock On the {G}auss law and global charge for quantum chromodynamics.
\newblock {\em J. Math. Phys.}, 43(4):1796--1808, 2002.

\bibitem{Kikuchi-1987}
F.~Kikuchi.
\newblock Mixed and penalty formulations for finite element analysis of an
  eigenvalue problem in electromagnetism.
\newblock {\em Comput. Methods Appl. Mech. Engrg.}, 64(1-3):509--521, 1987.

\bibitem{Liang-Xu-2022}
Q.~Liang and X.~Xu.
\newblock A two-level preconditioned {H}elmholtz-{J}acobi-{D}avidson method for
  the {M}axwell eigenvalue problem.
\newblock {\em Math. Comp.}, 91(334):623--657, 2022.

\bibitem{Liang-Xu-2023}
Q.~Liang and X.~Xu.
\newblock A two-level preconditioned {H}elmholtz subspace iterative method for
  {M}axwell eigenvalue problems.
\newblock {\em SIAM J. Numer. Anal.}, 61(2):642--674, 2023.

\bibitem{Liu-Tobon-Tang-Liu-2015}
N.~Liu, L.~Tob\'{o}n, Y.~Tang, and Q.~Liu.
\newblock Mixed spectral element method for 2{D} {M}axwell's eigenvalue
  problem.
\newblock {\em Commun. Comput. Phys.}, 17(2):458--486, 2015.

\bibitem{MarkowichRinghofer2012}
P.~A. Markowich, C.~A. Ringhofer, and C.~Schmeiser.
\newblock {\em Semiconductor Equations}.
\newblock Springer Science \& Business Media, 2012.

\bibitem{Monk2003}
P.~Monk.
\newblock {\em Finite {E}lement {M}ethods for {M}axwell's {E}quations}.
\newblock Numerical Mathematics and Scientific Computation. Oxford University
  Press, New York, 2003.

\bibitem{Mund-Rehren-Schroer-2020}
J.~Mund, K-H. Rehren, and B.~Schroer.
\newblock Gauss' law and string-localized quantum field theory.
\newblock {\em J. High Energy Phys.}, no.(1):001, 26pp, 2020.

\bibitem{Munz-Omnes-Schneider-2000}
C.-D. Munz, P.~Omnes, R.~Schneider, E.~Sonnendr\"{u}cker, and U.~Vo\ss.
\newblock Divergence correction techniques for {M}axwell solvers based on a
  hyperbolic model.
\newblock {\em J. Comput. Phys.}, 161(2):484--511, 2000.

\bibitem{Nedelec-1980}
J.-C. N\'{e}d\'{e}lec.
\newblock Mixed finite elements in {${\bf R}^{3}$}.
\newblock {\em Numer. Math.}, 35(3):315--341, 1980.

\bibitem{Nedelec-1986}
J.-C. N\'{e}d\'{e}lec.
\newblock A new family of mixed finite elements in {${\bf R}^3$}.
\newblock {\em Numer. Math.}, 50(1):57--81, 1986.

\bibitem{Orfanidis2002}
S.~J. Orfanidis.
\newblock {\em Electromagnetic Waves and Antennas}.
\newblock Rutgers University New Brunswick, NJ, 2002.

\bibitem{Pazner-Kolev-Dohrmann-2023}
W.~Pazner, T.~Kolev, and C.~R. Dohrmann.
\newblock Low-order preconditioning for the high-order finite element de {R}ham
  complex.
\newblock {\em SIAM J. Sci. Comput.}, 45(2):A675--A702, 2023.

\bibitem{Sanders-Dappiaggi-Hack-2014}
K.~Sanders, C.~Dappiaggi, and T.-P. Hack.
\newblock Electromagnetism, local covariance, the {A}haronov-{B}ohm effect and
  {G}auss' law.
\newblock {\em Comm. Math. Phys.}, 328(2):625--667, 2014.

\bibitem{Shen-1994b}
J.~Shen.
\newblock Efficient spectral-{G}alerkin method {I}. {D}irect solvers for
  second- and fourth-order equations by using {L}egendre polynomials.
\newblock {\em SIAM J. Sci. Comput.}, 15:1489--1505, 1994.

\bibitem{Strassen-1969}
V.~Strassen.
\newblock Gaussian elimination is not optimal.
\newblock {\em Numer. Math.}, 13:354--356, 1969.

\bibitem{SunZhou2017}
J.~Sun and A.~Zhou.
\newblock {\em Finite Element Methods for Eigenvalue Problems}.
\newblock Monographs and Research Notes in Mathematics. CRC Press, Boca Raton,
  FL, 2017.

\bibitem{Taflove-1980}
A.~Taflove.
\newblock Application of the finite-difference time-domain method to sinusoidal
  steady-state electromagnetic-penetration problems.
\newblock {\em IEEE Trans. Elect. Compat.}, EMC-22(3):191--202, 1980.

\bibitem{Yee-1966}
K.~Yee.
\newblock Numerical solution of initial boundary value problems involving
  {M}axwell's equations in isotropic media.
\newblock {\em IEEE Trans. Antennas and Prop.}, 14(3):302--307, 1966.

\bibitem{Yousept-Zou-2017}
I.~Yousept and J.~Zou.
\newblock Edge element method for optimal control of stationary {M}axwell
  system with {G}auss law.
\newblock {\em SIAM J. Numer. Anal.}, 55(6):2787--2810, 2017.

\bibitem{Zaglmayr2006}
S.~Zaglmayr.
\newblock {\em High Order Finite Element Methods for Electromagnetic Field
  Computation (Ph.D. thesis)}.
\newblock KEPLER, JOHANNES, 2006.

\bibitem{Zhang-Wang-Yang-2020}
J.~Zhang, W.~M. Wang, X.~H. Yang, and et~al.
\newblock Double-cone ignition scheme for inertial confinement fusion.
\newblock {\em Phil. Trans. Royal Soc. A}, 378(2184):20200015, 2020.

\bibitem{Zhang-2015}
Z.~Zhang.
\newblock How many numerical eigenvalues can we trust?
\newblock {\em J. Sci. Comput.}, 65(2):455--466, 2015.

\end{thebibliography}

\end{document}